\documentclass{article}

\usepackage{arxiv}
\usepackage{comment}
\usepackage{amsthm}
\usepackage[utf8]{inputenc} 
\usepackage[T1]{fontenc}    
\usepackage{hyperref}       
\usepackage{url}            
\usepackage{booktabs}       
\usepackage{amsfonts}       
\usepackage{amssymb}
\usepackage{nicefrac}       
\usepackage{microtype}      
\usepackage{lipsum}
\usepackage{graphicx}
\graphicspath{ {./images/} }
\usepackage{amsmath}
\usepackage[numbers,sort&compress]{natbib}
\usepackage{enumitem}
\usepackage{xcolor}
\usepackage[nameinlink,capitalise]{cleveref}
\crefname{assumption}{assumption}{assumptions}
\crefname{appendix}{Appendix}{Appendices}
\Crefname{appendix}{Appendix}{Appendices}
\usepackage{subfigure}
\DeclareMathOperator*{\argmin}{argmin} 
\usepackage{algorithm}
\usepackage{algorithmic}

\usepackage{cases}
\usepackage{multirow}

\newtheorem{proposition}{Proposition}
\newtheorem{theorem}{Theorem}
\newtheorem{assumption}{Assumption}

\newtheorem{remark}{Remark}

\usepackage{tikz}

\title{Invariant Manifolds of Discrete-time Dynamical Systems with Nonlinear Exosystems via Hybrid Physics-Informed Neural Networks}

\author{
\textbf{Dimitrios G. Patsatzis\textcolor{blue}{$^{1}$}, Nikolaos Kazantzis\textcolor{blue}{$^{2}$}, Ioannis G. Kevrekidis\textcolor{blue}{$^{3, 4, 5}$}, Lucia Russo\textcolor{blue}{$^{6}$}, Constantinos Siettos\textcolor{blue}{$^{7,}$}\thanks{Corresponding author, email: \texttt{constantinos.siettos@unina.it}}}
{}\\
\textcolor{blue}{$^{(1)}$}Modelling Engineering Risk \& Complexity, \emph{Scuola Superiore Meridionale}, Naples, Italy \\
\textcolor{blue}{$^{(2)}$}Dept. of Chemical Engineering, \emph{Worcester Polytechnic Institute}, Worcester, USA \\
\textcolor{blue}{$^{(3)}$}Dept. of Chemical \& Biomolecular Engineering, \emph{Johns Hopkins University}, Baltimore, USA \\
\textcolor{blue}{$^{(4)}$}Dept. of Applied Mathematics and Statistics, \emph{Johns Hopkins University}, Baltimore, USA\\
\textcolor{blue}{$^{(5)}$}Medical School, Department of Urology, \emph{Johns Hopkins University}, Baltimore, USA\\
\textcolor{blue}{$^{(6)}$}Institute of Science and Technology for Energy and Sustainable Mobility, \\
\emph{Consiglio Nazionale delle Ricerche}, Naples, Italy\\
\textcolor{blue}{$^{(7)}$}Dipartimento di Matematica e Applicazioni ``Renato Caccioppoli'', \\ \emph{University of Naples Federico II}, Naples, Italy\\
}

\begin{document}
\maketitle
\begin{abstract}
We propose a hybrid physics-informed machine learning framework to approximate invariant manifolds (IMs) of discrete-time dynamical systems driven by exogenous autonomous dynamics (exosystems). Such hierarchical/coupled systems arise in applications ranging from control theory to modelling of collective behaviour of multi-agent systems (e.g., bird flocks, fish schools, crowd and traffic dynamics) under hierarchical leadership. We pose the learning problem of the IM as the solution of a system of nonlinear functional equations obtained from the invariance equation, which characterizes the manifold as a functional relationship between the exogenous state and the state variables, via a hybrid physics-informed machine learning method. In particular, the proposed scheme integrates polynomial series with shallow neural networks, taking advantage of the complementary advantages of both methods. We target low- to medium-dimensional manifolds where polynomial series expansion complexity remains tractable. In this manner, one exploits the interpretability and convergence properties of the polynomial series near the equilibrium where the manifold originates, while farther away, learning of the IM relies on the universal approximation power of neural networks, which are not constrained by a local radius or domain of convergence. The learning problem includes a continuity penalty to match the polynomial and neural network approximations at the regime boundary, and is solved with analytically derived derivatives feeding the Levenberg–Marquardt algorithm. Naturally, depending on the dimensionality of the input-driven system, one may also employ a purely neural network-based IM approximation, for which we also establish a universal approximation theorem based on certain assumptions on system dynamics. We assess the efficiency of the proposed framework using two benchmark examples, namely an enzymatic bioreactor problem and an agent-based car-following model with a leader-follower structure, analyzing convergence behavior, numerical approximation accuracy, and computational training cost. Furthermore, we compare IM approximations obtained using standalone neural networks, standalone polynomial expansions, and the hybrid scheme, demonstrating that the latter hybrid approach outperforms standalone schemes in terms of numerical approximation accuracy. 
\end{abstract}


\section{Introduction}
Invariant manifolds (IMs) capture essential features of the dynamics--such as long-term behavior near equilibria or slow evolution in multiscale systems--and provide a natural framework for reducing the dimensionality of complex systems. By restricting the dynamics on IMs, one can reformulate high-dimensional systems in a lower-dimensional setting without losing critical information \cite{arnold2012geometrical,carr2012applications,guckenheimer2013nonlinear,wiggins2003introduction,kuehn2015multiple}, facilitating stability, bifurcation analysis, construction of ROMs \cite{wiggins2003introduction,guckenheimer2009computing,guckenheimer2013nonlinear,arnold2012geometrical,gorban2005invariant} and design of controllers \cite{castillo1993nonlinear,astolfi2003immersion,hamzi2005controlled,kazantzis2001invariant,kazantzis2003,kazantzis2005singular,kazantzis2005model,byrnes1997output,colonius2012dynamics}.
For continuous-time systems, the theoretical foundations of IMs are well-established. Geometric singular perturbation theory (GSPT) \cite{fenichel1979geometric,jones1995geometric,roberts1989utility,goussis2006efficient,kuehn2015multiple} and spectral submanifold theory \cite{haller2016nonlinear,hirsch1970invariant} provide rigorous tools for slow/fast systems and normally hyperbolic manifolds, with extensions to randomly perturbed systems \cite{benedetti2026numerical}.~A variety of computational/numerical methods based on well established theoretical approaches have been developed \cite{maas1992simplifying,dellnitz1997subdivision,lam1989understanding,cabre2003parameterization,gear2005projecting} such as the Computational Singular Perturbation (CSP) \cite{lam1989understanding,goussis1992study,goussis2006efficient,valorani2005higher},  the Zero-Derivative Principle (ZDP)  \cite{zagaris2009analysis}, the Flow Curvature (FC) \cite{ginoux2008slow}, and Polynomial Chaos for random ordinary differential equations\cite{breden2020computing}. Such approximations enable successful applications in diverse domains, such as chemical kinetics \cite{goussis2011model,valorani2015dynamical,manias2016mechanism,tingas2018chemical,kuehn2025fast}, biochemical networks \cite{patsatzis2016asymptotic,zagaris2009analysis}, structural dynamics \cite{haller2016nonlinear,touze2021model,vizzaccaro2021direct}, as well as in multiscale and complex systems modelling for the construction of coarse-grained models \cite{roberts2014dynamical}. Data-driven approaches for learning low-dimensional stable, unstable and center IMs from high-dimensional microscopic simulators/legacy codes within the Equation-free approach \cite{kevrekidis2003equation} for multiscale computations have also been proposed \cite{gear2005projecting,siettos2014equation,siettos2022numerical,zagaris2004fast,zagaris2009analysis,chin2022enabling,evangelou2023double,patsatzis2023data}. The core of these approaches is the construction of a ``coarse'' discrete map of the emergent dynamics, based on short simulations of microscale simulations. Although classical methods are well established, direct machine-learning-based identification of IMs is still limited to a small and relatively recent body of work. In this direction, the methods proposed to date are primarily data-driven, and include approaches based on SINDy \cite{delgado2026sindy} and deep neural networks \cite{chui2018deep,linot2020deep,ghadami2022deep,chen2024deep,koronaki2024nonlinear}.

Here, building on our previous recent work \cite{siettos2022numerical,patsatzis2024slow,patsatzis2024slowB} for the computation of IMs in the continuous-time domain via physics-informed neural networks (PINNs) \cite{raissi2019physics,karniadakis2021physics}, and on established theory for their existence \cite{kazantzis2001invariant,kazantzis2005model}, we propose a hybrid physics-informed (PI) machine learning framework for approximating IMs in the discrete-time domain focusing on nonlinear systems driven by exogenous systems (i.e. realized by state space representations otherwise known as skew-product systems).~Such systems naturally arise especially in control theory, such as in nonlinear output regulation problems in the presence of appropriately modeled exosystem-generated reference inputs  and/or disturbances \cite{huang2004nonlinear,byrnes2000output,chen2005robust,bin2019adaptive,jiang2020cooperative,kohler2021constrained}. Such dynamics arise also in complex system dynamics conforming to a leader-follower structure \cite{couzin2005effective,cucker2008flocking,ni2010leader,maffettone2026bio}, for instance in problems of pedestrian and traffic dynamics \cite{gazis1961nonlinear,argall2002rigorous,aw2002derivation,cristiani2014multiscale,tordeux2018traffic,chiarello2021multiscale,zhang2024car}, and the collective motion of animals for example in bird flocks \cite{cucker2007emergent,shen2008cucker,gu2009leader,li2010cucker,dalmao2011cucker,cristiani2021all}, where follow-the-leader mechanisms are used to describe the influence of a small number of ``informed'' agents on a larger population of followers. In this setting, the leaders play the role of an exogenous or reference subsystem, while the followers adjust their motion in response to the leaders’ trajectories and to local interactions within the group. For example, in animal groups, the leaders are assumed to possess information about a preferred direction, corresponding, for instance, to the location of a known resource or to a portion of a migration pathway \cite{couzin2005effective}. The assumption that the dynamics are driven by ``informed'' autonomous leaders implies the existence of a low- to medium-dimensional slow invariant manifold, with dimension controlled by the small number of leaders. This makes the proposed hybrid scheme particularly suitable for such systems, as the combinatorial complexity of polynomial approximations is still tractable: near the equilibrium, where the manifold emanates, the scheme employs polynomial series to exploit the corresponding convergence properties, while farther away it relies on the universal approximation power of neural networks, which are not constrained by a local radius or region of convergence. Naturally, for higher-dimensional manifolds, a standalone PINN approximation may also be employed. For this case, we also establish the assumptions required for a universal approximation theorem and derive the corresponding convergence rates of the PINN.
The learning task is formulated as the solution of a system of nonlinear functional equations derived from the associated invariance equation.  In this setting, the evolution of the full (large-scale) state is effectively slaved to a small number of driving coordinates, so that the long-term behavior can be represented on a reduced set of variables. This perspective aligns with the central assumption of reduced-order modeling: despite the apparent high dimensionality and complexity of the governing dynamics, the relevant trajectories often concentrate near low-dimensional structures that provide an accurate and computationally efficient representation of the system’s effective dynamics.

In this setting, the performance of the proposed PI hybrid scheme is illustrated and its performance is assessed via two illustrative benchmark examples: an enzymatic bioreactor problem \cite{kazantzis2005model} and a car-following problem with an autonomous leader vehicle \cite{bando1995dynamical}.~An additional validation example with known analytic IM is provided in the Appendix.~Comparative analysis of the IM approximations provided by the hybrid scheme demonstrates high approximation accuracy over those provided by the traditional series expansion approach \cite{kazantzis2005model}.~More importantly, the hybrid scheme outperforms the IM approximations constructed purely by polynomial series expansions or PINNs, particularly far from or near to the region of the equilibrium point, respectively.

The manuscript is organized as follows:~we first discuss the conditions of IM existence in discrete-time systems in \cref{sec:meth} and then present the PI hybrid scheme for the IM approximations. We first present the PINN component, which may also be used as a standalone approximation method, and state the corresponding universal approximation theorem along with the assumptions under which it holds. Then, we briefly present the polynomial series approximation of the manifold around the equilibrium and introduce the optimization problem of the PI hybrid scheme. Details regarding the implementation of the full hybrid scheme are presented in the Appendix, where we also provide analytical expressions for all derivatives required to solve the corresponding optimization problem using the Levenberg–Marquardt method.~In \cref{sec:Bench}, we briefly present the two illustrative examples, and then we assess the performance of the hybrid scheme as well as the other schemes based on the numerical results in \cref{sec:NumR}.~We finally offer a few concluding remarks focusing on the proposed method's comparative advantages and limitations in \cref{sec:Con}, providing also directions for future work.

\section{Methodology}
\label{sec:meth}
Let us consider the non-linear discrete-time dynamical system with a state space representation:
\begin{equation}
    \mathbf{x}(k+1)  = \mathbf{F}(\mathbf{x}(k),\mathbf{y}(k)), \qquad
    \mathbf{y}(k+1) =  \mathbf{G}(\mathbf{y}(k)), \label{eq:genDDS1} 
\end{equation}
where the state vector $\mathbf{x}\in\mathbb{R}^N$,  $\mathbf{F}:\mathbb{R}^N \times \mathbb{R}^M\rightarrow \mathbb{R}^N$ is driven by the leader/informed agents dynamics of $\mathbf{y}\in\mathbb{R}^M$ generated by the ``exosystem/leading system'' $\mathbf{G}:\mathbb{R}^M \rightarrow \mathbb{R}^M$ .~Let $\mathbf{x}(k)=[x_1(k),\ldots,x_N(k)]^\top$ denote the driven/follower and $\mathbf{y}(k) = [y_1(k),\ldots,y_M(k)]^\top$ denote the driving/leading state variables at the $k\in\mathbb{N}$ discrete time step.~Without loss of generality, we assume that the origin $(\mathbf{x}_0,\mathbf{y}_0) = (\mathbf{0}^N,\mathbf{0}^M)$ is an equilibrium point of the system in \cref{eq:genDDS1}, and in line with \cite{kazantzis2005model}, we make the following assumption to allow for polynomial expansions around the equilibrium:
\begin{assumption} \label{ass1}
    The matrix $\mathbf{A}=\partial_\mathbf{y} \mathbf{G} (\mathbf{y}_0) \in \mathbb{R}^{M\times M}$ has non-zero eigenvalues, denoted by $k_i$ for $i=1,\ldots,M$, all of which lie either strictly inside or outside the unit disc.~This assumption also ensures that the $\mathbf{G}(\mathbf{y})$ is locally invertible in a neighborhood of $\mathbf{y}_0$.
\end{assumption}
The system in \cref{eq:genDDS1} can be rewritten in the form:
\begin{equation}
    \mathbf{x}(k+1)  = \mathbf{B} \mathbf{x}(k) + \mathbf{C} \mathbf{y}(k) + \mathbf{f}(\mathbf{x}(k),\mathbf{y}(k)), \qquad
    \mathbf{y}(k+1)  =  \mathbf{A} \mathbf{y}(k) + \mathbf{g}(\mathbf{y}(k)), \label{eq:genDDS2} 
\end{equation}
where $\mathbf{B} =\partial_\mathbf{x} \mathbf{F} (\mathbf{x}_0,\mathbf{y}_0) \in \mathbb{R}^{N\times N}$ and $\mathbf{C} =\partial_\mathbf{y} \mathbf{F}(\mathbf{x}_0,\mathbf{y}_0) \in \mathbb{R}^{N\times M}$ are constant matrices derived from the linearization of the system around the equilibrium.~The non-linear terms $\mathbf{f}:\mathbb{R}^N \times \mathbb{R}^M\rightarrow \mathbb{R}^N$ and $\mathbf{g}:\mathbb{R}^M \rightarrow \mathbb{R}^M$ are real-valued functions, which satisfy $\mathbf{f}(\mathbf{x}_0,\mathbf{y}_0) = \mathbf{0}^N$ and $\mathbf{g}(\mathbf{y}_0)=\mathbf{0}^M$, as well as $\partial_\mathbf{x} \mathbf{f} (\mathbf{x}_0,\mathbf{y}_0) = \partial_\mathbf{y} \mathbf{f} (\mathbf{x}_0,\mathbf{y}_0) = \partial_\mathbf{y} \mathbf{g}(\mathbf{y}_0) = \mathbf{0}$ with consistent zero-matrix dimensions. Based on the above assumption, it has been shown that the system in \cref{eq:genDDS1} exhibits a unique locally \emph{invariant manifold} (IM) under the assumptions of Theorem 1 in \cite{kazantzis2005model}, which we restate below (for the proof, see \cite{kazantzis2005model}):
\begin{theorem}[Invariant Manifold existence and uniqueness, \cite{kazantzis2005model}] \label{th:th1Kaz}
    Consider the non-linear discrete dynamical system in \cref{eq:genDDS1}, where \Cref{ass1} is satisfied.~Additionally, assume that the eigenvalues of the matrix $\mathbf{A}=\partial_\mathbf{y} \mathbf{G}(\mathbf{y}_0)$, denoted by $k_i$ for $i=1,\ldots,M$ are not related to the eigenvalues of the matrix $\mathbf{B}=\partial_\mathbf{x} \mathbf{F}(\mathbf{x}_0,\mathbf{y}_0)$, denoted by $\lambda_j$ for $j=1,\ldots,N$, through any resonance condition of the form:
    \begin{equation}
        \prod_{i=1}^M k_i^{d_i} = \lambda_j,
    \end{equation}
    where $d_i$ are non-negative integers satisfying $\sum_{i=1}^M d_i>0$.~Then, there exists a compact neighborhood $\mathcal{V}\subset \mathbb{R}^M$ around the equilibrium $\mathbf{y}_0$, and a unique, locally analytic mapping $\boldsymbol{\pi}:\mathcal{V}\rightarrow\mathbb{R}^N$ such that:
    \begin{equation}
    \mathcal{M} = \{(\mathbf{x},\mathbf{y}) \in \mathbb{R}^N \times \mathcal{V}: \mathbf{x} = \boldsymbol{\pi} (\mathbf{y}), \boldsymbol{\pi}(\mathbf{y}_0) = \mathbf{x}_0 \}, 
    \label{eq:IM}
    \end{equation}
    defines an analytic local invariant manifold (IM) of the system in \cref{eq:genDDS1}.~This manifold passes through the equilibrium point $(\mathbf{x}_0,\mathbf{y}_0) = (\mathbf{0}^N,\mathbf{0}^M)$.~Furthermore, the mapping $\boldsymbol{\pi} (\mathbf{y})$ is the unique solution of the system of nonlinear functional equations (NFEs) associated with the dynamic equations \cref{eq:genDDS2}:
    \begin{equation}
    \boldsymbol{\pi} (\mathbf{A} \mathbf{y} + \mathbf{g}(\mathbf{y})) = \mathbf{B} \boldsymbol{\pi} (\mathbf{y}) + \mathbf{C} \mathbf{y} + \mathbf{f}(\boldsymbol{\pi} (\mathbf{y}),\mathbf{y}), \qquad \text{where} \qquad \boldsymbol{\pi} (\mathbf{y}_0) = \mathbf{x}_0.
    \label{eq:NFEs_IM}
\end{equation}
\end{theorem}
The system of NFEs in \cref{eq:NFEs_IM} arises from the invariance of the manifold $\mathcal{M}$ in \cref{eq:IM}.~Specifically, when the solution of the discrete system in \cref{eq:genDDS2} at the $k$-th step is on $\mathcal{M}$ (i.e., when $\mathbf{x}(k) = \boldsymbol{\pi}(\mathbf{y}(k))$), then the invariance ensures that the system remains on $\mathcal{M}$ at the next $(k+1)$-th step.~This implies $\mathbf{x}(k+1) = \boldsymbol{\pi}(\mathbf{y}(k+1))$,  which, upon substitution of the system dynamics from \cref{eq:genDDS2}, yields \cref{eq:NFEs_IM}.

According to \cref{th:th1Kaz}, the functional map $\mathbf{x} = \boldsymbol{\pi}(\mathbf{y})$ of the IM $\mathcal{M}$ in \cref{eq:IM} can be computed through the solution of the system of NFEs in \cref{eq:NFEs_IM}.~Due to their implicit and  nonlinear nature, exact solutions are generally intractable, necessitating approximate solutions.~Traditional methods \cite{haro2016parameterization,kazantzis2001invariant,kazantzis2005model} use local power series expansions (PSE) around $(\mathbf{x}_0,\mathbf{y}_0)$, where the coefficients are determined by solving a linear system of algebraic equations; see \cref{app:IM_PS} for a detailed formulation.~While polynomial approximations are rigorous, interpretable, and can achieve exponential convergence for analytic IMs, they face two critical limitations: 
\begin{enumerate}[label=(\roman*)]
    \item The number of coefficients grows combinatorially with the system dimension and degree of approximation. Although modern symbolic and automatic-differentiation tools can construct high-dimensional multivariate polynomial approximations efficiently, they do not eliminate the underlying combinatorial growth in the number of coefficients: a total-degree-$h$ polynomial representation in $M$ variables requires $\binom{M+h}{h}$ basis terms. Therefore, PSEs for IM approximations still encounter a ``curse of dimensionality"-type of bottleneck unless sparsity or other low-complexity structure is imposed \cite{lokshtanov2017beating}. 
    \item Also importantly, their effectiveness is limited by the analyticity domain of the target manifold. Indeed, the radius of convergence of Taylor expansions and the Bernstein-ellipse analyticity regions of spectral approximations such as Chebyshev and Legendre polynomials, can cause PSEs to fail outside these regions, which are not known a-priori \cite{trefethen2019approximation}.
\end{enumerate}
These challenges also arise in continuous-time fast/slow dynamical systems in the approximation of slow IMs.~In such systems, the invariance equation comes in a form of a partial differential equation (PDE), the solution of which may be approximated by asymptotic expansions \cite{gorban2003method,kuehn2015multiple,roussel1991geometry}, or by advanced numerical methods \cite{gear2005projecting,goussis1992study,goussis2006efficient,maas1992simplifying}.~In prior work \cite{patsatzis2024slow,patsatzis2024slowB}, we demonstrated that PINNs provide very accurate slow IM approximations, especially near the boundaries of the manifold, where the asymptotic series diverges.

Motivated by the above challenges and earlier work, we propose here a physics-informed \textit{hybrid scheme} for discrete-time systems driven by exogenous autonomous input dynamics that combines the complementary advantages of both approximation approaches, i.e., multivariate polynomial series and PINNs. In particular, both components are universal approximators of continuous functions on compact domains, as guaranteed by the Weierstrass Approximation Theorem \cite{rudin1964principles,weierstrass1885analytische} for polynomials and the universal approximation theorems for NNs \cite{cybenko1989approximation,hornik1989multilayer}. Polynomial approximations are interpretable and often computationally effective in low- and medium-dimensional settings and, for analytic invariant manifolds, can achieve geometric or exponential convergence with respect to polynomial degree \cite{rivlin1981introduction,trefethen2019approximation}. Their effectiveness, however, is limited by the combinatorial growth of coefficients with dimensionality and approximation order, as well as by the analyticity domain of the target manifold: for local power-series expansions, this is reflected in a finite radius of convergence, while for Chebyshev and Legendre approximations it is governed by the largest admissible Bernstein ellipse, which is typically unknown a priori. By contrast, NNs scale more favorably with dimension and apply under substantially weaker regularity assumptions, albeit at the cost of a more challenging NP-hard training problem even for shallow architectures with fixed input dimension \cite{froese2023training}.

We reiterate that our focus here is on IMs of relatively low to moderate dimension, in the spirit of the assumptions that commonly underlie the construction of low-dimensional reduced-order models (ROMs) for high-dimensional systems. Accordingly, the proposed hybrid scheme is particularly well-suited for this standard setting.
Of course, depending on the problem under consideration and the region where the IM is sought, the two approximation strategies may also be employed as standalone schemes.

\begin{figure}[!h]
    \centering
    \includegraphics[width=0.95\linewidth]{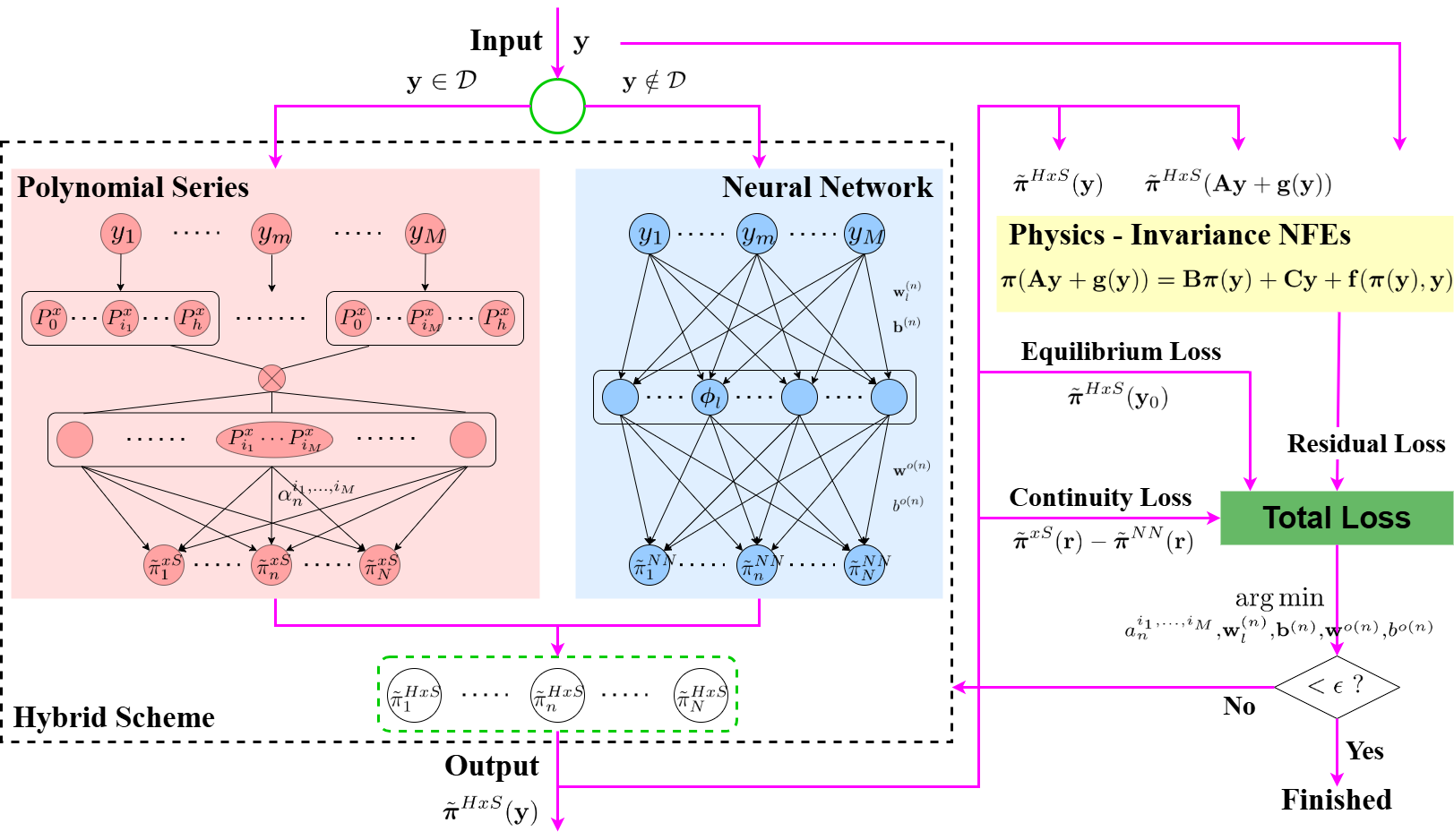}
    \caption{Schematic of the PI hybrid scheme for approximating IM functionals of discrete-time dynamical systems in the form of \cref{eq:genDDS1}.~The scheme combines polynomial series inside a hyperrectangle $\mathcal{D}$ of radius $\mathbf{r}$ centered at the equilibrium and NNs outside $\mathcal{D}$.~The hybrid IM scheme is trained to minimize a composite loss function that enforces the system of NFEs in \cref{eq:NFEs_IM}, including residuals for manifold invariance, equilibrium constraints and continuity constraints at $\partial \mathcal{D}$. The coefficients of the polynomial series and the parameters of the NN are jointly optimized to minimize the loss.}
    \label{fig:Sch1}
\end{figure}

A schematic of the proposed scheme is shown in \cref{fig:Sch1}.
In particular, we approximate the IM as  
\begin{equation}
\mathbf{x} = \tilde{\boldsymbol{\pi}}^{HxS}(\mathbf{y},\boldsymbol{a},\mathbf{p};\mathcal{H}_1,\mathcal{H}_2,\mathbf{r}) \in \mathbb{R}^N,
\label{eq:IM_hybrid0}
\end{equation}
the elements of which for $n=1,\ldots,N$ are expressed as:
\begin{equation}
    \tilde{\pi}^{HxS}_n(\mathbf{y},\boldsymbol{a}_n,\mathbf{p}_n;\mathcal{H}_1,\mathcal{H}_2,\mathbf{r}) = \tilde{\pi}^{xS}_n(\mathbf{y},\boldsymbol{a}_n;\mathcal{H}_1)^{H(\mathbf{y},\mathbf{r})} \cdot  \tilde{\pi}^{NN}_n(\mathbf{y},\mathbf{p}_n;\mathcal{H}_2)^{1-H(\mathbf{y},\mathbf{r})}, 
    \label{eq:IM_hybrid}
\end{equation}
where $\tilde{\pi}^{xS}_n(\mathbf{y},\boldsymbol{a}_n;\mathcal{H}_1)$ and $\tilde{\pi}^{NN}_n(\mathbf{y},\mathbf{p}_n;\mathcal{H}_2)$ are the $n$-th outputs of a polynomial series and a NN, as functions of $\mathbf{y}$, respectively.~These components include the parameters $\boldsymbol{a}_n$ and $\mathbf{p}_n$ and the hyperparameters $\mathcal{H}_1$ and $\mathcal{H}_2$, which will be discussed next.
~The function $H:\mathbb{R}^M\times\mathbb{R}^M\rightarrow\mathbb{R}$ in \cref{eq:IM_hybrid} is a product of Heaviside functions defined as:
\begin{equation}
    H(\mathbf{y},\mathbf{r}) = \prod_{i=1}^M H(r_i^2-y_i^2) = \begin{cases}
        1, & \lvert y_i\rvert < r_i, ~~\forall i=1,\ldots,M \\
        0, & \text{otherwise}.
    \end{cases},
    \label{eq:Heaviside}
\end{equation}
which takes the values $H(\mathbf{y},\mathbf{r})=1$ when $\mathbf{y}\in \mathcal{D}$ and $H(\mathbf{y},\mathbf{r})=0$ when $\mathbf{y}\notin \mathcal{D}$.~This ensures that the hybrid scheme approximation $\mathbf{x} = \tilde{\boldsymbol{\pi}}^{HxS}(\mathbf{y},\boldsymbol{a},\mathbf{p};\mathcal{H}_1,\mathcal{H}_2,\mathbf{r})$ in \cref{eq:IM_hybrid} reduces to 
a polynomial series (power, Legendre or Chebyshev) approximation within $\mathcal{D}$ and a NN approximation outside $\mathcal{D}$.~We hereby note that the hybrid approximation depends 
also on the radius $\mathbf{r}$ of the hyperrectangle, which is considered in \cref{eq:IM_hybrid} as an additional hyperparameter.  

In the ensuing sections, we first discuss the IM approximation by a NN, establishing the assumptions required for a universal approximation theorem, and deriving the corresponding convergence rates. Subsequently, we present the polynomial approximation scheme and, finally, the associated optimization problem.

\subsection{Physics-informed neural network approximation of the IM} 
The NN component of the proposed PI hybrid scheme is activated in the region where the uniform analyticity assumption of the power series/polynomial approximation fails (outside $\mathcal{D}$; see \cref{fig:Sch1}). Therein, we consider a feedforward artificial NN $\tilde{\boldsymbol{\pi}}^{NN}(\mathbf{y},\mathbf{p};\mathcal{H}_2) \in \mathbb{R}^N$ with inputs $\mathbf{y}\in \mathbb{R}^M$, where $\mathbf{p}$ is a vector including the parameters of the NN (weights and biases of each layer) and $\mathcal{H}_2$ includes the NN hyperparameters, such as the number of hidden layers, the type and parameters of activation function, etc..~For simplicity, we assume a single layer feedforward NN (SLFNN) with $L$ neurons in the hidden layer (more expressive architectures were tested and did not show significant numerical accuracy improvement).~Each output of the SLFNN, $\tilde{\pi}^{NN}_n(\mathbf{y},\mathbf{p}_n;\mathcal{H}_2)$ for $n=1,\ldots,N$, is expressed as:
\begin{equation}
    \tilde{\pi}^{NN}_n(\mathbf{y},\mathbf{p}_n;\mathcal{H}_2) = \mathbf{w}^{o(n)\top} \boldsymbol{\phi} \left( \mathbf{W}^{(n)} \mathbf{y}  + \mathbf{b}^{(n)} \right)  + b^{o(n)},
    \label{eq:IM_NN}
\end{equation}
where $\mathbf{p}_n = [\mathbf{w}^{o(n)},b^{o(n)},\mathbf{W}^{(n)},\mathbf{b}^{(n)}]^\top\in\mathbb{R}^{L(M+2)+1}$ includes: (i) the output weights $\mathbf{w}^{o(n)} = [w^{o(n)}_1, \ldots, w^{o(n)}_L]^\top\in \mathbb{R}^{L}$ connecting the hidden to the output layer, (ii) the bias $b^{o(n)} \in \mathbb{R}$ of the output layer, (iii) the internal weights $\mathbf{W}^{(n)}\in\mathbb{R}^{L\times M}$, where each column $\mathbf{w}^{(n)}_l = [w_{l,1}^{(n)}, \ldots, w_{l,M}^{(n)}]^\top\in \mathbb{R}^{M}$ corresponds to the weights between the input neurons and the $l$-th neuron in the hidden layer, and (iv) the internal biases  $\mathbf{b}^{(n)} = [b_1^{(n)}, \ldots, b_L^{(n)}]^\top\in\mathbb{R}^{L}$ of the neurons in the hidden layer.~The outputs of the activated neurons $\boldsymbol{\phi}_l\left( \mathbf{W}^{(n)} \mathbf{y}  + \mathbf{b}^{(n)} \right)$ for $l=1,\ldots,L$ are included in the column vector $\boldsymbol{\phi}(\cdot) \in\mathbb{R}^L$.~Here, the logistic sigmoid function is used as the activation function due to its universal approximation properties \cite{cybenko1989approximation} and its suitability for symbolic differentiation. 

At this point, we first present the following  Theorem.
\begin{theorem}
\label{th:NN_NFEs}
Let $\mathcal{V} \subset \mathbb{R}^M$ be the compact ball $\mathcal V:=B_R=\{\mathbf{z}\in\mathbb{R}^M:\|\mathbf{z}\|\le R\}$ for some $R>0$, and let $\boldsymbol{\pi}(\mathbf{y}):=\big({\pi_1}(\mathbf{y}),\dots,{\pi_N}(\mathbf{y})\big)^\top$ be the vector-valued IM of the input-driven discrete-time system in \cref{eq:genDDS2} for $\mathbf{y}\in \mathcal{V}$. Let $\mu$ be a probability measure on $\mathcal{V}$ (hence a finite Borel measure with $\mu(\mathcal{V})=1$, following \cite{barron1993universal}) and assume that:
\begin{enumerate}[label=(H\arabic*)]
\item (\emph{Barron regularity}\cite{barron1993universal}) each component $\pi_n$ belongs to the Barron class $\mathcal{B}(\mathcal{V})$, with norm $\|\pi_n\|_{\mathcal{B}}:=\int_{\mathbb{R}^M} \|\boldsymbol{\omega}\|\,|\widehat \pi_n(\boldsymbol{\omega})|\,d\boldsymbol{\omega}<\infty$, where $\widehat \pi_n$ is the Fourier transform of a suitable extension of $\pi_n$ to $\mathbb{R}^M$ for $n=1,2,\dots N$, and define $\|\boldsymbol{\pi}\|_{\mathcal{B}}:=\sum_{n=1}^N \|\pi_n\|_{\mathcal{B}}$,
\item (\emph{Lipschitz continuity of $\mathbf{f}$ in $\mathbf{x}$}) for each $\mathbf{y} \in\mathcal V$, there exists $L_{\mathbf{f}}(\mathbf{y})\ge 0$ such that for all $\mathbf{x}_1,\mathbf{x}_2\in\mathbb{R}^{N}$:
\begin{equation}
\label{eq:Lip_f}
\|\mathbf{f}(\mathbf{x}_1,\mathbf{y})-\mathbf{f}(\mathbf{x}_2,\mathbf{y})\|\le L_{\mathbf{f}}(\mathbf{y})\|\mathbf{x}_1-\mathbf{x}_2\|,
\end{equation}
and define $\bar{L}_{\mathbf{f}}:= \sup_{\mathbf{y} \in \mathcal{V}} L_{\mathbf{f}}(\mathbf{y})$, which is finite by compactness of $\mathcal{V}$ and continuity of $\mathbf{f}$, 
\item (\emph{Bounded composition or invariance  of $\mu$ with respect to $\mathbf{G}$ on $L^2(\mu)$}) the measure $\mu$ is chosen according to the dynamics of the exosystem $\mathbf{G}$, such that:
\begin{equation}
\|\boldsymbol{\varphi}\circ \mathbf{G}\|_{L^2(\mu;\mathbb{R}^N)}\le L_{\mathbf{G},\mu}\, \|\boldsymbol{\varphi}\|_{L^2(\mu;\mathbb{R}^N)},\, \quad \forall \boldsymbol{\varphi}\in L^2(\mu;\mathbb{R}^N),
\end{equation}
for some $L_{\mathbf{G},\mu}\ge 0$; in the special case $\mu$ is $\mathbf{G}$-invariant, we have $L_{\mathbf{G},\mu}=1$.
\end{enumerate}
Let $\mathcal{N}_L$ denote the class of SLFNNs of width $L$ with sigmoidal activation functions $\sigma \in C^1(\mathbb{R})$ and a single output. For each $n = 1,\ldots,N$, let $\pi_n^L \in \mathcal{N}_L$ be defined as:
\begin{equation}
\pi_n^L(\mathbf{y})=\sum_{k=1}^{L} w^{(o)n}_k\,\sigma((\mathbf{w}_k^{(n)})^\top \mathbf{y}+\mathbf{b}^{(n)}_k)+b^{o(n)}, 
\end{equation}
aligned with the parameter definition in \cref{eq:IM_NN}. Then, under hypotheses (H1)-(H3), there exists a vector-valued network $\boldsymbol{\pi}^L = (\pi_1^L, \dots, \pi_N^L)^\top \in C^1(\mathcal{V}; \mathbb{R}^N)$ such that 
\begin{enumerate}[label=\Alph*)]
    \item The approximation error for the IM can be made arbitrarily small, with the explicit convergence rate:
\begin{equation}
\label{eq:L2_bar_rate1}
\|\boldsymbol{\pi}^{L}-\boldsymbol{\pi}\|_{L^2(\mu;\mathbb{R}^{N})} \le C_{\sigma,\mathcal V}\,\frac{\|\boldsymbol{\pi}\|_{\mathcal B}}{\sqrt{L}}.
\end{equation}
\item The approximation to the solution of the invariance NFEs is bounded as: 
\begin{equation}
     \left(
\int_{\mathcal V}
\left\|\boldsymbol{\pi}^{L} \big(\mathbf{G}(\mathbf{y})\big) - \mathbf{B} \boldsymbol{\pi}^{L} (\mathbf{y}) - \mathbf{C} \mathbf{y} + \mathbf{f}(\boldsymbol{\pi}^{L} (\mathbf{y}),\mathbf{y})\right\|_{\mathbb{R}^{N}}^{2}
\, d\mu(\mathbf{y})
\right)^{1/2}\leq \left(L_{\mathbf{G},\mu}+\|\mathbf{B}\|_2+\bar{L}_{\mathbf{f}}\right) C_{\sigma,\mathcal V} \frac{\|\boldsymbol{\pi}\|_{\mathcal{B}}}{\sqrt{L}},
\end{equation}
where $\|\mathbf{B}\|_2$ denotes the matrix norm induced by the Euclidean vector norm.
\end{enumerate}
\end{theorem}
\begin{proof}
(A) By hypothesis (H1), for each component $\pi_n:\mathcal{V}\to\mathbb{R}$, there exists a scalar-output network $\pi_n^L\in \mathcal{N}_L$ of width $L$ such that \cite{barron1993universal}:
\begin{equation}
\label{eq:scalar_barron_rate}
\|\pi_n^{L}-\pi_n\|_{L^2(\mu)} \le C_{n,\sigma,\mathcal{V}}\,\dfrac{\|\pi_n\|_{\mathcal{B}}}{\sqrt{L}},
\end{equation}
where $C_{n,\sigma,\mathcal{V}} >0$ is a constant depending on $n$, the activation function $\sigma$, and the domain $\mathcal{V}$.

Let us define $\mathbf{e}: \mathcal V\to\mathbb{R}^N$, as $\mathbf{e}(\mathbf{y})=\boldsymbol{\pi}^{L}(\mathbf{y})-\boldsymbol{\pi}(\mathbf{y})$. Then:
\begin{align}
\|\mathbf{e}\|_{L^2(\mu;\mathbb{R}^{N})}
\le \sum_{n=1}^N \|e_n\|_{L^2(\mu)}
= \sum_{n=1}^N \|\pi_n^{L}-\pi_n\|_{L^2(\mu)}\le \sum_{n=1}^N C_{n,\sigma,\mathcal V}\,\frac{\|\pi_n\|_{\mathcal{B}}}{\sqrt{L}}=\\
\dfrac{1}{\sqrt{L}}\sum_{n=1}^N C_{n,\sigma,\mathcal V}\,\|\pi_n\|_{\mathcal{B}}\le
\dfrac{1}{\sqrt{L}}
\left(\max_{1\le n\le N} C_{n,\sigma,\mathcal V}\right)
\sum_{n=1}^N \|\pi_n\|_{\mathcal{B}}=C_{\sigma,\mathcal{V}}\,\dfrac{\|\boldsymbol{\pi}\|_{\mathcal B}}{\sqrt{L}},
\end{align} 
where $C_{\sigma,\mathcal V}=\max_{1\le n\le N} C_{n,\sigma,\mathcal V}$.

(B) Let's define now the residual for the NFEs of the IM (see \cref{eq:NFEs_IM}) as:
\begin{equation}
\label{eq:residual_def_IM}
\mathbf{r}(\mathbf{y}):=\boldsymbol{\pi}^{L}\big(\mathbf{G}(\mathbf{y})\big)-\mathbf{B}\boldsymbol{\pi}^{L}(\mathbf{y})-\mathbf{C}\mathbf{y}-\mathbf{f}(\boldsymbol{\pi}^{L}(\mathbf{y}),\mathbf{y}).
\end{equation}
Since $\boldsymbol{\pi}(\mathbf{y})$ satisfies exactly the NFEs of the IM for all $\mathbf{y}\in\mathcal V$, we have that:
\begin{align*}
\mathbf{r}(\mathbf{y})
&=\Big(\boldsymbol{\pi}^{L}\big(\mathbf{G}(\mathbf{y})\big)-\mathbf{B}\boldsymbol{\pi}^{L}(\mathbf{y})
-\mathbf{C}\mathbf{y}-\mathbf{f}(\boldsymbol{\pi}^{L}(\mathbf{y}),\mathbf{y})\Big)
-\Big(\boldsymbol{\pi}\big(\mathbf{G}(\mathbf{y})\big)-\mathbf{B}\boldsymbol{\pi}(\mathbf{y})
-\mathbf{C}\mathbf{y}-\mathbf{f}(\boldsymbol{\pi}(\mathbf{y}),\mathbf{y})\Big)\\
&=\Big(\boldsymbol{\pi}^{L}\big(\mathbf{G}(\mathbf{y})\big)-\boldsymbol{\pi}\big(\mathbf{G}(\mathbf{y})\big)\Big)
-\mathbf{B}\big(\boldsymbol{\pi}^{L}(\mathbf{y})-\boldsymbol{\pi}(\mathbf{y})\big)
-\big(\mathbf{f}(\boldsymbol{\pi}^{L}(\mathbf{y}),\mathbf{y})-\mathbf{f}(\boldsymbol{\pi}(\mathbf{y}),\mathbf{y})\big)\\
&= \mathbf{e}\big(\mathbf{G}(\mathbf{y})\big) - \mathbf{B}\mathbf{e}(\mathbf{y})
-\big(\mathbf{f}(\boldsymbol{\pi}^{L}(\mathbf{y}),\mathbf{y})-\mathbf{f}(\boldsymbol{\pi}(\mathbf{y}),\mathbf{y})\big).
\end{align*}
Hence:
\begin{equation}
\|\mathbf{r}(\mathbf{y})\|
\le \|\mathbf{e}\big(\mathbf{G}(\mathbf{y})\big)\| + \|\mathbf{B}\|\,\|\mathbf{e}(\mathbf{y})\|
     + \|\mathbf{f}(\boldsymbol{\pi}^{L}(\mathbf{y}),\mathbf{y})-\mathbf{f}(\boldsymbol{\pi}(\mathbf{y}),\mathbf{y})\|.
     \label{eq:rnorm}
\end{equation}
By the Lipschitz continuity hypothesis (H2) of $\mathbf{f}$, we get:
\begin{equation}
\|\mathbf{f}(\boldsymbol{\pi}^{L}(\mathbf{y}),\mathbf{y})-\mathbf{f}(\boldsymbol{\pi}(\mathbf{y}),\mathbf{y})\|\le L_\mathbf{f}(\mathbf{y})\,\|\boldsymbol{\pi}^{L}(\mathbf{y})-\boldsymbol{\pi}(\mathbf{y})\|=L_\mathbf{f}(\mathbf{y})\,\|\mathbf{e}(\mathbf{y})\|.
\end{equation}
Taking the $L^2(\mu;\mathbb{R}^N)$ norm on both sides of  \cref{eq:rnorm} yields:
\begin{equation}
\|\mathbf{r}\|_{L^2(\mu;\mathbb{R}^{N})}
\le \|\mathbf{e}\circ\mathbf{G}\|_{L^2(\mu;\mathbb{R}^{N})} + (\|\mathbf{B}\|_2+\bar{L}_{\mathbf{f}})\,\|\mathbf{e}\|_{L^2(\mu;\mathbb{R}^{N})}.
\label{eq:rnorm1}
\end{equation}
By hypothesis (H3) about the bounded composition of $\mu$ with respect to $\mathbf{G}$, we have:
\begin{equation*}
\|\mathbf{e}\circ \mathbf{G}\|_{L^2(\mu;\mathbb{R}^N)}\le L_{\mathbf{G},\mu}\, \|\mathbf{e}\|_{L^2(\mu;\mathbb{R}^N)},
\end{equation*}
which upon substitution this into \cref{eq:rnorm1} yields:
\begin{equation}
\|\mathbf{r}\|_{L^2(\mu;\mathbb{R}^{N})}
\le (L_{\mathbf{G},\mu}+\|\mathbf{B}\|_2+\bar{L}_{\mathbf{f}})\,\|\mathbf{e}\|_{L^2(\mu;\mathbb{R}^{N})}.
\end{equation}
Finally, applying the bound from part (A) yields (B).
\end{proof}

\subsection{Polynomial series approximation of the IM} 
The polynomial series component of the proposed PI hybrid scheme falls entirely within \cref{th:th1Kaz}, implying local manifold analyticity near equilibrium.~To model this component, we adopt a NN-like architecture inspired by \cite{tang2023physics} (see \cref{fig:Sch1}).~We consider a multivariate polynomial series for the external variables $\mathbf{y}\in \mathbb{R}^M$ of degree $h$, denoted as $\tilde{\boldsymbol{\pi}}^{xS}(\mathbf{y},\boldsymbol{a};\mathcal{H}_1) \in \mathbb{R}^N$, where $\boldsymbol{a}$ contains the polynomial coefficients and $\mathcal{H}_1$ includes both the degree $h$ and polynomial type $x$ (indicated by superscript $xS$).~Each output dimension $n$ (for $n=1,\ldots,N$) of the vector valued $\tilde{\boldsymbol{\pi}}^{xS}$ is given by:
\begin{equation}
    \tilde{\pi}^{xS}_n(\mathbf{y},\boldsymbol{a}_n;\mathcal{H}_1) =  \sum_{i_1,\ldots,i_M=0}^h \alpha_n^{i_1,\ldots,i_M} P^x_{i_1}(y_1) \cdots P^x_{i_M}(y_M) + \mathcal{O}(\mathbf{P}^x(\mathbf{y})^{h+1}),
    \label{eq:IM_polySexp}
\end{equation}
where $\alpha_n^{i_1,\ldots,i_M}$ are the expansion coefficients, $P^x_{i_j}(y_j)$ represents the $i_j$-th degree polynomial of type $x$ for the variable  $y_j$, and the reminder term $\mathcal{O}(\mathbf{P}^x(\mathbf{y})^{h+1})$ captures higher-order contributions. Although the summation includes terms beyond degree $h$ in its formulation (chosen to mimic a NN-like structure, as shown in \cref{fig:Sch1}), we enforce truncation by setting $\alpha_n^{i_1,\ldots,i_M}=0$ for $i_1 + i_2 + \ldots + i_M>h$.~The coefficient vector $\boldsymbol{a}_n$ collects all $\binom{h+M}{M}$ coefficients $\alpha_n^{i_1,\ldots,i_M}$ up to degree $h$, with the complete parameter vector $\boldsymbol{a}=[\boldsymbol{a}_1,\ldots,\boldsymbol{a}_N]^\top$ in $\tilde{\boldsymbol{\pi}}^{xS}(\mathbf{y},\boldsymbol{a};\mathcal{H}_1)$ aggregating all output dimensions $n=1,\ldots,N$.

We implement three polynomial variants in this framework: power series, Legendre and Chebyshev (of the second kind) polynomial series, denoted by the superscript $x=P,L,C$ in \cref{eq:IM_polySexp}.~The power series polynomials match the traditional PSE approaches \cite{kazantzis2005model} with monomials $P_l^{P}(y_m) = y_m^l$ of degree $l\in\mathbb{N}$ for the variable $y_m$ ($m=1,\ldots,M$).~We further consider Legendre and Chebyshev polynomials due to their increased accuracy and reduced computational costs in PINN forward problems \cite{mall2017single,tang2023physics,yang2020neural}.~Both polynomial types are orthogonal with respect to a weight function $\rho(y_m)$ and are defined on the interval $y_m\in[-1,1]$ for every $m=1,\ldots,M$.~Specifically, any two univariate Legendre/Chebyshev polynomials of degrees $i,j\in\mathbb{N}$ with $i\neq j$ satisfy the orthogonality relation:
\begin{equation}
\langle P_i(y_m) \; , \; P_j(y_m) \rangle = \int_{-1}^{1} P_i(y_m)P_j(y_m) \dfrac{dy_m}{\rho(y_m)} = 0.
\end{equation}
For Legendre polynomials with weight function $\rho(y_m)=1$, the $l$-th degree monomial is given by the recursive formula:
\begin{equation}
    P^{L}_0(y_m) = 1, \quad P^{L}_1(y_m) = y_m, \quad P^{L}_{l+1}(y_m) = \dfrac{(2l+1)y_m P^{L}_l(y_m)-lP^{L}_{l-1}(y_m)}{l+1},
    \label{eq:LS_rec}
\end{equation}
and for Chebyshev polynomials of the second kind with weight function $\rho(y_m)=\sqrt{1-y_m^2}$, the $l$-th degree monomial is given by:
\begin{equation}
    P^{C}_0(y_m) = 1, \quad P^{C}_1(y_m) = 2y_m, \quad P^{C}_{l+1}(y_m) = 2y_mP^{C}_l(y_m)-P^{C}_{l-1}(y_m). 
    \label{eq:CS_rec}
\end{equation}

\subsection{The physics-informed optimization problem} 
\label{sb:PI_opt_hyb}
To approximate the invariant manifold (IM) functional defined in \cref{eq:IM}, we seek a hybrid scheme, given by \cref{eq:IM_hybrid0,eq:IM_hybrid}, that satisfies the system of NFEs in \cref{eq:NFEs_IM}, in accordance with \cref{th:th1Kaz}. This is achieved by training the scheme in a physics-informed manner, ensuring that the discovered functional is indeed an invariant manifold.~Let $\mathbf{y}^{(q)}\in\Omega$ be a set of $q=1,\ldots,Q$ collocation points, selected from the domain $\Omega\subseteq\mathcal{V}$ where the IM approximation is sought.~Additionally, let $\mathbf{y}^{(r)}\in \partial \mathcal{D}$ be a set of $r=1,\ldots,R$ collocation points on the boundary of the hyperrectangle $\mathcal{D}$ of radius $\mathbf{r}$.~The parameters $\boldsymbol{a}$ and $\mathbf{p}$ of the hybrid approximation $\tilde{\boldsymbol{\pi}}^{HxS}(\mathbf{y},\boldsymbol{a},\mathbf{p})$ (the dependence on hyper-parameters $\mathcal{H}_1$, $\mathcal{H}_2$ and $\mathbf{r}$ is dropped here for conciseness) are then determined by solving the following optimization problem:
\begin{equation}
    \argmin_{\boldsymbol{a},\mathbf{p}} \mathcal{L}(\boldsymbol{a},\mathbf{p}) =  \argmin_{\boldsymbol{a},\mathbf{p}} \left( \omega_{\Omega} \sum_{q=1}^Q\lVert \mathcal{F}^{HxS}_{q}(\boldsymbol{a},\mathbf{p})\rVert^2 + \omega_{0} \lVert \mathcal{F}^{HxS}_{0}(\boldsymbol{a})\rVert^2 + \omega_{\partial \mathcal{D}} \sum_{r=1}^R\lVert \mathcal{F}^{HxS}_{r}(\boldsymbol{a},\mathbf{p})\rVert^2\right) 
    \label{eq:min_hybrid}
\end{equation}
where the residuals are defined as:
\begin{align}
  \mathcal{F}^{HxS}_{q}(\boldsymbol{a},\mathbf{p}) & =  \tilde{\boldsymbol{\pi}}^{HxS}(\mathbf{A} \mathbf{y}^{(q)} + \mathbf{g}(\mathbf{y}^{(q)}),\boldsymbol{a},\mathbf{p})  - \mathbf{B} \tilde{\boldsymbol{\pi}}^{HxS}(\mathbf{y}^{(q)},\boldsymbol{a},\mathbf{p}) - \mathbf{C} \mathbf{y}^{(q)}- \mathbf{f}(\tilde{\boldsymbol{\pi}}^{HxS}(\mathbf{y}^{(q)},\boldsymbol{a},\mathbf{p}),\mathbf{y}^{(q)}), \nonumber \\
  \mathcal{F}^{HxS}_{0}(\boldsymbol{a}) & = \tilde{\boldsymbol{\pi}}^{HxS}(\mathbf{0}^{M},\boldsymbol{a},\mathbf{p}) = \tilde{\boldsymbol{\pi}}^{xS}(\mathbf{0}^{M},\boldsymbol{a}), \nonumber \\
  \mathcal{F}^{HxS}_{r, \partial \mathcal{D}}(\boldsymbol{a},\mathbf{p}) & = \tilde{\boldsymbol{\pi}}^{xS}(\mathbf{y}^{(r)},\boldsymbol{a}) -  \tilde{\boldsymbol{\pi}}^{NN}(\mathbf{y}^{(r)},\mathbf{p}).
  \label{eq:nonlin_res}
\end{align} 
The first term of the loss function in \cref{eq:min_hybrid} is a physics-informed term enforcing the satisfaction of the NFEs by $\tilde{\boldsymbol{\pi}}^{HxS}(\cdot)$, while the second term ensures that $\tilde{\boldsymbol{\pi}}^{HxS}(\cdot)$ passes through the equilibrium (which can be viewed as a loss function at the boundary).~Note that the hybrid scheme approximation at the equilibrium consists only of the polynomial series component, since $\mathbf{0}^M \in \mathcal{D}$.~The third term is a $C^0$ continuity in the boundary $\partial \mathcal{D}$ ($C^k$ continuity conditions can be also imposed), as it matches the polynomial series component with the NN component on the boundary collocation points $\mathbf{y}^{(r)} \in \partial \mathcal{D}$.~The second and third terms in the loss function serve for pinning the polynomial series and the NN components respectively, which are further balanced by $\omega_{\Omega}, \omega_0, \omega_{\partial \mathcal{D}} \in \mathbb{R}$.

\begin{remark}
\label{rm:1}
A key feature of the hybrid scheme in \cref{eq:min_hybrid} is that the IM approximation $\tilde{\boldsymbol{\pi}}^{HxS}(\cdot)$ at $\mathbf{y}^{(q)}$ may be evaluated by the NN, while at $\mathbf{A}\mathbf{y}^{(q)} + \mathbf{g}(\mathbf{y}^{(q)})$ by the polynomial series and vice-versa, depending on whether these points lie inside or outside $\mathcal{D}$. The radius of the hyperrectangle $\mathcal{D}$ is an additional hyperparameter that needs tuning, especially in the case of power series polynomials, since in the case of Legendre or Chebyshev polynomials the natural choice is $r_m=1$ for $m=1,\ldots,M$.
\end{remark}

The solution of the optimization problem in \cref{eq:min_hybrid} is obtained by minimizing the nonlinear residuals defined in \cref{eq:nonlin_res}, which corresponds to solving a nonlinear least-squares problem; the detailed formulation is provided in \cref{app:OptPIHS_LM}.~To solve this least-squares problem, we employ the Levenberg–Marquardt (LM) gradient-based optimization algorithm \cite{hagan1994training}.~Importantly, we have derived closed-form analytical expressions for all derivatives of the residuals with respect to the parameters of the hybrid scheme, as required by the LM algorithm (see \cref{app:LMalg}).~While the LM algorithm can also be implemented using finite-difference (FD) approximations or automatic differentiation (AD), we emphasize that providing the necessary derivatives analytically yields a more reliable and efficient \textit{numerical analysis-informed} procedure, particularly useful in regimes where the residual is stiff, or highly sensitive to perturbations.~In our previous work \cite{patsatzis2024slow}, where we systematically compared analytical derivatives, AD, and FD in the identification of slow invariant manifolds for deterministic stiff ODE systems, we found that the analytical derivatives consistently outperformed both AD and FD in terms of accuracy and computational cost.~Motivated by these findings, we adopt the same strategy here: exact derivative formulas facilitate the LM optimization and lead to more robust and scalable training/parameter identification within the proposed hybrid framework. 

 
Finally, we note that the optimization framework in \cref{eq:min_hybrid} readily accommodates standalone polynomial series and NN schemes, as the required derivative expressions are already available from the hybrid scheme formulation.~In both cases, the third term in the loss function (which couples the polynomial series and NNs on the boundary $\partial \mathcal{D}$) is omitted by setting $\omega_{\partial \mathcal{D}}=0$.~For the standalone NN scheme, the second term becomes $\mathcal{F}^{NN}_0(\mathbf{p})=\tilde{\boldsymbol{\pi}}(\mathbf{0}^M,\mathbf{p})$, ensuring that the NN approximation passes through the equilibrium.~Importantly, for the standalone NN case, we provide approximation error and convergence guarantees in \cref{th:NN_NFEs}.~For completeness, the solution of the standalone NN optimization problem is discussed in \cref{app:OptPIHS_LM}, and the implementation of the LM algorithm for both the hybrid and standalone NN schemes is detailed in \cref{app:LMalg}, along with a pseudo-code.

\begin{remark}
    \label{rm:2}
    We emphasize that the PI optimization problem in \cref{eq:min_hybrid} is inherently non-linear, even within the hyperrectangle $\mathcal{D}$ where the polynomial component of the hybrid scheme acts, due to the non-linearity of the NFEs.~However, in a regression setting (where the hybrid scheme approximates a function using labeled targets), the polynomial counterpart reduces to a linear residual minimization problem.~Unlike non-linear optimization, this problem does not require iterative methods (e.g., the LM algorithm) to determine the polynomial coefficients.~Instead, it should be solved directly using linear techniques such as the Moore-Penrose pseudo-inverse, as iterative optimization introduces unnecessary numerical errors in the coefficient estimates; we provide such a polynomial regression example in \cref{app:reg_LvsNL}.
\end{remark}

\begin{remark}
    \label{rm:3}
    Building on \cref{rm:2}, one might consider linearizing the PI optimization problem in \cref{eq:min_hybrid} by locally approximating the polynomial counterpart via Taylor series expansions around equilibrium, following \cite{kazantzis2001invariant,kazantzis2005model}.~However, this would necessitate a two-step procedure for the hybrid scheme; first, solve the linearized problem within $\mathcal{D}$ to obtain the polynomial coefficients and then solve a non-linear optimization problem to train the NN for approximating the IM outside $\mathcal{D}$.~This approach has several drawbacks: (i) it requires Taylor expansions of the non-linear functions $\mathbf{f}$ and $\mathbf{g}$ in \cref{eq:genDDS2}, and as such it does not provide flexibility on the radius of the hyperrectangle $\mathcal{D}$, since the latter is constrained by the radius of convergence of the Taylor expansions, (ii) it is problem-dependent, as the linearized polynomial series varies for each dynamical system, and (iii) it becomes computationally intractable as the system dimension or polynomial degree increases.~In contrast, our method--while sacrificing some accuracy in the polynomial coefficients (see \cref{rm:2})--provides a one-step, plug-and-play solution that avoids these limitations. 
\end{remark}

\begin{remark}
\label{rm:nn_th}
For the standalone NN scheme, we establish approximation error and convergence guarantees in \cref{th:NN_NFEs}.~This result demonstrates that under Barron regularity assumptions on the true IM, the NN approximation of the NFEs converges at a rate $\mathcal{O}(L^{-1/2})$ in the $L^2$ sense.~While extending such guarantees to the hybrid scheme is nontrivial due to the piecewise architecture, the theorem provides theoretical support for the NN component used in our hybrid scheme. 
\end{remark}

\section{Illustrative examples}
\label{sec:Bench}
To demonstrate the efficiency of the proposed PI hybrid scheme, we consider two illustrative examples.~The first example focuses on an enzymatic bioreactor problem studied in \cite{kazantzis2005model}, in which an $M=1$-dim.~autonomous input drives an $N=1$-dim.~state, with the PSE of the IM known from \cite{kazantzis2005model}.~The second example considers a car-following problem with an autonomous leader vehicle \cite{bando1995dynamical}, in which the leader's $M=2$-dim.~input dynamics drives the $N=20$-dim.~state, acting as a pacemaker (e.g., a police car) in highway traffic flow.~This example serves for demonstrating the generalization of our method to approximate IMs of higher dimensions.~An additional validation example involving a dynamical system constructed to possess an exact (rather than approximate) IM is provided in \cref{app:lnEx}, offering a known analytical solution to validate our approach.~Below, we briefly describe these examples.  

\subsection{Case Study A. Enzymatic bioreactor problem}
\label{sb:Kax05ex}
We first consider the enzymatic bioreactor problem studied in \cite{kazantzis2005model}, which models a continuous stirred tank reactor where an enzyme converts a substrate into product via a ping-pong bi-mechanism.~Applying Euler discretization with a time step $\delta = 0.01$ and introducing deviations from the equilibrium $(S_0, 0)$ as $x = S - S_0$ (substrate deviation) and $y = E$ (enzyme concentration), we obtain the following non-linear discrete dynamical system in the form of \cref{eq:genDDS2}:
\begin{align}
    x(k+1) & = (1-\delta v_r)x(k) + \dfrac{\delta k_1 S_0}{1-k_2 S_0} y(k) + \dfrac{\delta k_1 y(k) x(k)}{(1 - k_2 S_0) (1 - k_2 S_0 - k_2x(k))}, \nonumber \\ 
    y(k+1) & = (1-\delta k_{d1}) y(k). \label{eq:ER_DDSeq}
\end{align}
Here, $k_1$, $k_2$ and $k_{d1}$ are kinetic parameters, $v_r$ is the ratio of the flow rate of substrate to the reactor volume and $S_0$ is the substrate concentration in the feed stream.~Following \cite{kazantzis2005model}, we adopt the parameter values $k_1=0.082$, $k_2=0.59$, $k_{d1}=0.0034$, $v_r=2$ and $S_0=3.4$.

The system in \cref{eq:ER_DDSeq} satisfies all the assumptions regarding the non-linear functions at the equilibrium $(x_0,y_0)=(0,0)$.~Additionally, since $\mathbf{A}=1-\delta k_{d1}$ and $\mathbf{B}=1-\delta v_r$, the assumptions of \cref{th:th1Kaz} are satisfied, as the non-resonance condition $(1-\delta k_{d1})^d \neq 1-\delta v_r$ holds for any $d\in\mathbb{N}$.~This implies the existence of an analytic $N=1$-dim. IM $\mathcal{M}$ in the form of \cref{eq:IM}, the mapping of which $x=\pi(y)$ can be approximated by solving the NFEs in \cref{eq:NFEs_IM}.~Following \cref{app:IM_PS}, we obtained the PSE approximation of the IM functional $\tilde{x}=\tilde{\pi}^{PSE}(y)$ in the form of \cref{eq:IM_PSexp}.~The analytically determined coefficients $\pi^1,\ldots,\pi^{1,\ldots,1}$ are in excellent agreement with those reported in \cite{kazantzis2005model} for a Taylor expansion of degree $h=10$.

\subsection{Case Study B. Car-following problem with an autonomous leader}
\label{sb:CFM}
We consider a platoon of $N_c$ vehicles moving along a single lane of highway, where each vehicle follows the one directly ahead.~The leading vehicle is an autonomous pacemaker with prescribed dynamics that is unaffected by the vehicles behind it.~This pacemaker can be interpreted as a police car or a connected autonomous vehicle (CAV) regulating traffic flow, or as a leader who observes stop-and-go wave patterns upstream.

To model this problem, we employ the optimal velocity model (OVM) \cite{bando1995dynamical} for the $N_c$ followers, according to which each vehicle adjusts its acceleration based on the relative velocity and spacing with respect to the vehicle ahead.~Following standard practice, we describe the follower dynamics in terms of headway $h_i(k)$ and velocity $v_i(k)$ for each follower $i=1,\ldots,N_c$, where $h_i(k)$ denotes the distance to the vehicle ahead.~Applying Euler discretization with a time step $\delta=0.05$, the discrete system for the followers is:
\begin{equation}
    h_i(k+1) = h_i(k) + \delta \bigl( v_{i+1}(k) - v_i(k) \bigr), \qquad
    v_i(k+1)= v_i(k) + \delta \tau^{-1} \bigl( V(h_i(k)) - v_i(k) \bigr), \label{eq:CFM_fol1}
\end{equation}
where $\tau$ is the adaptation time, and $V(h)$ is the optimal velocity function:
\begin{equation}
    V(h) = v_0 \dfrac{\tanh(\gamma h - \beta) + \tanh(\beta)}{1 + \tanh(\beta)},
    \label{eq:OVF}
\end{equation}
where $v_0$ is the desired speed (maximum speed assuming no vehicle upstream), $\gamma$ is the inverse of transition width, and $\beta$ is a form factor parameter \cite{bando1995dynamical}.~Note that the $N_c$-th follower in \cref{eq:CFM_fol1} depends on the velocity of the autonomous leader $v_{N_c+1}(k)\equiv v_{\ell}(k)$. 

The autonomous leader follows its own prescribed dynamics independent of the followers, serving as a driver input to the system of followers.~To model a pacemaker that induces oscillatory behavior (e.g., representing stop-and-go waves), we consider a damped oscillator targeting a desired velocity $v_{des}$.~Introducing an auxiliary variable $z(k)$, the discrete dynamics of the leader is:
\begin{equation}
z(k+1) = z(k) + \delta \left( v_{\ell}(k) - v_{des} \right), \qquad 
v_{\ell}(k+1) = v_{\ell}(k) - \delta \left( \tau_{\ell}^{-1} (v_{\ell}(k) - v_{des}) + \mu z(k) \right), \label{eq:CFM_lead1}
\end{equation}
where $\tau_{\ell}$ is the adaptation time of the leader and $\mu$ controls the damping.

From \cref{eq:CFM_fol1,eq:CFM_lead1}, the equilibrium state is given by $h_{i,0}=V^{-1}(v_{des})$, $v_{i,0}=v_{des}$, $z_0=0$ and $v_{\ell,0}=v_{des}$, where $V^{-1}$ is the inverse of the optimal velocity function in \cref{eq:OVF}.~Collecting the  followers' state $\hat{\mathbf{x}}(k)=[h_1(k), \ldots, h_{N_c}(k), v_1(k), \ldots, v_{N_c}(k)]^\top \in \mathbb{R}^{2N_c}$ and the leader's state $\hat{\mathbf{y}}(k) = [z(k),v_{\ell}(k)]^\top \in \mathbb{R}^2$, and introducing the deviations from equilibrium $\mathbf{x}=\hat{\mathbf{x}}-[V^{-1}(v_{des}), \ldots, V^{-1}(v_{des}), v_{des}, \ldots, v_{des}]^\top$ and $\mathbf{y}=\hat{\mathbf{y}}-[0, v_{des}]^\top$, we obtain a nonlinear discrete dynamical system in the form of \cref{eq:genDDS2} with $N = 2N_c$ and $M = 2$; the analytic expressions are omitted here for brevity.

In our simulations, we considered $N_c=10$ followers.~Following \cite{treiber2012traffic}, for the followers we adopt typical highway parameter values: $\tau=0.65$, $\gamma=1/15$, $\beta=1.5$ and $v_0=33.3$ (corresponding to $120$ km/h).~For the leader, we set a much higher adaptation time $\tau_{\ell} = 10$, $\mu=1$ and the desired velocity $v_{des}=v_0/2$ (corresponding to $60$ km/h); thus, the leader slowly regulates the platoon's velocity with a damped oscillation.~Under this parameter setup, the model satisfies all the assumptions of \cref{th:th1Kaz}, implying the existence of an analytic $N=20$-dim. IM $\mathcal{M}$ in the form of \cref{eq:IM}, whose mapping $\mathbf{x}=\boldsymbol{\pi}(\mathbf{y})$ can be approximated by solving the NFEs in \cref{eq:NFEs_IM}.~Following \cref{app:IM_PS}, we obtained the PSE approximation of the IM functional $\tilde{\mathbf{x}}=\tilde{\boldsymbol{\pi}}^{PSE}(\mathbf{y})$ in the form of \cref{eq:IM_PSexp} using a Taylor expansion of degree $h=5$.

\section{Numerical Results}
\label{sec:NumR}
This section evaluates the numerical accuracy of the proposed PI hybrid scheme in learning IM approximations for the two examples described in \cref{sec:Bench}.~We compare the IM approximations provided by the PI hybrid schemes with those obtained by using standalone PINNs and the PSE expansions derived by the symbolic approach proposed in \cite{kazantzis2001invariant,kazantzis2005model}.~Training and evaluation of the PI schemes followed the procedures outlined in \cref{sbsb:train,sbsb:numAc}.~All simulations were carried out with a CPU Intel(R) Core(TM) CPU i7-13700H @ 2.40GHz, RAM 32.0 GB using MATLAB R2024b.

\subsection{Case Study A. Enzymatic bioreactor problem}
Here, we approximate the IM in the domain $\Omega = [0,4]$, where the system in \cref{eq:ER_DDSeq} is biologically plausible (i.e., for non-negative enzyme concentrations $y \ge 0$).~For the proposed PI hybrid and standalone PINN schemes, we considered polynomials of degree $h = 10$, and NNs with $L = 10$ neurons in the hidden layer.~Following the procedure detailed in \cref{sbsb:train}, we uniformly sampled $Q = 620$ collocation points along $\Omega$.~All hyperparameters are set as described in \cref{app:Imp}.~To explore the role of the radius of convergence $r$, we trained the PI hybrid scheme with power series (PI-HPS) for $r = 0.5, 1, 2, 4$.~Note that choosing $r = 4$ makes the hyperrectangle $\mathcal{D}$ cover $\Omega$, and thus the hybrid scheme is equivalent to considering the polynomial series over the entire domain. We also trained the PI hybrid schemes with Legendre and Chebyshev polynomials (PI-HLS, PI-HCS) with $r = 1$, and the standalone PINNs.~Convergence results for the seven PI approximations are provided in \cref{tb:conv_Kaz05ex} of Supplement \ref{supp1}, obtained over 100 training realizations with different initial parameters; see \cref{sbsb:train} for details.  


To evaluate the numerical approximation accuracy of the learned IM approximations, we constructed a testing data set from $n_{IC}=10$ system trajectories, removing the first $k_{trans}=8$ steps as transients to ensure the data lie on the manifold.~Following the procedure in \cref{sbsb:numAc}, we collected $S=10,000$ data points $(x^{(s)},y^{(s)})$ in $\Omega=[0,4]$ by selecting initial conditions with random $x(0)\in[-1,1]$ and fixed $y(0)=4.3$.~From this testing set, we computed the relative errors $\lVert x^{(s)}-\tilde{\pi}(y^{(s)})\rVert/ \lVert x^{(s)}\rVert$ for the IM approximations $\tilde{\pi}(y)$ of the trained PI schemes.~The means and 5–95\% percentiles of the $L^1$, $L^2$, and $L^\infty$ norms of these errors, calculated over the 100 parameter sets obtained during training, are reported in \cref{tb:acc_Kaz05ex}.~For comparison, the same metrics are included for the power series expansion (PSE) of degree $h=10$ derived in \cite{kazantzis2005model}.
\begin{table}[!h]
\centering
\caption{Numerical approximation accuracy of IM approximations $\tilde{\pi}(y)$ for the enzymatic bioreactor problem in \cref{sb:Kax05ex} on the test data.~Relative errors ($L^1$, $L^2$, and $L^\infty$ norms) are reported for the PSE with degree $h=10$ from \cite{kazantzis2005model}, the standalone PINN and the PI schemes (PI-HxS with $x=P,L,C$ denoting power, Legendre, and Chebyshev polynomial series); for hyperparameters see \cref{tb:conv_Kaz05ex}.~Mean values and 5–95\% percentiles are computed over the 100 parameter sets obtained during training of the PI schemes.~Testing data errors are evaluated over $S=10,000$ points.}
\resizebox{\textwidth}{!}{
\begin{tabular}{l c c c c c c}
\textbf{Scheme} & \multicolumn{2}{c}{$\lVert x^{(s)}-\tilde{\pi}(y^{(s)})\rVert_1/\lVert x^{(s)}\rVert_1$} & \multicolumn{2}{c}{$\lVert x^{(s)}-\tilde{\pi}(y^{(s)})\rVert_2/ \lVert x^{(s)}\rVert_2$} & \multicolumn{2}{c}{$\lVert x^{(s)}-\tilde{\pi}(y^{(s)})\rVert_\infty/ \lVert x^{(s)}\rVert_\infty$} \\
& mean & 5--95\% percentiles & mean & 5--95\% percentiles & mean & 5--95\% percentiles \\
\noalign{\smallskip}\hline\noalign{\smallskip}
PSE	&	$3.15\times 10^{-3}$	&				&	$1.11\times 10^{-2}$	&				&	$3.85\times 10^{-2}$ &				\\
PINN	        &	$2.04\times 10^{-4}$	& [$3.11\times 10^{-5}-4.62\times 10^{-4}$] &	$2.58\times 10^{-4}$	& [$6.39\times 10^{-5}-3.62\times 10^{-4}$] &	$1.47\times 10^{-3}$	& [$4.80\times 10^{-4}-2.04\times 10^{-3}$] \\	
PI-HPS, $r=0.5$	&	$6.55\times 10^{-5}$	& [$1.64\times 10^{-5}-1.44\times 10^{-4}$] &	$1.76\times 10^{-4}$	& [$6.05\times 10^{-5}-3.46\times 10^{-4}$] &	$1.14\times 10^{-3}$	& [$4.62\times 10^{-4}-2.08\times 10^{-3}$] \\
PI-HPS, $r=1$	&	$6.13\times 10^{-5}$	& [$1.85\times 10^{-5}-1.19\times 10^{-4}$] &	$1.87\times 10^{-4}$	& [$7.14\times 10^{-5}-3.26\times 10^{-4}$] &	$1.21\times 10^{-3}$	& [$5.31\times 10^{-4}-1.99\times 10^{-3}$] \\
PI-HPS, $r=2$	&	$4.12\times 10^{-5}$	& [$1.46\times 10^{-5}-5.99\times 10^{-5}$] &	$1.56\times 10^{-4}$	& [$6.62\times 10^{-5}-2.13\times 10^{-4}$] &	$1.04\times 10^{-3}$	& [$5.01\times 10^{-4}-1.37\times 10^{-3}$] \\
PI-HPS, $r=4$	&	$1.66\times 10^{-4}$	& [$1.66\times 10^{-4}-1.66\times 10^{-4}$] &	$1.48\times 10^{-4}$	& [$1.48\times 10^{-4}-1.48\times 10^{-4}$] &	$8.55\times 10^{-4}$	& [$8.55\times 10^{-4}-8.56\times 10^{-4}$] \\
PI-HLS, $r=1$	&	$6.22\times 10^{-5}$	& [$1.52\times 10^{-5}-1.15\times 10^{-4}$] &	$1.84\times 10^{-4}$	& [$5.78\times 10^{-5}-3.18\times 10^{-4}$] &	$1.19\times 10^{-3}$	& [$4.47\times 10^{-4}-1.95\times 10^{-3}$] \\
PI-HCS, $r=1$	&	$7.15\times 10^{-5}$	& [$1.72\times 10^{-5}-1.45\times 10^{-4}$] &	$1.97\times 10^{-4}$	& [$6.75\times 10^{-5}-3.46\times 10^{-4}$] &	$1.27\times 10^{-3}$	& [$5.06\times 10^{-4}-2.10\times 10^{-3}$] \\
\noalign{\smallskip}\hline
\end{tabular}
}
\label{tb:acc_Kaz05ex}
\end{table}
\Cref{tb:acc_Kaz05ex} shows that the PSE approximation, derived via expansion around the equilibrium, provides modest accuracy.~On the other hand, the PI hybrid schemes achieve much higher approximation accuracy, with mean $L^2$ relative errors on the order of $10^{-5}$.~As $r$ increases, accuracy improves slightly across the testing set.~Overall, the PI hybrid schemes achieve slightly higher global accuracy than the PINN scheme, and both outperform the local PSE approximation by one to two orders of magnitude.

\begin{figure}[!h]
\centering
\centering
    \subfigure[Degree $h=10$]{\includegraphics[width=0.47\textwidth]{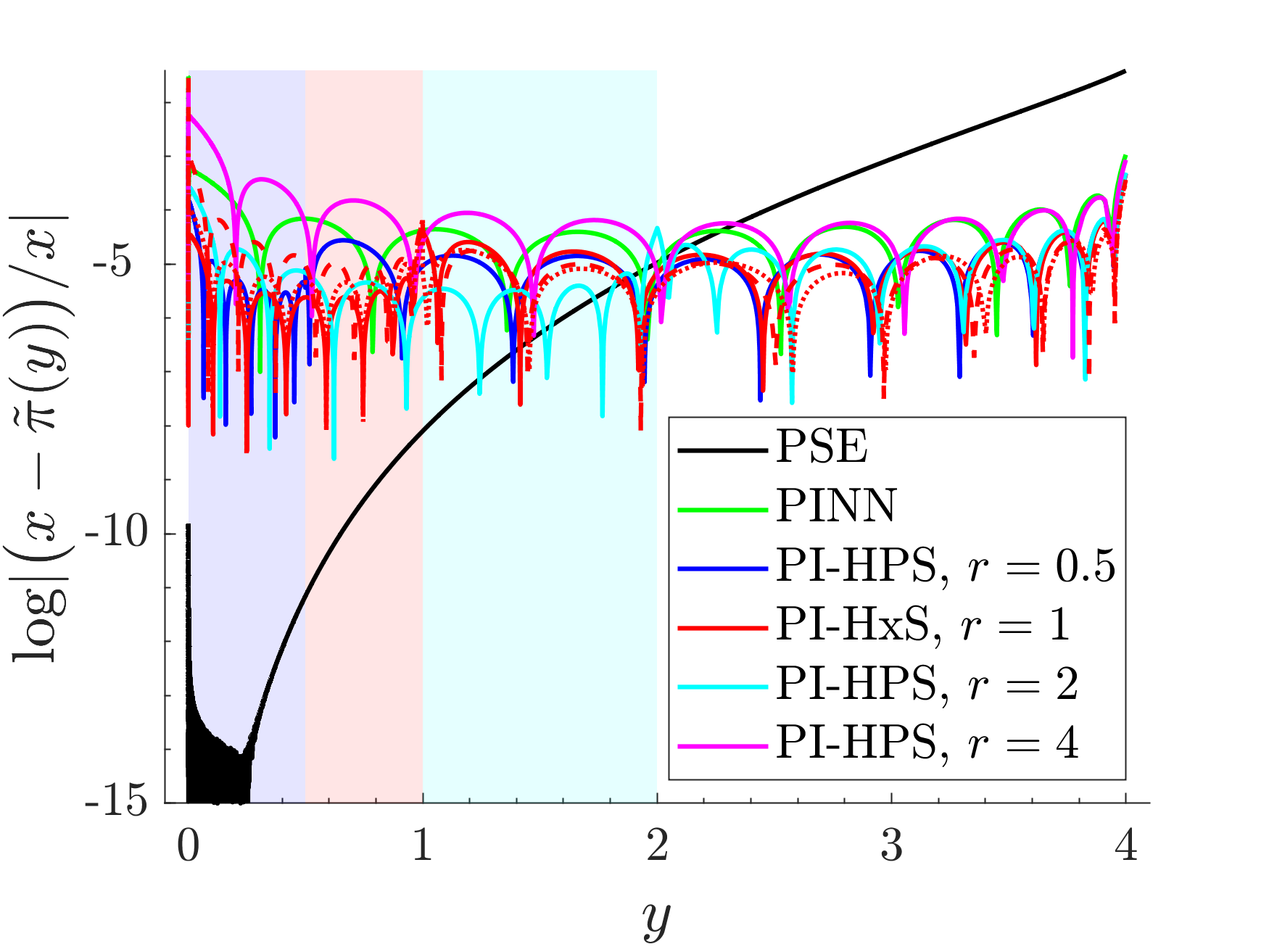}} \hspace{10pt}
    \subfigure[Degree $h=5$]{\includegraphics[width=0.47\textwidth]{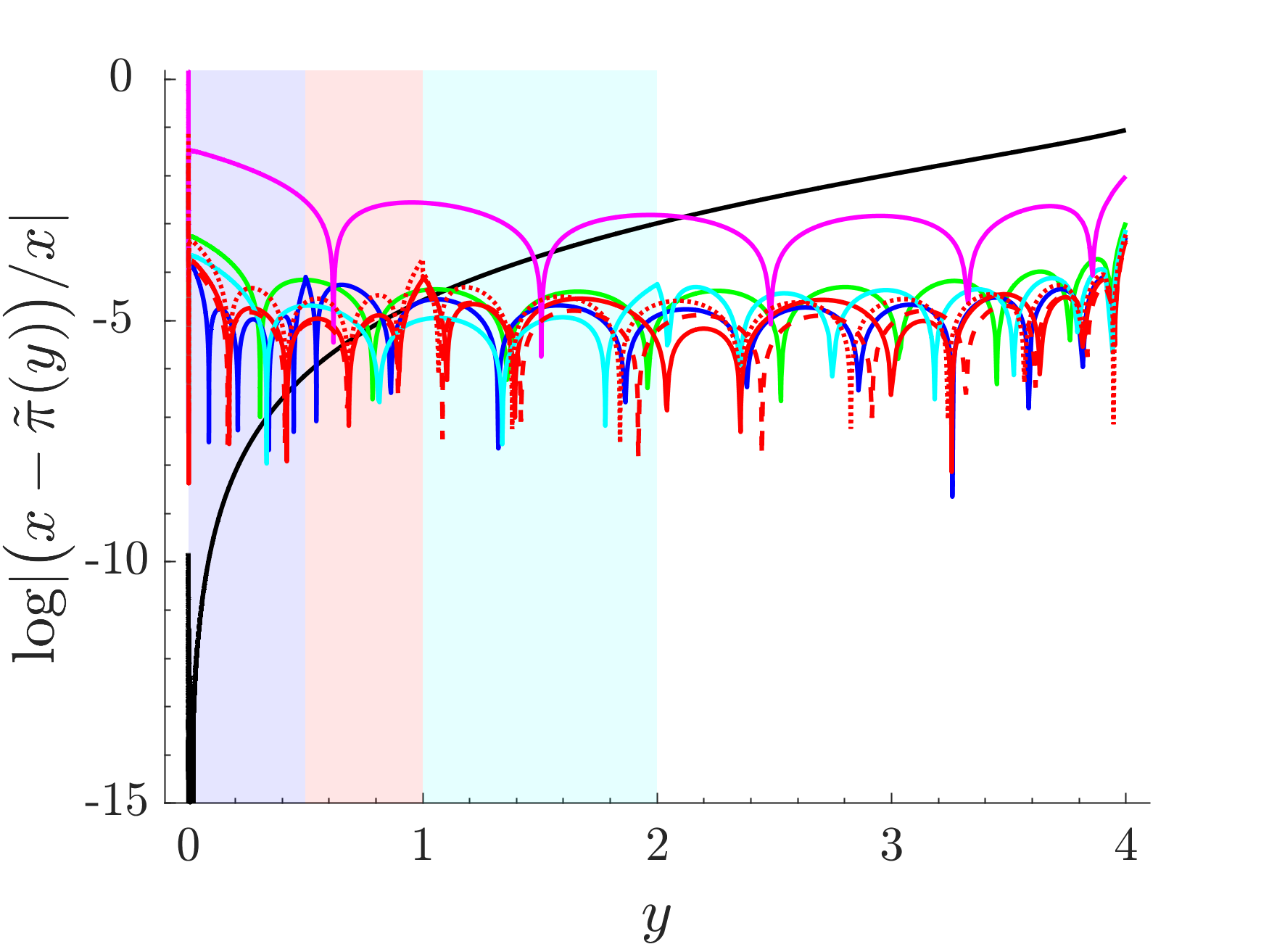}}
\caption{Relative errors of IM approximations $\tilde{\pi}(y)$ for the enzymatic bioreactor problem in \cref{sb:Kax05ex} on the test data lying on the IM.~The IM approximations include the PSE in \cite{kazantzis2005model} (black), the PINN (green) and the PI hybrid schemes with different radii $r$ (blue, red, cyan, magenta).~Power, Legendre, and Chebyshev polynomial series ($x=P,L,C$) are distinguished by solid, dashed, and dotted curves, respectively.~Background colors (red, blue, cyan) correspond to the hyperrectangles $\mathcal{D}$ of each hybrid scheme.~Panels (a) and (b) correspond to polynomial degrees $h=10$ and $h=5$, respectively, and $L=10$ neurons are used in the NNs.}
\label{fig:acc_Kaz05ex}
\end{figure}

As the above error metrics are global on the testing set, we examine local accuracy by depicting the point-wise relative errors $\lvert x^{(s)}-\tilde{\pi}(y^{(s)})\rvert/ \lvert x^{(s)} \rvert$ in \cref{fig:acc_Kaz05ex}, 
where the PSE approximation derived in \cite{kazantzis2005model}, the PINN, and the PI hybrid scheme approximations for different radii $r$ are compared.~\Cref{fig:acc_Kaz05ex}a) demonstrates that the PSE provides high accuracy near the equilibrium but exhibits large errors farther away, as reflected in the global accuracy metrics in \cref{tb:acc_Kaz05ex}.~On the other hand, the PINN scheme generates a homogeneously distributed accuracy profile across the entire domain $\Omega$ on the order of $10^{-5}$.~The PI hybrid scheme combines the regions of high accuracy of the above approximations.~In particular, as indicated by the highlighted backgrounds in \cref{fig:acc_Kaz05ex}a), the hybrid scheme delivers one to two orders of magnitude higher accuracy than the PINN within $[0,r]$ (where the polynomial series are active) and comparable (or even higher) accuracy for $y>r$.~
In addition, the Legendre and Chebyshev polynomial series provide lower accuracy near the equilibrium compared to the PSE.~Very similar trends are observed for $h=5$ in \cref{fig:acc_Kaz05ex}b), with the PSE and the PI hybrid scheme with $r=4$ being less accurate; all other differences compared to the $h=10$ schemes are negligible.

\subsection{Case Study B. Car-following problem with an autonomous leader}
In this case study, we approximate the IM in the domain $\Omega=[-5,5]\times[-5,5]$ via the proposed PI hybrid and standalone PINN schemes, using polynomials of degree $h=3$ and NNs of $L=20$ neurons in the hidden layer.~Following the procedure detailed in \cref{sbsb:train}, we sampled $Q=1620$ collocation points along $\Omega$ by sampling data from $n_{IC}=200$ system trajectories, after removing $k_{trans}=800$ steps as transients.~For each initial condition, the followers' headways were randomly sampled from $h_i(0)\in\mathcal{U}[35,45]$ with velocities computed via the optimal velocity function in \cref{eq:OVF}, while the leader was placed $50$ m ahead of the $N_c$-th follower with a random $z(0)\in\mathcal{U}[-10,10]$ and velocity $v_{\ell}=27.7$ (corresponding to $100$ km/h).~From these $Q$ collocation points, we trained the PI hybrid schemes with power, Legendre, and Chebyshev polynomials (PI-HPS, PI-HLS, PI-HCS), setting the radius to $r=1$, and the standalone PINN scheme. Convergence results for these four PI approximations are provided in \cref{tb:conv_CFM} of \cref{supp1}, obtained over 100 training realizations with different initial parameters; see \cref{sbsb:train} for details.

To evaluate the numerical approximation accuracy of the learned IM approximations, we constructed a testing data set following the same procedure used for the collocation points, this time sampling $S=10,000$ data points $(\mathbf{x}^{(s)},\mathbf{y}^{(s)})$ randomly in $\Omega$; see \cref{sbsb:numAc}.~From this testing set, we computed, for each component $n=1,\ldots,N$, the $L^1$, $L^2$, and $L^\infty$ norms of the relative errors $\lVert x_n^{(s)}-\tilde{\pi}_n(\mathbf{y}^{(s)})\rVert/ \lVert x_n^{(s)}\rVert$ over all data points.~These norms were then averaged across the $N=20$ components, yielding a single $L^1$, $L^2$, and $L^\infty$ error per training realization.~The means and 5-95\% percentiles of these errors over the 100 parameter sets obtained during training are reported in \cref{tb:acc_CFM}.~For comparison, the same metrics are included for the PSE of degree $h=3$.
\begin{table}[!h]
\centering
\caption{Numerical approximation accuracy of IM approximations $\tilde{\boldsymbol{\pi}}(\mathbf{y})$ for the car-following problem with an autonomous leader in \cref{sb:CFM} on the test data.~Errors ($L^1$, $L^2$, and $L^\infty$ norms) are reported for the PSE with degree $h=3$, the standalone PINN and the PI schemes (PI-HxS with $x=P,L,C$ denoting power, Legendre, and Chebyshev polynomial series); for hyperparameters see \cref{tb:conv_CFM}.~Mean values and 5–95\% percentiles are computed over the 100 parameter sets obtained during training of the PI schemes.~Testing data errors are evaluated over $S=10,000$ points.}
\resizebox{\textwidth}{!}{
\begin{tabular}{l c c c c c c}
\textbf{Scheme} & \multicolumn{2}{c}{$\text{mean}_{n=1,\ldots,N} \lVert x_n^{(s)}-\tilde{\pi}_n(\mathbf{y}^{(s)})\rVert_1/ \lVert x_n^{(s)}\rVert_1$} & \multicolumn{2}{c}{$\text{mean}_{n=1,\ldots,N} \lVert x_n^{(s)}-\tilde{\pi}_n(\mathbf{y}^{(s)})\rVert_2/ \lVert x_n^{(s)}\rVert_2$} & \multicolumn{2}{c}{$\text{mean}_{n=1,\ldots,N}  \lVert x_n^{(s)}-\tilde{\pi}_n(\mathbf{y}^{(s)})\rVert_\infty/\lVert x_n^{(s)}\rVert_\infty$} \\
& mean & 5--95\% percentiles & mean & 5--95\% percentiles & mean & 5--95\% percentiles \\
\noalign{\smallskip}\hline\noalign{\smallskip}
PSE	&	$7.98\times 10^{-2}$	&				&	$1.81\times 10^{-1}$	&				&	$5.96\times 10^{-1}$	&				\\
PINN	   &	$4.18\times 10^{-3}$	& [$2.70\times 10^{-3}-5.45\times 10^{-3}$] &	$6.30\times 10^{-3}$	& [$3.41\times 10^{-3}-8.30\times 10^{-3}$] &	$2.81\times 10^{-2}$	& [$1.52\times 10^{-2}-3.60\times 10^{-2}$] \\
PI-HPS, $r=1$	&	$2.65\times 10^{-3}$	& [$1.80\times 10^{-3}-3.82\times 10^{-3}$] &	$4.38\times 10^{-3}$	& [$2.88\times 10^{-3}-6.45\times 10^{-3}$] &	$2.05\times 10^{-2}$	& [$1.35\times 10^{-2}-3.01\times 10^{-2}$] \\
PI-HLS, $r=1$	&	$3.03\times 10^{-3}$	& [$2.10\times 10^{-3}-4.22\times 10^{-3}$] &	$5.01\times 10^{-3}$	& [$2.91\times 10^{-3}-7.78\times 10^{-3}$] &	$2.30\times 10^{-2}$	& [$1.27\times 10^{-2}-3.52\times 10^{-2}$] \\
PI-HCS, $r=1$	&	$2.98\times 10^{-3}$	& [$1.68\times 10^{-3}-4.14\times 10^{-3}$] &	$4.92\times 10^{-3}$	& [$2.68\times 10^{-3}-6.85\times 10^{-3}$] &	$2.28\times 10^{-2}$	& [$1.21\times 10^{-2}-3.09\times 10^{-2}$] \\
\noalign{\smallskip}\hline
\end{tabular}
}
\label{tb:acc_CFM}
\end{table}
\Cref{tb:acc_CFM} shows that the PSE approximation, derived via expansion around the equilibrium, provides poor accuracy.~On the other hand, the PI hybrid schemes achieve higher approximation accuracy, with mean $L^2$ relative errors across all $N=20$ components and over all testing points on the order of $10^{-3}$.~Importantly, the PI hybrid schemes achieve slightly higher global accuracy than the PINN scheme, both outperforming the local PSE approximation.

\begin{figure}[!h]
\centering
    \subfigure[$x_1$ error, PSE]{\includegraphics[width=0.32\textwidth]{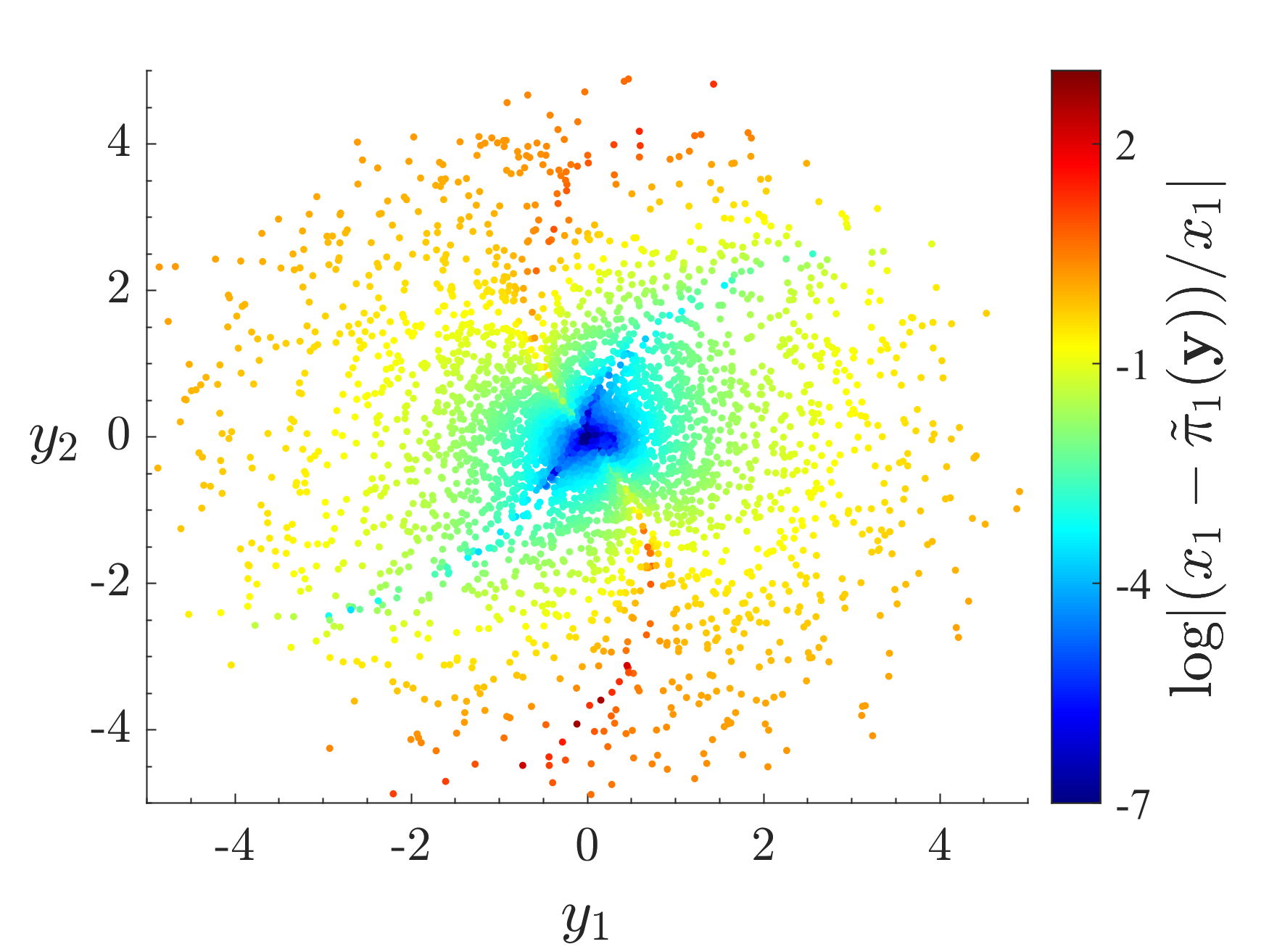}} \hspace{1pt}
    \subfigure[$x_1$ error, PINN]{\includegraphics[width=0.32\textwidth]{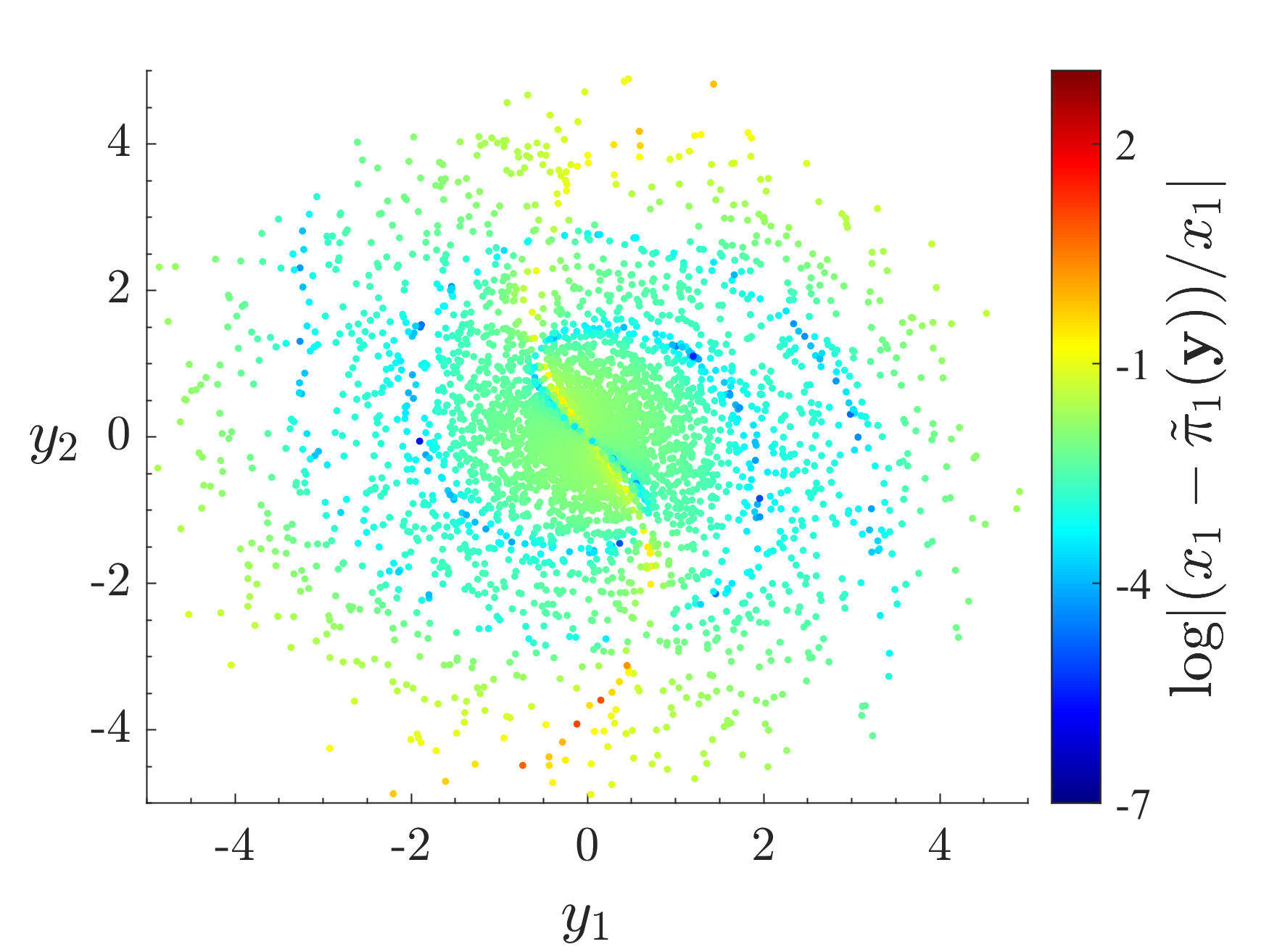}} \hspace{1pt} 
    \subfigure[$x_1$ error, PI-HPS]{\includegraphics[width=0.32\textwidth]{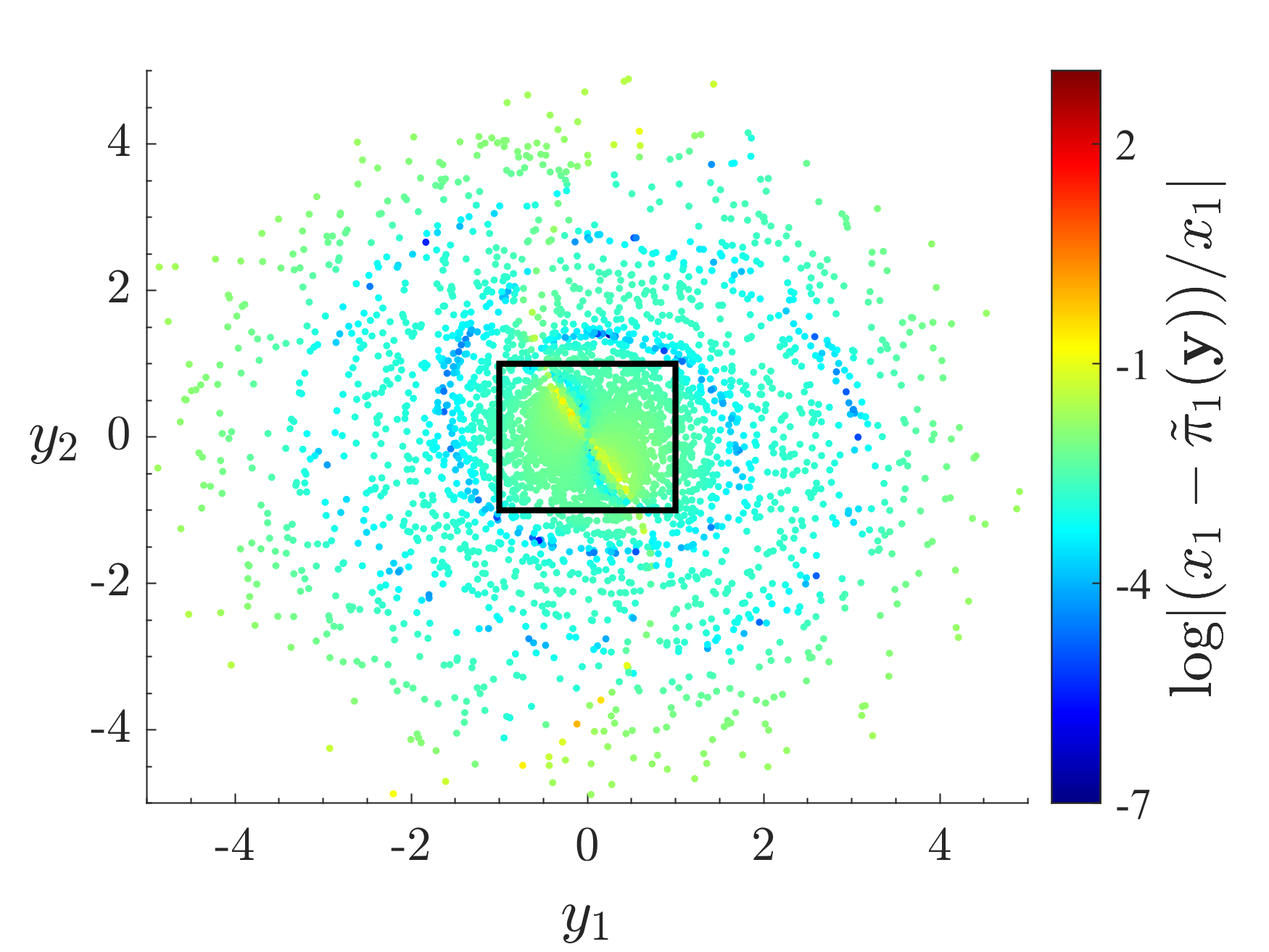}} \\
    \subfigure[$x_5$ error, PSE]{\includegraphics[width=0.32\textwidth]{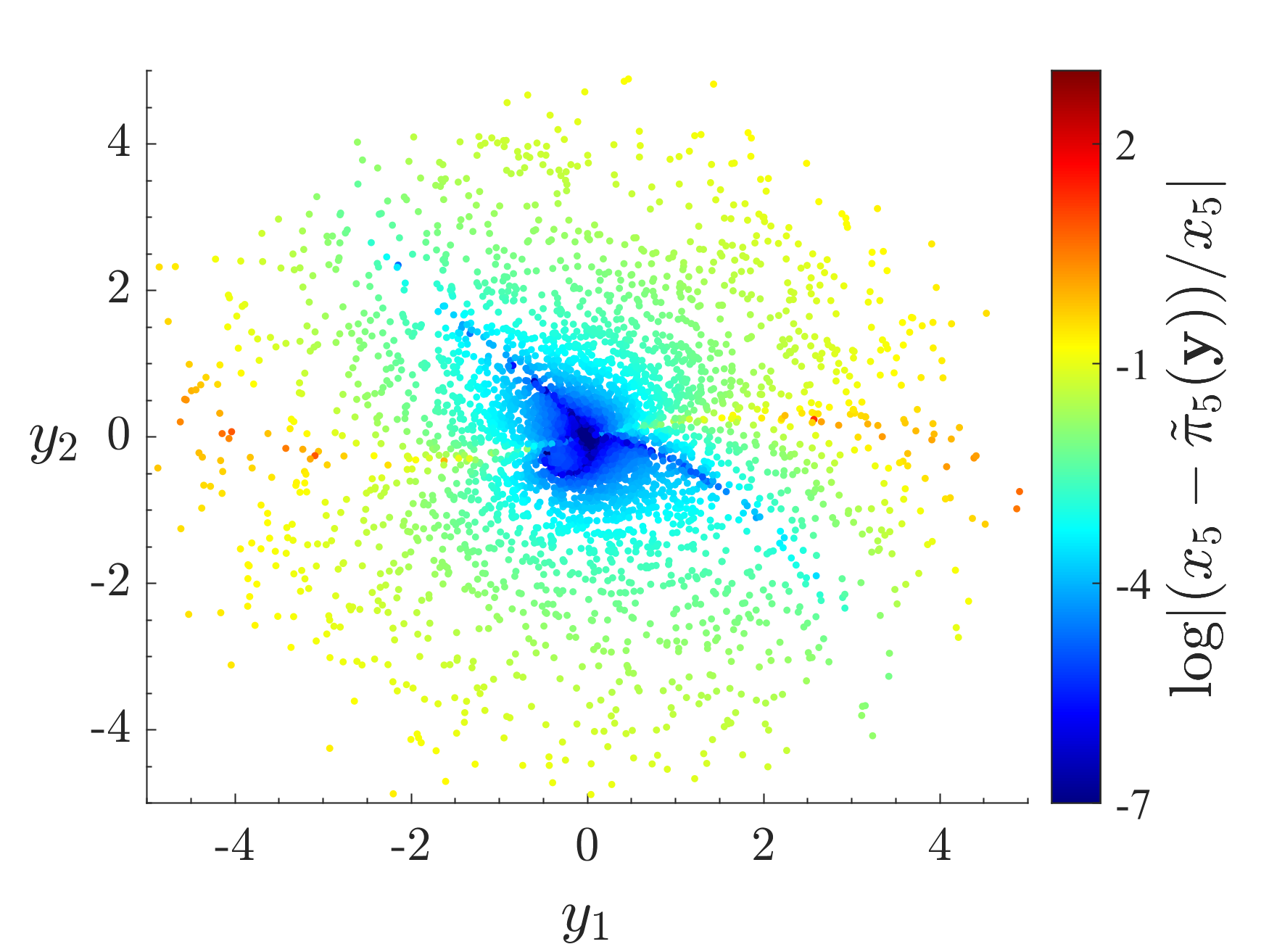}} \hspace{1pt}
    \subfigure[$x_5$ error, PINN]{\includegraphics[width=0.32\textwidth]{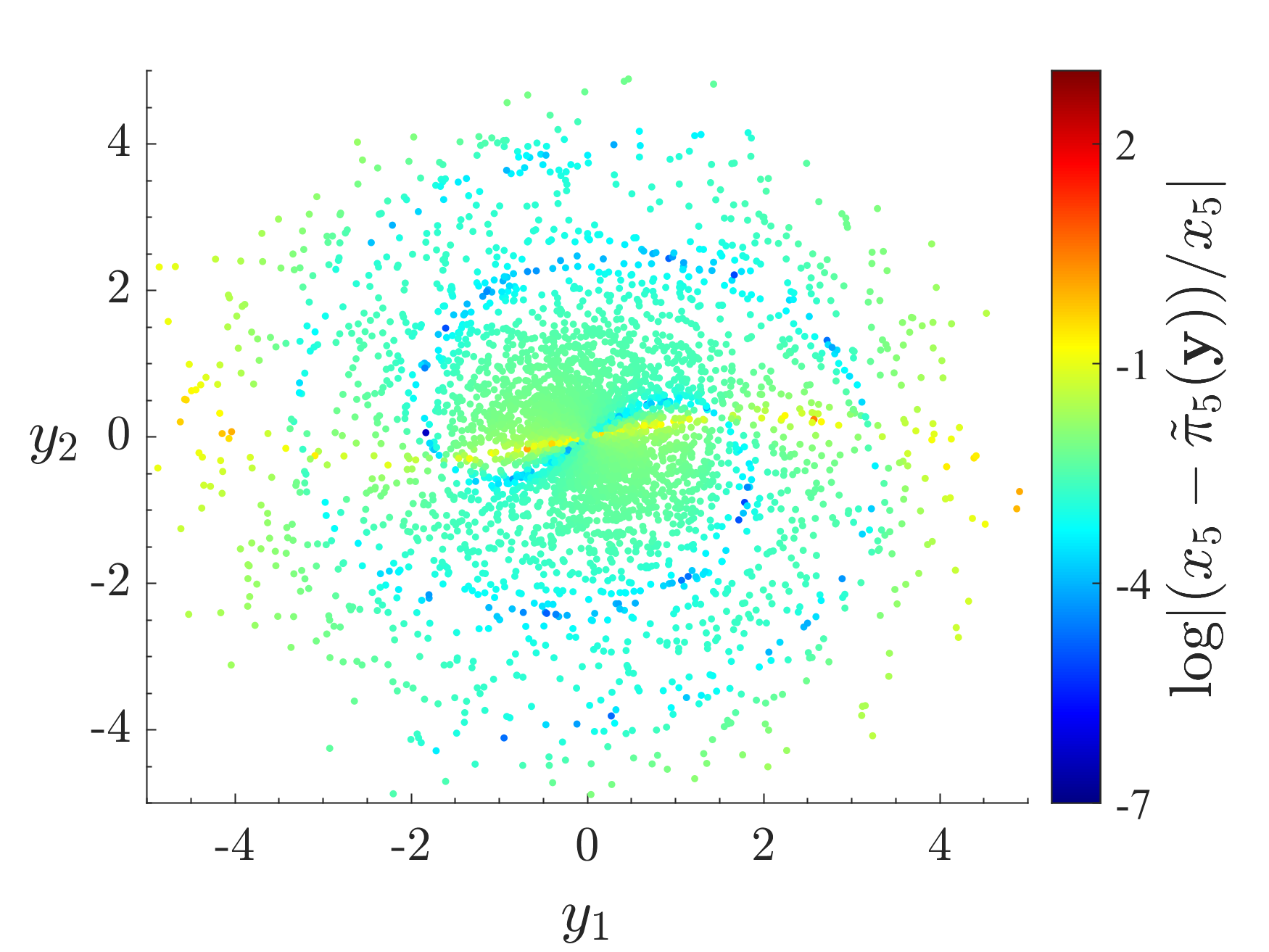}} \hspace{1pt} 
    \subfigure[$x_5$ error, PI-HPS]{\includegraphics[width=0.32\textwidth]{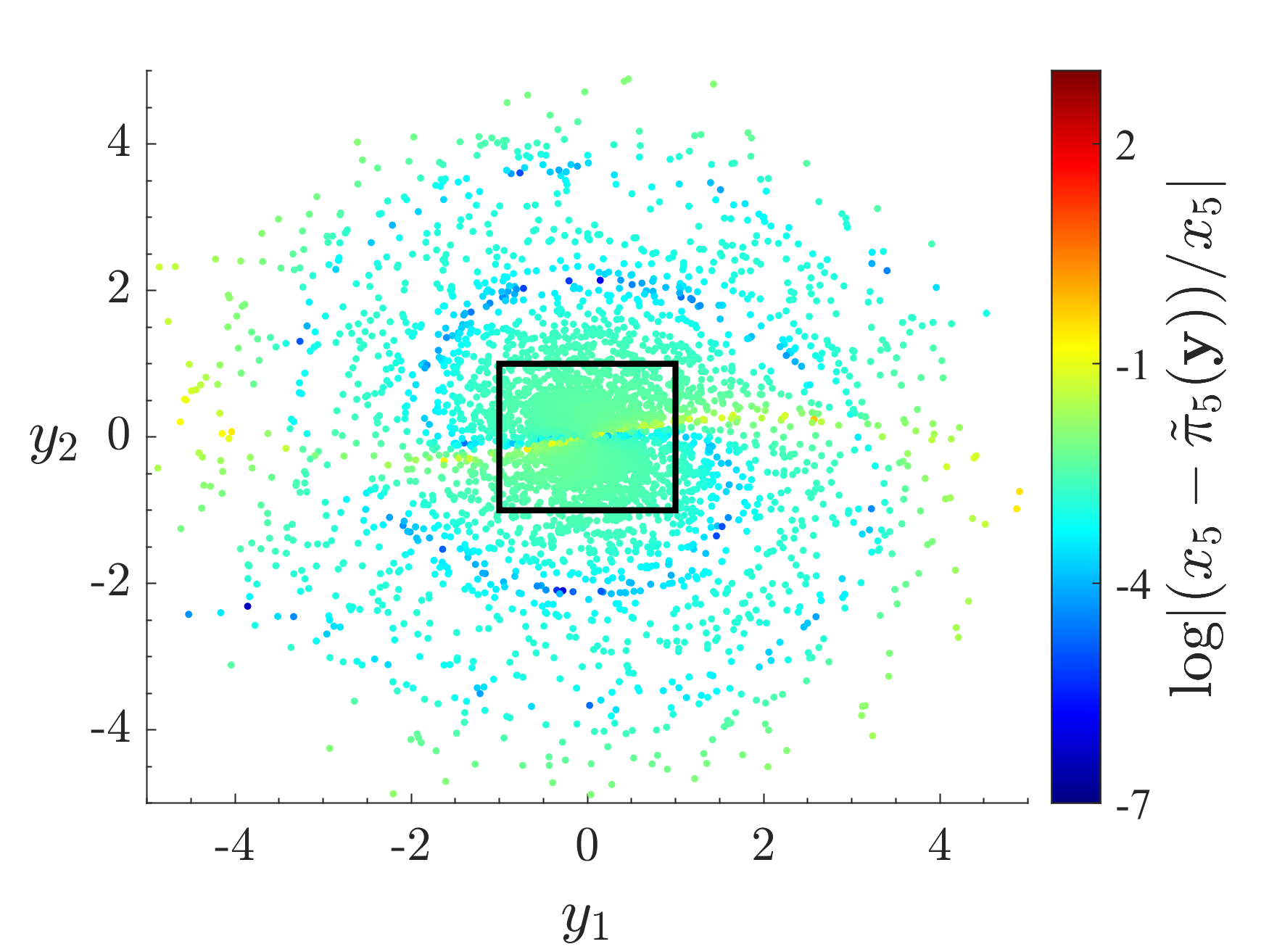}} 
\caption{Relative errors of IM approximations $\tilde{\boldsymbol{\pi}}(\mathbf{y})$ for the car-following problem with an autonomous leader in \Cref{sb:CFM}, on the test data lying on the IM.~Top/bottom row panels show the errors of the $x_1$/$x_5$ component.~Panels show (a,d) the PSE, (b,e) the PINN approximations, and (c,f) the PI hybrid scheme approximations with power series.~The black rectangle in panels (c,f) denote the hyperrectangles $\mathcal{D}$ of the PI hybrid schemes.~Polynomial degree is $h=3$, and $L=20$ neurons are used in the NNs.~Full figures are provided in \cref{fig:acc_CFM_x1_full,fig:acc_CFM_x5_full} of Supplement~\ref{supp2}.}
\label{fig:acc_CFM}
\end{figure}

To appreciate the local accuracy provided by the physics-informed ML schemes, we depict the point-wise relative errors $\lvert \mathbf{x}^{(s)} - \tilde{\boldsymbol{\pi}}(\mathbf{y}^{(s)})\rvert/\lvert \mathbf{x}^{(s)} \rvert$ for two components ($x_1$ and $x_5$) across the testing set in \cref{fig:acc_CFM}, which compares the PSE, the PINN, and the hybrid physics-informed (PI) scheme with power series; a detailed comparison of all IM approximations is provided in \cref{fig:acc_CFM_x1_full,fig:acc_CFM_x5_full} of Supplement~\ref{supp2}.
\Cref{fig:acc_Kaz05ex}a,d) demonstrate that the PSE provides high accuracy near the equilibrium but exhibits large errors farther away, as reflected in the global accuracy metrics in \cref{tb:acc_CFM}.~In contrast, both the PINN and the PI hybrid scheme produce a homogeneous error distribution across the entire domain, with no systematic degradation far from the equilibrium; see \cref{fig:acc_CFM}b–e).~The relative errors for the PINN and the PI hybrid schemes are comparable, with the hybrid schemes providing slightly lower errors near the equilibrium; the difference is less pronounced here in comparison to case study A in \cref{fig:acc_Kaz05ex}.~Nonetheless, both PI approximations of the IM outperform the PSE, despite $x_1$ exhibiting larger relative errors than $x_5$; a difference in magnitudes that is reflected in the global metrics of \cref{tb:acc_CFM}.

\section{Conclusions}
\label{sec:Con}
We present a hybrid physics-informed machine learning scheme for learning low dimensional (smooth) invariant manifolds (IM) of nonlinear discrete-time dynamical systems with nonlinear exosystems. For relatively low-dimensional IM, the proposed scheme exploits the interplay and synergistic advantages of polynomial series approximations of solutions to the associated invariance equations representing the IM-map around a fixed point, and those of NNs. Indeed, polynomial series approximations and NNs display their own relative pros and cons when used to approximate nonlinear multivariate functions. Polynomials—via Taylor or Chebyshev expansions—offer geometric (exponential) convergence rates around a point, delivering extremely rapid error decay as the degree or truncation order increases; however, as discussed in the introduction, their effectiveness is limited by the analyticity domain of the target invariant manifold. Also, while univariate Chebyshev polynomials are well studied, for their multivariate extensions (e.g.\ tensor-product Chebyshev bases) the number of coefficients explodes combinatorially with the number of variables, while generalizing their alternance and convergence properties for larger numbers of variables and higher polynomial degrees is not an easy task either \cite{mason1980near}. In practice, one must truncate the multivariate series expansion by selecting a subset of multi-indices, yet deciding which terms to keep is not trivial; thus, neglecting high-order cross terms ad-hoc can significantly degrade accuracy. Furthermore, large multivariate expansions are numerically sensitive: small errors in the estimation of the coefficients, especially of the higher-degree terms, can lead to large deviations and numerical instabilities in the reconstructed function, especially at the boundaries of the domain of convergence. On the other hand, neural networks rovide uniform approximations over broader domains where a polynomial series might diverge, and they remain both theoretically and practically scalable in higher-dimensional settings. The downside of networks is that their training is an NP‐hard problem \cite{froese2023training}. Consequently, finding a global minimum of the loss function is not guaranteed, and practical performance depends heavily on initialization, architecture, and density of training data.  In ongoing and future research work, we plan to rigorously analyze the convergence properties of such a hybrid schemes.  Another aspect is that of the incorporation of ``parsimonious'' tensor‐product Chebyshev polynomial series for multivariate manifolds, and that of treating the associated domain of convergence as a tunable hyperparameter.

\bibliographystyle{unsrt} 
\bibliography{references}  

@article{kazantzis2001invariant,
  title={On invariant manifolds of nonlinear discrete-time input-driven dynamical systems},
  author={Kazantzis, Nikolaos},
  journal={Physics Letters A},
  volume={292},
  number={1-2},
  pages={107--114},
  year={2001},
  publisher={Elsevier}
}

@article{barron1993universal,
  title={Universal approximation bounds for superpositions of a sigmoidal function},
  author={Barron, Andrew R},
  journal={IEEE Transactions on Information Theory},
  volume={39},
  number={3},
  pages={930--945},
  year={1993},
  publisher={IEEE}
}

@article{kazantzis2005model,
  title={A model-based characterization of the long-term asymptotic behavior of nonlinear discrete-time processes using invariance functional equations},
  author={Kazantzis, Nikolaos and Huynh, Nguyen and Good, Theresa A},
  journal={Computers \& Chemical Engineering},
  volume={29},
  number={11-12},
  pages={2346--2354},
  year={2005},
  publisher={Elsevier}
}

@article{kazantzis2003,
  title={Invariance Inducing Control of Nonlinear Discrete-Time Dynamical Systems},
  author={Kazantzis, Nikolaos},
  journal={Journal of Nonlinear Science},
  volume={13},
  pages={579-601},
  year={2003},
  publisher={Elsevier}
}

@article{castillo1993nonlinear,
  title={Nonlinear regulation for a class of discrete-time systems},
  author={Castillo, B and Di Gennaro, Stefano and Monaco, Salvatore and Normand-Cyrot, D},
  journal={Systems \& Control Letters},
  volume={20},
  number={1},
  pages={57--65},
  year={1993},
  publisher={Elsevier}
}

@article{astolfi2003immersion,
  title={Immersion and invariance: A new tool for stabilization and adaptive control of nonlinear systems},
  author={Astolfi, Alessandro and Ortega, Romeo},
  journal={IEEE Transactions on Automatic Control},
  volume={48},
  number={4},
  pages={590--606},
  year={2003},
  publisher={IEEE}
}

@article{hamzi2005controlled,
  title={The controlled center dynamics},
  author={Hamzi, Boumediene and Kang, Wei and Krener, Arthur J},
  journal={Multiscale Modeling \& Simulation},
  volume={3},
  number={4},
  pages={838--852},
  year={2005},
  publisher={SIAM}
}

@article{kazantzis2005singular,
  title={Singular control-invariance PDEs for nonlinear systems},
  author={Kazantzis, Nikolaos and Demetriou, Michael A},
  journal={Multiscale Modeling \& Simulation},
  volume={3},
  number={4},
  pages={731--748},
  year={2005},
  publisher={SIAM}
}

@book{kuehn2015multiple,
  title={Multiple time scale dynamics},
  author={Kuehn, Christian},
  volume={191},
  year={2015},
  publisher={Springer}
}

@book{arnold2012geometrical,
  title={Geometrical methods in the theory of ordinary differential equations},
  author={Arnold, Vladimir Igorevich},
  volume={250},
  year={2012},
  publisher={Springer Science \& Business Media}
}

@book{gorban2005invariant,
  title={Invariant manifolds for physical and chemical kinetics},
  author={Gorban, Aleksandr Nikolaevich and Karlin, Ilya V},
  volume={660},
  year={2005},
  publisher={Springer}
}

@book{byrnes1997output,
  title={Output Regulation of Uncertain Nonlinear Systems},
  author={Byrnes, Christopher I and Priscoli, Francesco Delli and Isidori, Alberto},
    volume = {1},
  publisher={Birkhauser},
  year={1997}
}

@book{carr2012applications,
  title={Applications of centre manifold theory},
  author={Carr, Jack},
  volume={35},
  year={2012},
  publisher={Springer Science \& Business Media}
}

@book{colonius2012dynamics,
  title={The dynamics of control},
  author={Colonius, Fritz and Kliemann, Wolfgang},
  year={2012},
  publisher={Springer Science \& Business Media}
}

@article{gorban2003method,
  title={Method of invariant manifold for chemical kinetics},
  author={Gorban, Alexander N and Karlin, Iliya V},
  journal={Chemical Engineering Science},
  volume={58},
  number={21},
  pages={4751--4768},
  year={2003},
  publisher={Elsevier}
}

@article{roussel1991geometry,
  title={On the geometry of transient relaxation},
  author={Roussel, Marc R and Fraser, Simon J},
  journal={The Journal of Chemical Physics},
  volume={94},
  number={11},
  pages={7106--7113},
  year={1991},
  publisher={American Institute of Physics}
}

@inproceedings{goussis1992study,
  title={A study of homogeneous methanol oxidation kinetics using CSP},
  author={Goussis, Dimitris A and Lam, Sau-Hai},
  booktitle={Symposium (International) on Combustion},
  volume={24},  
  pages={113--120},
  year={1992},
  organization={Elsevier}
}

@article{goussis2006efficient,
  title={An efficient iterative algorithm for the approximation of the fast and slow dynamics of stiff systems},
  author={Goussis, Dimitris A and Valorani, Mauro},
  journal={Journal of Computational Physics},
  volume={214},
  number={1},
  pages={316--346},
  year={2006},
  publisher={Elsevier}
}

@article{maas1992simplifying,
  title={Simplifying chemical kinetics: intrinsic low-dimensional manifolds in composition space},
  author={Maas, Ulrich and Pope, Stephen B},
  journal={Combustion and Flame},
  volume={88},
  number={3-4},
  pages={239--264},
  year={1992},
  publisher={Elsevier}
}

@article{gear2005projecting,
  title={Projecting to a slow manifold: Singularly perturbed systems and legacy codes},
  author={Gear, C William and Kaper, Tasso J and Kevrekidis, Ioannis G and Zagaris, Antonios},
  journal={SIAM Journal on Applied Dynamical Systems},
  volume={4},
  number={3},
  pages={711--732},
  year={2005},
  publisher={SIAM}
}

@article{weierstrass1885analytische,
  title={{\"U}ber die analytische Darstellbarkeit sogenannter willk{\"u}rlicher Functionen einer reellen Ver{\"a}nderlichen},
  author={Weierstrass, Karl},
  journal={Sitzungsberichte der K{\"o}niglich Preu{\ss}ischen Akademie der Wissenschaften zu Berlin},
  volume={2},
  number={633-639},
  pages={364},
  year={1885}
}

@book{rudin1964principles,
  title={Principles of mathematical analysis},
  author={Rudin, Walter and others},
  volume={3},
  year={1964},
  publisher={McGraw-hill New York}
}

@article{tang2023physics,
  title={Physics-informed neural networks combined with polynomial interpolation to solve nonlinear partial differential equations},
  author={Tang, Siping and Feng, Xinlong and Wu, Wei and Xu, Hui},
  journal={Computers \& Mathematics with Applications},
  volume={132},
  pages={48--62},
  year={2023},
  publisher={Elsevier}
}

@article{yang2020neural,
  title={Neural network algorithm based on Legendre improved extreme learning machine for solving elliptic partial differential equations},
  author={Yang, Yunlei and Hou, Muzhou and Sun, Hongli and Zhang, Tianle and Weng, Futian and Luo, Jianshu},
  journal={Soft Computing},
  volume={24},
  pages={1083--1096},
  year={2020},
  publisher={Springer}
}

@article{mall2017single,
  title={Single layer Chebyshev neural network model for solving elliptic partial differential equations},
  author={Mall, Susmita and Chakraverty, Snehashish},
  journal={Neural Processing Letters},
  volume={45},
  pages={825--840},
  year={2017},
  publisher={Springer}
}

@article{hagan1994training,
  title={Training feedforward networks with the Marquardt algorithm},
  author={Hagan, Martin T and Menhaj, Mohammad B},
  journal={IEEE Transactions on Neural Networks},
  volume={5},
  number={6},
  pages={989--993},
  year={1994},
  publisher={IEEE}
}

@article{cybenko1989approximation,
  title={Approximation by superpositions of a sigmoidal function},
  author={Cybenko, George},
  journal={Mathematics of Control, Signals and Systems},
  volume={2},
  number={4},
  pages={303--314},
  year={1989},
  publisher={Springer}
}

@article{hornik1989multilayer,
  title={Multilayer feedforward networks are universal approximators},
  author={Hornik, Kurt and Stinchcombe, Maxwell and White, Halbert},
  journal={Neural Networks},
  volume={2},
  number={5},
  pages={359--366},
  year={1989},
  publisher={Elsevier}
}

@inproceedings{glorot2010understanding,
  title={Understanding the difficulty of training deep feedforward neural networks},
  author={Glorot, Xavier and Bengio, Yoshua},
  booktitle={Proceedings of the thirteenth international conference on artificial intelligence and statistics},
  pages={249--256},
  year={2010},
  organization={JMLR Workshop and Conference Proceedings}
}

@article{jones1995geometric,
  title={Geometric singular perturbation theory},
  author={Jones, Christopher KRT},
  journal={Dynamical Systems},
  pages={44--118},
  year={1995},
  publisher={Springer}
}

@article{fenichel1979geometric,
  title={Geometric singular perturbation theory for ordinary differential equations},
  author={Fenichel, Neil},
  journal={Journal of Differential Equations},
  volume={31},
  number={1},
  pages={53--98},
  year={1979},
  publisher={Academic Press}
}

@article{hirsch1970invariant,
  title={Invariant manifolds},
  author={Hirsch, Morris W and Pugh, Charles Chapman and Shub, Michael},
  journal={Bulletin of the American Mathematical Society},
  volume={76},
  number={5},
  pages={1015--1019},
  year={1970}
}

@article{cabre2003parameterization,
  title={The parameterization method for invariant manifolds I: manifolds associated to non-resonant subspaces},
  author={Cabr{\'e}, Xavier and Fontich, Ernest and De la Llave, Rafael},
  journal={Indiana University mathematics journal},
  pages={283--328},
  year={2003},
  publisher={JSTOR}
}

@article{dellnitz1997subdivision,
  title={A subdivision algorithm for the computation of unstable manifolds and global attractors},
  author={Dellnitz, Michael and Hohmann, Andreas},
  journal={Numerische Mathematik},
  volume={75},
  pages={293--317},
  year={1997},
  publisher={Springer}
}

@inproceedings{lam1989understanding,
  title={Understanding complex chemical kinetics with computational singular perturbation},
  author={Lam, Sau-Hai and Goussis, Dimitris A},
  booktitle={Symposium (International) on Combustion},
  volume={22},
  pages={931--941},
  year={1989},
  organization={Elsevier}
}

@incollection{goussis2011model,
  title={Model reduction for combustion chemistry},
  author={Goussis, Dimitris A and Maas, Ulrich},
  booktitle={Turbulent Combustion Modeling},
  pages={193--220},
  year={2011},
  publisher={Springer}
}

@article{zagaris2009analysis,
  title={Analysis of the accuracy and convergence of equation-free projection to a slow manifold},
  author={Zagaris, Antonios and Gear, C William and Kaper, Tasso J and Kevrekidis, Yannis G},
  journal={ESAIM: Mathematical Modelling and Numerical Analysis},
  volume={43},
  number={4},
  pages={757--784},
  year={2009},
  publisher={EDP Sciences}
}

@article{haller2016nonlinear,
  title={Nonlinear normal modes and spectral submanifolds: existence, uniqueness and use in model reduction},
  author={Haller, George and Ponsioen, Sten},
  journal={Nonlinear Dynamics},
  volume={86},
  pages={1493--1534},
  year={2016},
  publisher={Springer}
}

@article{tingas2018chemical,
  title={Chemical kinetic insights into the ignition dynamics of n-hexane},
  author={Tingas, Efstathios-Al and Wang, Zhandong and Sarathy, S Mani and Im, Hong G and Goussis, Dimitris A},
  journal={Combustion and Flame},
  volume={188},
  pages={28--40},
  year={2018},
  publisher={Elsevier}
}

@article{manias2016mechanism,
  title={The mechanism by which CH2O and H2O2 additives affect the autoignition of CH4/air mixtures},
  author={Manias, Dimitris M and Tingas, Efstathios Al and Frouzakis, Christos E and Boulouchos, Konstantinos and Goussis, Dimitris A},
  journal={Combustion and Flame},
  volume={164},
  pages={111--125},
  year={2016},
  publisher={Elsevier}
}

@article{patsatzis2016asymptotic,
  title={Asymptotic analysis of a target-mediated drug disposition model: algorithmic and traditional approaches},
  author={Patsatzis, Dimitris G and Maris, Dimitris T and Goussis, Dimitris A},
  journal={Bulletin of Mathematical Biology},
  volume={78},
  pages={1121--1161},
  year={2016},
  publisher={Springer}
}

@article{haro2016parameterization,
  title={The parameterization method for invariant manifolds},
  author={Haro, Alex and Canadell, Marta and Figueras, Jordi-Lluis and Luque, Alejandro and Mondelo, Josep-Maria},
  journal={Applied Mathematical Sciences},
  volume={195},
  year={2016},
  publisher={Springer}
}

@book{wiggins2003introduction,
  title={Introduction to Applied Nonlinear Dynamical Systems and Chaos},
  author={Wiggins, Stephen},
  volume={2},
  series={Texts in Applied Mathematics},
  year={2003},
  publisher={Springer}
}

@book{guckenheimer2013nonlinear,
  title={Nonlinear oscillations, dynamical systems, and bifurcations of vector fields},
  author={Guckenheimer, John and Holmes, Philip},
  volume={42},
  year={2013},
  publisher={Springer Science \& Business Media}
}

@book{rivlin1981introduction,
  title={An introduction to the approximation of functions},
  author={Rivlin, Theodore J},
  year={1981},
  publisher={Courier Corporation}
}

@article{froese2023training,
  title={Training neural networks is np-hard in fixed dimension},
  author={Froese, Vincent and Hertrich, Christoph},
  journal={Advances in Neural Information Processing Systems},
  volume={36},
  pages={44039--44049},
  year={2023}
}

@article{mason1980near,
  title={Near-best multivariate approximation by Fourier series, Chebyshev series and Chebyshev interpolation},
  author={Mason, John C},
  journal={Journal of Approximation Theory},
  volume={28},
  number={4},
  pages={349--358},
  year={1980},
  publisher={Elsevier}
}

@article{koronaki2024nonlinear,
  title={Nonlinear dimensionality reduction then and now: AIMs for dissipative PDEs in the ML era},
  author={Koronaki, Eleni D and Evangelou, Nikolaos and Martin-Linares, Cristina P and Titi, Edriss S and Kevrekidis, Ioannis G},
  journal={Journal of Computational Physics},
  volume={506},
  pages={112910},
  year={2024},
  publisher={Elsevier}
}

@article{chin2022enabling,
  title={Enabling equation-free modeling via diffusion maps},
  author={Chin, Tracy and Ruth, Jacob and Sanford, Clayton and Santorella, Rebecca and Carter, Paul and Sandstede, Bj{\"o}rn},
  journal={Journal of Dynamics and Differential Equations},
  pages={1--20},
  year={2022},
  publisher={Springer}
}

@article{evangelou2023double,
  title={Double diffusion maps and their latent harmonics for scientific computations in latent space},
  author={Evangelou, Nikolaos and Dietrich, Felix and Chiavazzo, Eliodoro and Lehmberg, Daniel and Meila, Marina and Kevrekidis, Ioannis G},
  journal={Journal of Computational Physics},
  volume={485},
  pages={112072},
  year={2023},
  publisher={Elsevier}
}

@article{kevrekidis2003equation,
  year = {2003},
  publisher = {International Press of Boston},
  volume = {1},
  number = {4},
  pages = {715--762},
  author = {Kevrekidis, Ioannis G. and  Gear, C. William and Hyman, James M. and Kevrekidis, Panagiotis G. and Runborg, Olof and Theodoropoulos, Constantinos },
  title = {Equation-Free,  Coarse-Grained Multiscale Computation: Enabling Microscopic Simulators to Perform System-Level Analysis},
  journal = {Communications in Mathematical Sciences}
}

@article{patsatzis2023data,
  title={Data-driven control of agent-based models: An equation/variable-free machine learning approach},
  author={Patsatzis, Dimitrios G and Russo, Lucia and Kevrekidis, Ioannis G and Siettos, Constantinos},
  journal={Journal of Computational Physics},
volume={478},
  pages={111953},
  year={2023},
  publisher={Elsevier}
}

@article{karniadakis2021physics,
  title={Physics-informed machine learning},
  author={Karniadakis, George Em and Kevrekidis, Ioannis G and Lu, Lu and Perdikaris, Paris and Wang, Sifan and Yang, Liu},
  journal={Nature Reviews Physics},
  volume={3},
  number={6},
  pages={422--440},
  year={2021},
  publisher={Nature Publishing Group}
}

@article{raissi2019physics,
  title={Physics-informed neural networks: A deep learning framework for solving forward and inverse problems involving nonlinear partial differential equations},
  author={Raissi, Maziar and Perdikaris, Paris and Karniadakis, George E},
  journal={Journal of Computational Physics},
  volume={378},
  pages={686--707},
  year={2019},
  publisher={Elsevier}
}

@article{wang2021understanding,
  title={Understanding and mitigating gradient flow pathologies in physics-informed neural networks},
  author={Wang, Sifan and Teng, Yujun and Perdikaris, Paris},
  journal={SIAM Journal on Scientific Computing},
  volume={43},
  number={5},
  pages={A3055--A3081},
  year={2021},
  publisher={SIAM}
}

@article{patsatzis2024slow,
  title={Slow invariant manifolds of singularly perturbed systems via physics-informed machine learning},
  author={Patsatzis, Dimitrios and Fabiani, Gianluca and Russo, Lucia and Siettos, Constantinos},
  journal={SIAM Journal on Scientific Computing},
  volume={46},
  number={4},
  pages={C297--C322},
  year={2024},
  publisher={SIAM}
}

@article{patsatzis2024slowB,
  title={Slow Invariant Manifolds of Fast-Slow Systems of ODEs with Physics-Informed Neural Networks},
  author={Patsatzis, Dimitrios G and Russo, Lucia and Siettos, Constantinos},
  journal={SIAM Journal on Applied Dynamical Systems},
  volume={23},
  number={4},
  pages={3077--3122},
  year={2024},
  publisher={SIAM}
}

@article{siettos2022numerical,
  title={A numerical method for the approximation of stable and unstable manifolds of microscopic simulators},
  author={Siettos, Constantinos and Russo, Lucia},
  journal={Numerical Algorithms},
  volume={89},
  number={3},
  pages={1335--1368},
  year={2022},
  publisher={Springer}
}

@article{siettos2014equation,
  title={Equation-free computation of coarse-grained center manifolds of microscopic simulators},
  author={Siettos, Constantinos},
  journal={Journal of Computational Dynamics},
  volume={1},
  number={2},
  pages={377--389},
  year={2014},
  publisher={Journal of Computational Dynamics}
}

@article{zagaris2004fast,
  title={Fast and slow dynamics for the computational singular perturbation method},
  author={Zagaris, Antonios and Kaper, Hans G and Kaper, Tasso J},
  journal={Multiscale Modeling \& Simulation},
  volume={2},
  number={4},
  pages={613--638},
  year={2004},
  publisher={SIAM}
}

@article{roberts1989utility,
  title={The utility of an invariant manifold description of the evolution of a dynamical system},
  author={Roberts, Antony J},
  journal={SIAM Journal on Mathematical Analysis},
  volume={20},
  number={6},
  pages={1447--1458},
  year={1989},
  publisher={SIAM}
}

@article{benedetti2026numerical,
    author = {Benedetti, Kaio C. B. and Gonçalves, Paulo B. and Lenci, Stefano and Rega, Giuseppe},
    title = {Numerical computation of the stable and unstable manifolds of saddles of randomly perturbed dynamical systems: An operator approach},
    journal = {Chaos: An Interdisciplinary Journal of Nonlinear Science},
    volume = {36},
    number = {3},
    pages = {033121},
    year = {2026},
    month = {03},
    issn = {1054-1500},
    doi = {10.1063/5.0307783},
    url = {https://doi.org/10.1063/5.0307783},
    eprint = {https://pubs.aip.org/aip/cha/article-pdf/doi/10.1063/5.0307783/20937547/033121_1_5.0307783.pdf},
}

@book{trefethen2019approximation,
  title={Approximation theory and approximation practice, extended edition},
  author={Trefethen, Lloyd N},
  year={2019},
  publisher={SIAM}
}

@article{jiang2020cooperative,
  title={Cooperative adaptive optimal output regulation of nonlinear discrete-time multi-agent systems},
  author={Jiang, Yi and Fan, Jialu and Gao, Weinan and Chai, Tianyou and Lewis, Frank L},
  journal={Automatica},
  volume={121},
  pages={109149},
  year={2020},
  publisher={Elsevier}
}

@article{bin2019adaptive,
  title={Adaptive output regulation for linear systems via discrete-time identifiers},
  author={Bin, Michelangelo and Marconi, Lorenzo and Teel, Andrew R},
  journal={Automatica},
  volume={105},
  pages={422--432},
  year={2019},
  publisher={Elsevier}
}

@book{huang2004nonlinear,
  title={Nonlinear output regulation: theory and applications},
  author={Huang, Jie},
  year={2004},
  publisher={SIAM}
}

@article{valorani2005higher,
  title={Higher order corrections in the approximation of low-dimensional manifolds and the construction of simplified problems with the CSP method},
  author={Valorani, Mauro and Goussis, Dimitris A and Creta, Francesco and Najm, Habib N},
  journal={Journal of Computational Physics},
  volume={209},
  number={2},
  pages={754--786},
  year={2005},
  publisher={Elsevier}
}

@article{roberts2014dynamical,
  title={A dynamical systems approach to simulating macroscale spatial dynamics in multiple dimensions},
  author={Roberts, AJ and MacKenzie, Tony and Bunder, Judith E},
  journal={Journal of Engineering Mathematics},
  volume={86},
  number={1},
  pages={175--207},
  year={2014},
  publisher={Springer}
}

@article{ginoux2008slow,
  title={Slow invariant manifolds as curvature of the flow of dynamical systems},
  author={Ginoux, Jean-Marc and Rossetto, Bruno and Chua, Leon O},
  journal={International Journal of Bifurcation and Chaos},
  volume={18},
  number={11},
  pages={3409--3430},
  year={2008},
  publisher={World Scientific}
}

@article{valorani2015dynamical,
  title={Dynamical system analysis of ignition phenomena using the tangential stretching rate concept},
  author={Valorani, Mauro and Paolucci, Samuel and Martelli, Emanuele and Grenga, Temistocle and Ciottoli, Pietro P},
  journal={Combustion and Flame},
  volume={162},
  number={8},
  pages={2963--2990},
  year={2015},
  publisher={Elsevier}
}

@inproceedings{lokshtanov2017beating,
  title={Beating brute force for systems of polynomial equations over finite fields},
  author={Lokshtanov, Daniel and Paturi, Ramamohan and Tamaki, Suguru and Williams, Ryan and Yu, Huacheng},
  booktitle={Proceedings of the Twenty-Eighth Annual ACM-SIAM Symposium on Discrete Algorithms},
  pages={2190--2202},
  year={2017},
  organization={SIAM}
}

@article{bando1995dynamical,
  title={Dynamical model of traffic congestion and numerical simulation},
  author={Bando, Masako and Hasebe, Katsuya and Nakayama, Akihiro and Shibata, Akihiro and Sugiyama, Yuki},
  journal={Physical review E},
  volume={51},
  number={2},
  pages={1035},
  year={1995},
  publisher={APS}
}

@article{argall2002rigorous,
  title={A rigorous treatment of a follow-the-leader traffic model with traffic lights present},
  author={Argall, Brenna and Hinde, Colin and Cheleshkin, Eugene and Greenberg, James M and Lin, Pei-Jen},
  journal={SIAM Journal on Applied Mathematics},
  volume={63},
  number={1},
  pages={149--168},
  year={2002},
  publisher={SIAM}
}

@article{aw2002derivation,
  title={Derivation of continuum traffic flow models from microscopic follow-the-leader models},
  author={Aw, AATM and Klar, Axel and Rascle, Michel and Materne, Thorsten},
  journal={SIAM Journal on applied mathematics},
  volume={63},
  number={1},
  pages={259--278},
  year={2002},
  publisher={SIAM}
}

@article{tordeux2018traffic,
  title={From traffic and pedestrian follow-the-leader models with reaction time to first order convection-diffusion flow models},
  author={Tordeux, Antoine and Costeseque, Guillaume and Herty, Michael and Seyfried, Armin},
  journal={SIAM journal on applied mathematics},
  volume={78},
  number={1},
  pages={63--79},
  year={2018},
  publisher={SIAM}
}

@article{chiarello2021multiscale,
  title={Multiscale control of generic second order traffic models by driver-assist vehicles},
  author={Chiarello, Felisia Angela and Piccoli, Benedetto and Tosin, Andrea},
  journal={Multiscale Modeling \& Simulation},
  volume={19},
  number={2},
  pages={589--611},
  year={2021},
  publisher={SIAM}
}

@article{byrnes2000output,
  title={Output regulation for nonlinear systems: an overview},
  author={Byrnes, CI and Isidori, A},
  journal={International Journal of Robust and Nonlinear Control: IFAC-Affiliated Journal},
  volume={10},
  number={5},
  pages={323--337},
  year={2000},
  publisher={Wiley Online Library}
}

@article{chen2005robust,
  title={Robust output regulation with nonlinear exosystems},
  author={Chen, Zhiyong and Huang, Jie},
  journal={Automatica},
  volume={41},
  number={8},
  pages={1447--1454},
  year={2005},
  publisher={Elsevier}
}

@article{maffettone2026bio,
  title={Bio-inspired density control of multi-agent swarms via leader-follower plasticity},
  author={Maffettone, Gian Carlo and Boldini, Alain and di Bernardo, Mario and Porfiri, Maurizio},
  journal={Automatica},
  volume={188},
  pages={112924},
  year={2026},
  publisher={Elsevier}
}

@article{gu2009leader,
  title={Leader--follower flocking: algorithms and experiments},
  author={Gu, Dongbing and Wang, Zongyao},
  journal={IEEE Transactions on Control Systems Technology},
  volume={17},
  number={5},
  pages={1211--1219},
  year={2009},
  publisher={IEEE}
}

@article{cristiani2021all,
  title={An all-leader agent-based model for turning and flocking birds},
  author={Cristiani, Emiliano and Menci, Marta and Papi, Marco and Brafman, L{\'e}onard},
  journal={Journal of Mathematical Biology},
  volume={83},
  number={4},
  pages={45},
  year={2021},
  publisher={Springer}
}

@book{treiber2012traffic,
  title={Traffic flow dynamics: data, models and simulation},
  author={Treiber, Martin and Kesting, Arne},
  year={2012},
  publisher={Springer Science \& Business Media}
}

@article{couzin2005effective,
  title={Effective leadership and decision-making in animal groups on the move},
  author={Couzin, Iain D and Krause, Jens and Franks, Nigel R and Levin, Simon A},
  journal={Nature},
  volume={433},
  number={7025},
  pages={513--516},
  year={2005},
  publisher={Nature Publishing Group UK London}
}

@article{shen2008cucker,
  title={Cucker--Smale flocking under hierarchical leadership},
  author={Shen, Jackie},
  journal={SIAM Journal on Applied Mathematics},
  volume={68},
  number={3},
  pages={694--719},
  year={2008},
  publisher={SIAM}
}

@article{li2010cucker,
  title={Cucker--Smale flocking under rooted leadership with fixed and switching topologies},
  author={Li, Zhuchun and Xue, Xiaoping},
  journal={SIAM Journal on Applied Mathematics},
  volume={70},
  number={8},
  pages={3156--3174},
  year={2010},
  publisher={SIAM}
  }

@article{dalmao2011cucker,
  title={Cucker--Smale flocking under hierarchical leadership and random interactions},
  author={Dalmao, Federico and Mordecki, Ernesto},
  journal={SIAM Journal on Applied Mathematics},
  volume={71},
  number={4},
  pages={1307--1316},
  year={2011},
  publisher={SIAM}
}

@article{kohler2021constrained,
  title={Constrained nonlinear output regulation using model predictive control},
  author={K{\"o}hler, Johannes and M{\"u}ller, Matthias A and Allg{\"o}wer, Frank},
  journal={IEEE Transactions on Automatic Control},
  volume={67},
  number={5},
  pages={2419--2434},
  year={2021},
  publisher={IEEE}
}

@article{cucker2008flocking,
  title={Flocking with informed agents},
  author={Cucker, Felipe and Huepe, Cristi{\'a}n},
  journal={Mathematics in Action},
  volume={1},
  number={1},
  pages={1--25},
  year={2008}
}

@article{cucker2007emergent,
  title={Emergent behavior in flocks},
  author={Cucker, Felipe and Smale, Steve},
  journal={IEEE Transactions on automatic control},
  volume={52},
  number={5},
  pages={852--862},
  year={2007},
  publisher={IEEE}
}

@article{gazis1961nonlinear,
  title={Nonlinear follow-the-leader models of traffic flow},
  author={Gazis, Denos C and Herman, Robert and Rothery, Richard W},
  journal={Operations research},
  volume={9},
  number={4},
  pages={545--567},
  year={1961},
  publisher={INFORMS}
}

@article{ni2010leader,
  title={Leader-following consensus of multi-agent systems under fixed and switching topologies},
  author={Ni, Wei and Cheng, Daizhan},
  journal={Systems \& control letters},
  volume={59},
  number={3-4},
  pages={209--217},
  year={2010},
  publisher={Elsevier}
}

@article{delgado2026sindy,
  title={Sindy on attracting manifolds},
  author={Delgado-Cano, Diemen and Kracht, Erick and Fasel, Urban and Herrmann, Benjamin},
  journal={Nonlinear Dynamics},
  volume={114},
  number={1},
  pages={20},
  year={2026},
  publisher={Springer}
}

@article{linot2020deep,
  title={Deep learning to discover and predict dynamics on an inertial manifold},
  author={Linot, Alec J and Graham, Michael D},
  journal={Physical Review E},
  volume={101},
  number={6},
  pages={062209},
  year={2020},
  publisher={APS}
}

@article{ghadami2022deep,
  title={Deep learning for centre manifold reduction and stability analysis in nonlinear systems},
  author={Ghadami, Amin and Epureanu, Bogdan I},
  journal={Philosophical Transactions of the Royal Society A: Mathematical, Physical and Engineering Sciences},
  volume={380},
  number={2229},
  year={2022},
  publisher={The Royal Society}
}

@article{breden2020computing,
  title={Computing invariant sets of random differential equations using polynomial chaos},
  author={Breden, Maxime and Kuehn, Christian},
  journal={SIAM Journal on Applied Dynamical Systems},
  volume={19},
  number={1},
  pages={577--618},
  year={2020},
  publisher={SIAM}
}

@article{kuehn2025fast,
  title={Fast reactions and slow manifolds},
  author={Kuehn, Christian and Sulzbach, Jan-Eric},
  journal={Nonlinear Differential Equations and Applications NoDEA},
  volume={32},
  number={4},
  pages={72},
  year={2025},
  publisher={Springer}
}

@article{guckenheimer2009computing,
  title={Computing slow manifolds of saddle type},
  author={Guckenheimer, John and Kuehn, Christian},
  journal={SIAM Journal on Applied Dynamical Systems},
  volume={8},
  number={3},
  pages={854--879},
  year={2009},
  publisher={SIAM}
}

@article{chen2024deep,
  title={Deep neural network approximations for the stable manifolds of the Hamilton--Jacobi--Bellman equations},
  author={Chen, Guoyuan},
  journal={IEEE Transactions on Automatic Control},
  volume={69},
  number={10},
  pages={7239--7246},
  year={2024},
  publisher={IEEE}
}

@article{touze2021model,
  title={Model order reduction methods for geometrically nonlinear structures: a review of nonlinear techniques},
  author={Touz{\'e}, Cyril and Vizzaccaro, Alessandra and Thomas, Olivier},
  journal={Nonlinear Dynamics},
  volume={105},
  number={2},
  pages={1141--1190},
  year={2021},
  publisher={Springer}
}

@article{vizzaccaro2021direct,
  title={Direct computation of nonlinear mapping via normal form for reduced-order models of finite element nonlinear structures},
  author={Vizzaccaro, Alessandra and Shen, Yichang and Salles, Lo{\"\i}c and Blaho{\v{s}}, Ji{\v{r}}{\'\i} and Touz{\'e}, Cyril},
  journal={Computer Methods in Applied Mechanics and Engineering},
  volume={384},
  pages={113957},
  year={2021},
  publisher={Elsevier}
}

@article{chui2018deep,
  title={Deep nets for local manifold learning},
  author={Chui, Charles K and Mhaskar, Hrushikesh N},
  journal={Frontiers in Applied Mathematics and Statistics},
  volume={4},
  pages={12},
  year={2018},
  publisher={Frontiers Media SA}
}

@article{zhang2024car,
  title={Car-following models: A multidisciplinary review},
  author={Zhang, Tianya Terry and Jin, Peter J and McQuade, Sean T and Bayen, Alexandre and Piccoli, Benedetto},
  journal={IEEE Transactions on Intelligent Vehicles},
  year={2024},
  publisher={IEEE}
}

@book{cristiani2014multiscale,
  title={Multiscale modeling of pedestrian dynamics},
  author={Cristiani, Emiliano and Piccoli, Benedetto and Tosin, Andrea},
  year={2014},
  publisher={Springer}
}

\clearpage
\newpage
\section*{APPENDIX}
\appendix
\renewcommand{\theequation}{A.\arabic{equation}}
\renewcommand{\thefigure}{A.\arabic{figure}}
\setcounter{equation}{0}
\setcounter{figure}{0}
\section{Invariant Manifold approximations via Power Series Expansions}
\label[appendix]{app:IM_PS}
We hereby describe the approach proposed by \cite{kazantzis2001invariant,kazantzis2005model} for deriving approximations of the IM mapping in \cref{eq:IM} by approximating the solution of the system of NFEs in \cref{eq:NFEs_IM} via power series expansions (PSE).~We assume that the IM functional $\mathbf{x} = \boldsymbol{\pi}(\mathbf{y}) \in \mathbb{R}^N$ can be locally approximated by a multivariate Taylor series expansion around the equilibrium $(\mathbf{x}_0,\mathbf{y}_0)=(\mathbf{0}^N,\mathbf{0}^M)$ of degree $h$.~The PSE can be compactly written in the form:
\begin{equation}
    \tilde{\boldsymbol{\pi}}^{PSE}(\mathbf{y}) = [\tilde{\pi}_1^{PSE}(\mathbf{y}), \ldots, \tilde{\pi}_n^{PSE}(\mathbf{y}), \ldots, \tilde{\pi}_N^{PSE}(\mathbf{y})]^\top
\end{equation}
where each element for $n=1,\ldots,N$ can be expressed as:
\begin{equation}
    \tilde{\pi}_n^{PSE}(\mathbf{y}) = \pi_n^0 + \sum_{i_1=1}^M \pi_n^{i_1} y_{i_1} + \dfrac{1}{2!} \sum_{i_1,i_2=1}^M \pi_n^{i_1,i_2} y_{i_1} y_{i_2} + \ldots + \dfrac{1}{h!} \sum_{i_1,\ldots,i_h=1}^M \pi_n^{i_1,\ldots,i_h} y_{i_1} \cdots y_{i_h} + \mathcal{O}(\mathbf{y}^{h+1}). 
  %
    \label{eq:IM_PSexp}
\end{equation}
Each term in the expansion corresponds to the zeroth, first, second, $\ldots$ and $h$-th order approximations with respect to the components of $\mathbf{y}\in\mathbb{R}^M$, while the reminder term $\mathcal{O}(\mathbf{y}^{h+1})$ captures the higher order contributions, including products of the form $y_{i_1} \cdots y_{i_h} y_{i_{h+1}}$.

The coefficients $\pi_n^{i_1,\ldots,i_h}$ of the PSE in \cref{eq:IM_PSexp} represent the $h$-th derivatives of the unknown IM functional at $(\mathbf{x}_0,\mathbf{y}_0)$.~The zeroth order coefficient is set to $\pi_n^0=0$, since the IM approximation passes from the equilibrium; i.e., $\tilde{\boldsymbol{\pi}}^{PSE}(\mathbf{0}^M)=\mathbf{0}^N$.~To determine the remaining coefficients $\pi_n^{i_1,\ldots,i_h}$, the expansion in \cref{eq:IM_PSexp} is substituted into the NFEs in \cref{eq:NFEs_IM}, and the nonlinear functions $\mathbf{f}$ and $\mathbf{g}$ are expanded around the equilibrium.~By equating coefficients of the same order on both sides of \cref{eq:NFEs_IM}, the coefficients $\pi_n^{i_1,\ldots,i_h}$ can be computed recursively (i.e., the computation of the $h$-th order coefficients requires knowledge about the $(h-1)$-th order ones) through a system of linear algebraic equations \cite{kazantzis2001invariant,kazantzis2005model}.

This approach bears similarities to the derivation of center and slow invariant manifold approximations for continuous-time dynamical systems \cite{gear2005projecting,gorban2003method,kuehn2015multiple,roussel1991geometry}, which also employ (asymptotic) series expansions.~Similarly to the continuous case, this approach relies on the availability of the analytical expressions of $\mathbf{f}$ and $\mathbf{g}$.~By leveraging their Taylor expansions (which are accurate sufficiently far from singularities), the coefficients of the PSE in \cref{eq:IM_PSexp} can be determined through the solution of a linear system of algebraic equations \cite{kazantzis2001invariant,kazantzis2005model}.~For a tensorial notation of the above approach, see \cite{kazantzis2005model}.~We hereby note that as the dimension of the system $N+M$, and/or the desired order $h$ of the PSE increases, this approach becomes increasingly challenging to implement, even with the aid of symbolic computation software, due to its combinatorially increasing number of coefficients to determine.

\renewcommand{\theequation}{B.\arabic{equation}}
\renewcommand{\thefigure}{B.\arabic{figure}}
\setcounter{equation}{0}
\setcounter{figure}{0}
\section{Residuals and Derivatives for the Hybrid and standalone NN Scheme}
\label[appendix]{app:OptPIHS_LM}
In this appendix, we formulate the physics-informed optimization problem in \cref{eq:min_hybrid} as a nonlinear residual minimization problem.~We begin by deriving the residuals and their analytical derivatives for the hybrid scheme introduced in \cref{sb:PI_opt_hyb}.~We then show how these residuals and derivatives specialize to the standalone neural network (NN) scheme, which is obtained by omitting the polynomial series component and setting $\omega_{\partial \mathcal{D}}=0$.~This unified treatment provides the foundation for the Levenberg–Marquardt (LM) algorithm that follows in \cref{app:LMalg}.

\subsection{Hybrid Scheme}
\label[appendix]{sb:hyb_scheme_ders}
First, we collect all the components of the hybrid approximation into the $N$-dim. column vector:
\begin{equation}
    \tilde{\boldsymbol{\pi}}^{HxS}(\mathbf{y},\boldsymbol{a},\mathbf{p}) = [\tilde{\pi}^{HxS}_1(\mathbf{y},\boldsymbol{a}_1,\mathbf{p}_1), \ldots, \tilde{\pi}^{HxS}_N(\mathbf{y},\boldsymbol{a}_N,\mathbf{p}_N)]^\top,
\end{equation}
where each element $\tilde{\pi}^{HxS}_n(\mathbf{y},\boldsymbol{a}_n,\mathbf{p}_n)$ is given by \cref{eq:IM_hybrid} for $n=1,\ldots,N$.~In the ensuing text, we drop the dependence on the hyperparameters $\mathcal{H}_1$, $\mathcal{H}_2$ and the radius of the hyperrectangle $\mathbf{r}$ for the conciseness.~Let us further collect all tunable parameters into the column vector of dimensions $N\left( L(M+2) + \binom{M+h}{h} +1 \right)$, as:
\begin{equation}
    \boldsymbol{\nu} = [\boldsymbol{\nu}_1,\ldots,\boldsymbol{\nu}_n,\ldots,\boldsymbol{\nu}_N]^\top =[\boldsymbol{a}_1,\mathbf{p}_1, \ldots, \boldsymbol{a}_n,\mathbf{p}_n, \ldots, \boldsymbol{a}_N,\mathbf{p}_N ]^\top,  
    \label{eq:pars_PIHxS}
\end{equation}
where each pair $\boldsymbol{\nu}_n=[\boldsymbol{a}_n, \mathbf{p}_n]$ for $n=1,\ldots,N$ includes the coefficients of the $n$-th element of the polynomial series in \cref{eq:IM_polySexp} and the parameters of the $n$-th output of the NN in \cref{eq:IM_NN}, respectively; i.e., all the parameters of the $n$-th output of the hybrid scheme in \cref{eq:IM_hybrid}.

Next, consider the set of collocation points $\mathbf{y}^{(q)}\in\Omega$ with $q=1,\ldots,Q$ and the set of the ones on the boundary $\mathbf{y}^{(r)}\in\partial\mathcal{D}$ with $r=1,\ldots,R$.~The non-linear residuals in \cref{eq:nonlin_res} for the optimization of the loss function in \cref{eq:min_hybrid} are then written, with the inclusion of the balancing weights $\omega_{\Omega}, \omega_0, \omega_{\partial \mathcal{D}}$, as:
\begin{enumerate}[label=(\roman*)]
\item NFEs satisfying residuals $\mathcal{F}^{HxS}_{q,n}(\boldsymbol{\nu})$ for $q=1,\ldots,Q$ and $n=1,\ldots,N$ in $\Omega$:
\begin{multline}
    \mathcal{F}^{HxS}_{q,n}(\boldsymbol{\nu}) = \omega_{\Omega} \Big(\tilde{\pi}^{HxS}_n(\mathbf{A}\mathbf{y}^{(q)} + \mathbf{g}(\mathbf{y}^{(q)}),\boldsymbol{\nu}_n) - \sum_{i=1}^N B_{n,i} \tilde{\pi}^{HxS}_i(\mathbf{y}^{(q)},\boldsymbol{\nu}_i) - \\ 
    \sum_{j=1}^M C_{n,j} y^{(q)}_j  - f_n(\tilde{\boldsymbol{\pi}}^{HxS}(\mathbf{y}^{(q)},\boldsymbol{\nu}),\mathbf{y}^{(q)}) \Big), \label{eq:res_PIHxS1}
\end{multline}
where $B_{n,i}$ and $C_{n,j}$ are the $(n,i)$ and $(n,j)$ elements of $\mathbf{B}\in\mathbb{R}^{N\times N}$ and $\mathbf{C}\in\mathbb{R}^{N\times M}$, respectively, $y^{(q)}_j$ is the $j$-th element of $\mathbf{y}^{(q)}=[y^{(q)}_1,\ldots,y^{(q)}_M]^\top \in \mathbb{R}^M$, and $f_n$ is the $n$-th element of $\mathbf{f}=[f_1,\ldots,f_N]^\top$.
\item Equilibrium-passing residuals $\mathcal{F}^{HxS}_{0,n}(\boldsymbol{\nu}_n)$ for $n=1,\ldots,N$:
\begin{equation}
    \mathcal{F}^{HxS}_{0,n}(\boldsymbol{\nu}_n) = \omega_0 \tilde{\pi}^{HxS}_n(\mathbf{0}^M,\boldsymbol{\nu}_n)= \omega_0 \tilde{\pi}^{xS}_n(\mathbf{0}^M,\boldsymbol{a}_n),
    \label{eq:res_PIHxS2}
\end{equation}
where the latter equality follows from \cref{eq:IM_hybrid}, since only the polynomial series is active at equilibrium.
\item $C^0$ continuity residuals $\mathcal{F}^{HxS}_{r,n}(\boldsymbol{\nu}_n)$ for $r=1,\ldots,R$ and $n=1,\ldots,N$ at $\partial \mathcal{D}$:
\begin{equation}
    \mathcal{F}^{HxS}_{r,n}(\boldsymbol{\nu}_n) = \omega_{\partial \mathcal{D}} \Big( \tilde{\pi}^{xS}_n(\mathbf{y}^{(r)},\boldsymbol{a}_n) - \tilde{\pi}^{NN}_n(\mathbf{y}^{(r)},\mathbf{p}_n) \Big).
    \label{eq:res_PIHxS3}
\end{equation}
\end{enumerate}
Note that the residuals in \cref{eq:res_PIHxS1} are computed on the basis of all, in general, the coefficients in $\boldsymbol{\nu}$, since all components of $\mathbf{y}$ are contributing to the linear and non-linear terms, $B_{n,i} \cdot \tilde{\pi}^{xS}_i(\mathbf{y}^{(q)},\boldsymbol{\nu}_i)$ and $f_n(\tilde{\boldsymbol{\pi}}^{xS}(\mathbf{y}^{(q)},\boldsymbol{\nu}),\mathbf{y}^{(q)})$, respectively.~Using the expressions of the residuals in \cref{eq:res_PIHxS1,eq:res_PIHxS2,eq:res_PIHxS3}, we then collect all the residuals into the $N(Q+R+1)$-dim. column vector:
\begin{multline}
    \boldsymbol{\mathcal{F}}^{HxS} = [\mathcal{F}^{HxS}_{1,1},\ldots,\mathcal{F}^{HxS}_{Q,1},\ldots,\mathcal{F}^{HxS}_{q,n},\ldots,\mathcal{F}^{HxS}_{1,N},\ldots,\mathcal{F}^{HxS}_{Q,N},\\ \mathcal{F}^{HxS}_{0,1},\ldots,\mathcal{F}^{HxS}_{0,n},\ldots,\mathcal{F}^{HxS}_{0,N}, \mathcal{F}^{HxS}_{1,1},\ldots,\mathcal{F}^{HxS}_{R,1},\ldots,\mathcal{F}^{HxS}_{r,n},\ldots,\mathcal{F}^{HxS}_{1,N},\ldots,\mathcal{F}^{HxS}_{R,N}]^\top.
    \label{eq:res_HS}
\end{multline}

For the minimization of the above residuals with the gradient-based optimization via the LM algorithm, the required Jacobian matrix $\mathbf{J}=\nabla_{\boldsymbol{\nu}} \boldsymbol{\mathcal{F}}^{HxS}$ is computed via symbolic differentiation.~To compute the Jacobian, let $\nu_k$ denote any coefficient belonging to the vector of coefficients $\boldsymbol{\nu}_k=[\boldsymbol{a}_k, \mathbf{p}_k]$ of the $k$-th element of $\boldsymbol{\nu}$ in \cref{eq:pars_PIHxS} for $k=1,\ldots,N$.~Then, the derivatives of the residuals in \cref{eq:res_PIHxS1,eq:res_PIHxS2,eq:res_PIHxS3} w.r.t. $\nu_k$ can be computed by:
\begin{align}
    \dfrac{\partial \mathcal{F}^{HxS}_{q,n}(\boldsymbol{\nu})}{\partial \nu_k} & = \omega_{\Omega} \Bigg(\dfrac{\partial \tilde{\pi}^{HxS}_n(\mathbf{A}\mathbf{y}^{(q)} + \mathbf{g}(\mathbf{y}^{(q)}),\boldsymbol{\nu}_n)}{\partial \nu_k} \delta_{nk} -  \sum_{i=1}^N \left(B_{n,i} + \dfrac{\partial f_n}{\partial x_i}\right) \dfrac{\partial \tilde{\pi}^{HxS}_i(\mathbf{y}^{(q)},\boldsymbol{\nu}_i)}{\partial \nu_k} \delta_{ik}\Bigg), \label{eq:resHxS_ders1} \\
    \dfrac{\partial \mathcal{F}^{HxS}_{0,n}(\boldsymbol{\nu}_n)}{\partial \nu_k} & = \omega_0 \dfrac{\partial \tilde{\pi}^{xS}_n(\mathbf{0}^M,\boldsymbol{a}_n)}{\partial \nu_k} \delta_{nk}, \label{eq:resHxS_ders2} \\
    \dfrac{\partial \mathcal{F}^{HxS}_{r,n}(\boldsymbol{\nu}_n)}{\partial \nu_k} & = \omega_{\partial \mathcal{D}} \Bigg( \dfrac{\partial \tilde{\pi}^{xS}_n(\mathbf{y}^{(r)},\boldsymbol{a}_n)}{\partial \nu_k} - \dfrac{\partial \tilde{\pi}^{NN}_n(\mathbf{y}^{(r)},\mathbf{p}_n)}{\partial \nu_k} \Bigg) \delta_{nk},
    \label{eq:resHxS_ders3}
\end{align}
where $\delta_{ij}$ is the Kronecher delta, and the last term in \cref{eq:resHxS_ders1} includes contributions from all $\tilde{\pi}^{HxS}_i(\cdot)$ outputs for $i=1,\ldots,N$, since, in general, the gradient $\partial f_n / \partial x_i$ of the vector field is considered full.

We highlight here that the derivatives in \cref{eq:resHxS_ders1} are computed on the basis of the hybrid scheme $\tilde{\pi}^{HxS}_n(\cdot)$, while the ones in \cref{eq:resHxS_ders2,eq:resHxS_ders3} are computed on the basis of the hybrid scheme's counterparts; i.e., the polynomial series $\tilde{\pi}^{xS}_n(\cdot)$ and the NN $\tilde{\pi}^{NN}_n(\cdot)$, respectively.~For computing the derivatives of the hybrid scheme in \cref{eq:resHxS_ders1}, we exploit the structure of the hybrid scheme in \cref{eq:IM_hybrid}.~Since the Heaviside function $H(\mathbf{y},\mathbf{r})$ is piecewise constant, differentiation directly yields the following piecewise expressions:
\begin{subnumcases}{\dfrac{\partial \tilde{\pi}^{HxS}_n(\mathbf{z}^{(q)},\boldsymbol{\nu}_n)}{\partial \nu_k} =} 
    \dfrac{\partial \tilde{\pi}^{xS}_n(\mathbf{z}^{(q)},\boldsymbol{a}_n)}{\partial \nu_k}, & when $\mathbf{z}^{(q)} \in \mathcal{D}$
    \label{eq:hyb_der1_xS} \\
    \dfrac{\partial \tilde{\pi}^{NN}_n(\mathbf{z}^{(q)},\mathbf{p}_n)}{\partial \nu_k}, & when $\mathbf{z}^{(q)} \notin \mathcal{D}$ \label{eq:hyb_der1_NN}
\end{subnumcases} 
\begin{subnumcases}{\dfrac{\partial \tilde{\pi}^{HxS}_i(\mathbf{y}^{(q)},\boldsymbol{\nu}_i)}{\partial \nu_k} =}
    \dfrac{\partial \tilde{\pi}^{xS}_i(\mathbf{y}^{(q)},\boldsymbol{a}_i)}{\partial \nu_k}, & when  $\mathbf{y}^{(q)} \in \mathcal{D}$ \label{eq:hyb_der2_xS} \\
    \dfrac{\partial \tilde{\pi}^{NN}_i(\mathbf{y}^{(q)},\mathbf{p}_i)}{\partial \nu_k}, & when  $\mathbf{y}^{(q)} \notin \mathcal{D}$ \label{eq:hyb_der2_NN}
\end{subnumcases}
where $\mathbf{z}^{(q)}=\mathbf{A}\mathbf{y}^{(q)} + \mathbf{g}(\mathbf{y}^{(q)})$ and $\mathcal{D}$ is the hyperrectangle of radius $\mathbf{r}$ within which the polynomial series are active and outside of which the NN acts.

Given the expressions in (\ref{eq:hyb_der1_xS}), (\ref{eq:hyb_der1_NN}), (\ref{eq:hyb_der2_xS}) and (\ref{eq:hyb_der2_NN}), the computation of the derivatives of the hybrid scheme's residuals in \cref{eq:resHxS_ders1,eq:resHxS_ders2,eq:resHxS_ders3} is simplified by the computation of the derivatives of the polynomial series and the NN w.r.t. any parameter $\nu_k$ belonging in $\boldsymbol{\nu}$ of \cref{eq:pars_PIHxS}.~However, $\boldsymbol{\nu}$ includes both the polynomial series coefficients $\boldsymbol{a}_k$ and the NN parameters $\mathbf{p}_k$.~Hence, in the case where:
\begin{enumerate}
    \item $\nu_k \equiv \alpha_k^{i_1,\ldots,i_M}\in \boldsymbol{a}_k$, the derivatives in (\ref{eq:hyb_der1_xS}), (\ref{eq:hyb_der2_xS}) for the coefficient $\alpha_k^{i_1,\ldots,i_M}$ corresponding to the $i_1$-th,\ldots, $i_M$-th multivariate monomial of \cref{eq:IM_polySexp}, are computed as:
    \begin{equation}
        \dfrac{\partial \tilde{\pi}^{xS}_n(\mathbf{y}^{(q)},\boldsymbol{a}_n)}{\partial \alpha_k^{i_1,\ldots,i_M}}=P^x_{i_1}(y_1^{(q)}) \cdots P^x_{i_M}(y_M^{(q)}), \label{eq:derPS_coef}
    \end{equation}
    for $k=n$ and zero otherwise.~Since $\nu_k \equiv \alpha_k^{i_1,\ldots,i_M} \notin \mathbf{p}_n$, the derivatives in (\ref{eq:hyb_der1_NN}, \ref{eq:hyb_der2_NN}) are also zero.~Similarly, the derivatives in \cref{eq:resHxS_ders2} are computed by \cref{eq:derPS_coef}, while the ones in \cref{eq:resHxS_ders3} are zero. 
    \item $\nu_k \in \mathbf{p}_k$, the related derivatives are computed via the symbolic differentiation of the NN output, detailed in \cref{sb:NN_scheme_ders} below.~In particular, the derivatives in (\ref{eq:hyb_der1_NN}), (\ref{eq:hyb_der2_NN}) are computed by \cref{eq:derNN_pars}, while the ones in (\ref{eq:hyb_der1_xS}), (\ref{eq:hyb_der2_xS}) are zero, since now $\nu_k \notin \boldsymbol{a}_k$.~Similarly, the derivatives in \cref{eq:resHxS_ders3} are computed by \cref{eq:derNN_pars}, while the ones in \cref{eq:resHxS_ders2} are zero.
\end{enumerate}
In summary, substitution of \cref{eq:derPS_coef} and \cref{eq:derNN_pars} into (\ref{eq:hyb_der1_xS}), (\ref{eq:hyb_der1_NN}), (\ref{eq:hyb_der2_xS}) and (\ref{eq:hyb_der2_NN}), and then to \cref{eq:resHxS_ders1,eq:resHxS_ders2,eq:resHxS_ders3} provides all the related derivatives for the formation of the $N(Q+R+1) \times N\left( L(M+2) + \binom{M+h}{h} +1 \right)$ Jacobian matrix $\mathbf{J}=\nabla_{\boldsymbol{\nu}} \boldsymbol{\mathcal{F}}^{HxS}$.~A detailed presentation of the iterative procedure of the LM algorithm for the minimization of the residuals is provided in \cref{app:LMalg}, along with a detailed pseudo-code.

\subsection{Neural Network Scheme (special case of Hybrid Scheme)}
\label[appendix]{sb:NN_scheme_ders}
The standalone NN scheme is recovered as a special case of the hybrid formulation in \cref{sb:hyb_scheme_ders} by omitting the polynomial series component, thus setting  the NN component active everywhere (equivalently, taking $\mathcal{D} = \emptyset$).~This results in simplifying the parameter vector $\boldsymbol{\nu}$ in \cref{eq:pars_PIHxS} into the $N(L(M+2)+1)$-dim. column vector:
\begin{equation}
    \mathbf{p} = [\mathbf{p}_1, \ldots, \mathbf{p}_n, \ldots, \mathbf{p}_N ]^\top,  
    \label{eq:pars_PINN}
\end{equation}
where each $\mathbf{p}_n = [\mathbf{w}^{o(n)},b^{o(n)},\mathbf{W}^{(n)},\mathbf{b}^{(n)}]^\top$ includes the weights and biases of the hidden and output layers for the $n$-th element of the SLFNN in \cref{eq:IM_NN}.

For the optimization of the loss function in \cref{eq:min_hybrid} using the standalone NN scheme, the absence of the polynomial series component results in neglecting the $C^0$ continuity residuals in \cref{eq:res_PIHxS3} (equivalently, setting $\omega_{\partial \mathcal{D}}=0$), and instead considering the equilibrium-passing residuals in \cref{eq:res_PIHxS2} for the NN component.~Hence, the residuals in \cref{eq:res_PIHxS1,eq:res_PIHxS2} reduce to:
\begin{align}
    \mathcal{F}^{NN}_{q,n}(\mathbf{p}) & = \omega_{\Omega} \bigg(\tilde{\pi}^{NN}_n(\mathbf{A}\mathbf{y}^{(q)} + \mathbf{g}(\mathbf{y}^{(q)}),\mathbf{p}_n) - \sum_{i=1}^N B_{n,i} \tilde{\pi}^{NN}_i(\mathbf{y}^{(q)},\mathbf{p}_i) - \nonumber \\
    & \qquad \qquad \qquad \qquad \qquad \qquad \qquad \qquad  \sum_{j=1}^M C_{n,j} y^{(q)}_j - f_n(\tilde{\boldsymbol{\pi}}^{NN}(\mathbf{y}^{(q)},\mathbf{p}),\mathbf{y}^{(q)})\bigg), \label{eq:res_PINN1} \\
    \mathcal{F}^{NN}_{0,n}(\mathbf{p}_n) & = \omega_0 \tilde{\pi}^{NN}_n(\mathbf{0}^M,\mathbf{p}_n),
    \label{eq:res_PINN2}
\end{align}
where $\mathbf{y}^{(q)}$ are the $q=1,\ldots,Q$ collocation points.~These residuals are collected into the $N(Q+1)$-dim. column vector:
\begin{equation}
    \boldsymbol{\mathcal{F}}^{NN} = [\mathcal{F}^{NN}_{1,1},\ldots,\mathcal{F}^{NN}_{Q,1},\ldots,\mathcal{F}^{NN}_{q,n},\ldots,\mathcal{F}^{NN}_{1,N},\ldots,\mathcal{F}^{NN}_{Q,N}, \mathcal{F}^{NN}_{0,1},\ldots,\mathcal{F}^{NN}_{0,n},\ldots,\mathcal{F}^{NN}_{0,N},]^\top,
    \label{eq:res_NN}
\end{equation}
in analogy to \cref{eq:res_HS}. 

For the computation of the Jacobian matrix $\mathbf{J}=\nabla_{\mathbf{p}} \boldsymbol{\mathcal{F}}^{NN}$, let $p_k$ denote any coefficient belonging to the vector of coefficients $\mathbf{p}_k$ of the $k$-th element of $\mathbf{p}$ in \cref{eq:pars_PINN}.~Then, the derivatives of the residuals in \cref{eq:res_PINN1,eq:res_PINN2} w.r.t. $p_k$ are given by:
\begin{align}
    \dfrac{\partial \mathcal{F}^{NN}_{q,n}(\mathbf{p})}{\partial p_k} & = \omega_{\Omega} \left( \dfrac{\partial \tilde{\pi}^{NN}_n(\mathbf{A}\mathbf{y}^{(q)} + \mathbf{g}(\mathbf{y}^{(q)}),\mathbf{p}_n)}{\partial p_k} \delta_{nk} - \sum_{i=1}^N \left(B_{n,i} + \dfrac{\partial f_n}{\partial x_i}\right) \dfrac{\partial \tilde{\pi}^{NN}_i(\mathbf{y}^{(q)} ,\mathbf{p}_i)}{\partial p_k} \delta_{ik} \right), \label{eq:resNN_ders1} \\
    \dfrac{\partial \mathcal{F}^{NN}_{0,n}(\mathbf{p}_n)}{\partial p_k} & = \omega_0 \dfrac{\partial \tilde{\pi}^{NN}_n(\mathbf{0}^M,\mathbf{p}_n)}{\partial p_k} \delta_{nk},
    \label{eq:resNN_ders2}
\end{align}
where the last term in \cref{eq:resNN_ders1} includes contributions from all the SLFNN outputs, due to, in general, a dense gradient $\partial f_n / \partial x_i$.~For the computation of the derivatives of the SLFNN included in \cref{eq:resNN_ders1,eq:resNN_ders2}, we now exploit the simple structure of the SLFNN in \cref{eq:IM_NN} and the fact that the derivative of the sigmoid activation function is $\phi'(\cdot)=\phi(\cdot)(1-\phi(\cdot))$.~Letting the parameter $p_k$ be (i) the $l$-th element of the output weight of the $k$-th SLFNN output, $w_l^{o(k)}$, or (ii) the $k$-th element of the output bias, $b^{o(k)}$, or (iii) the $l$-the element of the internal weight of the $h$-th input corresponding to the $k$-th SLFNN output, $w_{l,h}^{(k)}$, or (iv) the $l$-th element of the input bias, $b_{l}^{(k)}$, the related derivatives for $k=n$ (they are zero $k\neq n$) can be computed as:
\begin{align}
   \dfrac{\partial \tilde{\pi}^{NN}_n(\mathbf{y}^{(q)},\mathbf{p}_n)}{\partial w_l^{o(k)}} & = -\phi_l(\cdot), & \dfrac{\partial \tilde{\pi}^{NN}_n(\mathbf{y}^{(q)},\mathbf{p}_n)}{\partial w_{lh}^{(k)}} & = -w_l^{o(k)} y^{(q)}_h \phi_l(\cdot)(1-\phi_l(\cdot)),  \nonumber  \\ 
   \dfrac{\partial \tilde{\pi}^{NN}_n(\mathbf{y}^{(q)},\mathbf{p}_n)}{\partial b^{o(k)}} & = -1,  & 
    \dfrac{\partial \tilde{\pi}^{NN}_n(\mathbf{y}^{(q)},\mathbf{p}_n)}{\partial b_{l}^{(k)}} &= - w_l^{o(r)} \phi_l(\cdot)(1-\phi_l(\cdot)),
\label{eq:derNN_pars}
\end{align}
where $\phi_l(\cdot) = \phi_l(\mathbf{w}_{l}^{(n)\top} \mathbf{y}^{(q)}  + b^{(n)}_l)$,  $y_h^{(q)}$ denotes the $h$-th element of the $q$-th collocation point, for $q=1,\ldots,Q$, $n=1,\ldots,N$, $l=1,\ldots,L$ and $h=1,\ldots,M$.~\Cref{eq:derNN_pars} provides all the necessary derivatives for the computation of the residual derivatives in \cref{eq:resNN_ders1,eq:resNN_ders2}, which can be then computed for the formation of the $N(Q+1) \times N(L(M+2)+1)$ Jacobian matrix $\mathbf{J}=\nabla_{\mathbf{p}} \boldsymbol{\mathcal{F}}^{NN}$.~The iterative procedure of the LM algorithm for the minimization of the above residuals is described in \cref{app:LMalg}.

For the standalone NN scheme, we provide approximation error and convergence guarantees in \cref{th:NN_NFEs}.

\renewcommand{\theequation}{C.\arabic{equation}}
\renewcommand{\thefigure}{C.\arabic{figure}}
\renewcommand{\thealgorithm}{C.\arabic{algorithm}}
\setcounter{equation}{0}
\setcounter{figure}{0}
\setcounter{algorithm}{0}
\section{The Levenberg-Marquardt algorithm for the non-linear residuals minimization}
\label[appendix]{app:LMalg}
Here, we briefly present the Levenberg-Marquardt algorithm \cite{hagan1994training}, a gradient-based optimization algorithm, for the solution of non-linear least squares problems.~As already discussed, we employ the LM algorithm to solve the optimization problem in \cref{eq:min_hybrid} resulting in the IM approximation provided by the hybrid and the standalone NN schemes.

\begin{algorithm}[!h] 
\footnotesize
\caption{Solution of the physics-informed schemes via the LM algorithm} 
\label{alg:PI_LM}
\begin{algorithmic}[1]
\REQUIRE Discrete dynamical system in \cref{eq:genDDS2} satisfying the assumptions of \cref{th:th1Kaz} \COMMENT{$\mathbf{A}$, $\mathbf{B}$, $\mathbf{C}$, $\mathbf{f}(\mathbf{x},\mathbf{y})$ and $\mathbf{g}(\mathbf{y})$}
\ENSURE Collocation points $\mathbf{y}^{(q)}\in \Omega$ for $q=1,\ldots,Q$ \COMMENT{Ensure $\Omega\subset \mathcal{V}$ where IM functional of \cref{eq:IM} exists}
\STATE Set type of IM approximation, $type \gets HxS$ or $NN$ \COMMENT{Hybrid scheme or NNs}
\item[] \textcolor{gray}{$\triangleright$~\textit{Set up the physics-informed loss function for the selected IM functional approximation}}
\IF{$type=HxS$}
    \STATE Set polynomial type, $polyType\gets PS,~LS$ or $CS$ \COMMENT{Power, Legendre or Chebyshev series}
    \STATE Set hyperparameters $\mathcal{H}_1$; degree $h$, polynomial expressions of $polyType$
        \STATE Set radius $\mathbf{r}$ of $\mathcal{D}$, and sample collocation points $\mathbf{y}^{(r)}\in \partial \mathcal{D}$ for $r=1,\ldots,R$
    \STATE Initialize the coefficients $\boldsymbol{a}_n$ of the $polyType$ series in \cref{eq:IM_polySexp}
    \COMMENT{See \cref{sbsb:train}}   
\ENDIF
    \STATE Set hyperparameters $\mathcal{H}_2$; number of layers $L$, activation function, etc
    \STATE Initialize the NN parameters $\mathbf{p}_n$ in \cref{eq:IM_NN} \COMMENT{See \cref{sbsb:train}}  
\STATE Set $\omega_{\Omega},\omega_0,\omega_{\partial \mathcal{D}}$ for balancing the residual terms in the loss function \COMMENT{See \cref{eq:min_hybrid}}
\item[] \textcolor{gray}{$\triangleright$~\textit{Employ the LM algorithm for solving the related physics-informed minimization problems}}
\IF{$type=HxS$} 
    \STATE Construct parameter vector $\boldsymbol{\nu}$ and set $\mathbf{R}_0\gets \boldsymbol{\nu}$ \COMMENT{\cref{eq:pars_PIHxS}}  
\ELSIF{$type=NN$}     
    \STATE Construct parameter vector $\mathbf{p}$ and set $\mathbf{R}_0\gets \mathbf{p}$ \COMMENT{\cref{eq:pars_PINN}}
\ENDIF
\STATE Set $k=0$, initial damping factor $\lambda_0$, tolerances $tol_F$, $tol_R$ and maximum iterations $k_{max}$
\item[] \textcolor{gray}{$\triangleright$~\textit{LM algorithm iterations}}
\REPEAT 
\IF{$type=HxS$}
\STATE Form the residual vector $\boldsymbol{\mathcal{F}}(\mathbf{R}_k)$ \COMMENT{$\boldsymbol{\mathcal{F}}^{HxS}(\boldsymbol{\nu})$ in \cref{eq:res_HS}}
\STATE Form the Jacobian matrix $\mathbf{J}(\mathbf{R}_k)$ \COMMENT{$\nabla_{\boldsymbol{\nu}} \boldsymbol{\mathcal{F}}^{HxS}(\boldsymbol{\nu})$ computed via \cref{eq:resHxS_ders1,eq:resHxS_ders2,eq:resHxS_ders3}}
\ELSIF{$type=NN$}
\STATE Form the residual vector $\boldsymbol{\mathcal{F}}(\mathbf{R}_k)$ \COMMENT{$\boldsymbol{\mathcal{F}}^{NN}(\mathbf{p})$ in \cref{eq:res_NN}}
\STATE Form the Jacobian matrix $\mathbf{J}(\mathbf{R}_k)$ \COMMENT{$\nabla_{\mathbf{p}} \boldsymbol{\mathcal{F}}^{NN}(\mathbf{p})$ computed via \cref{eq:resNN_ders1,eq:resNN_ders2}}
\ENDIF
\STATE Estimate Hessian matrix and find search direction $\mathbf{d}_k$ \COMMENT{Solve \cref{eq:LMsd}}
\STATE \textcolor{gray}{$\triangleright$~\textit{Perform the LM algorithm update}}
\STATE Compute the residuals at potential next iteration $\boldsymbol{\mathcal{F}} (\mathbf{R}_k+\mathbf{d}_k)$
\IF{$\lVert \boldsymbol{\mathcal{F}} (\mathbf{R}_k+\mathbf{d}_k)\rVert_{2}  <\lVert\boldsymbol{\mathcal{F}}(\mathbf{R}_k)\rVert_{2}$}  
\STATE Set $\mathbf{R}_{k+1} \gets \mathbf{R}_k+\mathbf{d}_k$ and $\lambda_{k+1} \gets \lambda_{k}/10$ \COMMENT{Reducing residuals}
\ELSE 
\STATE Set $\mathbf{R}_{k+1} \gets \mathbf{R}_k$ and $\lambda_{k+1} \gets 10\lambda_{k}$ \COMMENT{Increasing residuals}
\ENDIF
\STATE $k\gets k+1$
\UNTIL{$\lVert \boldsymbol{\mathcal{F}}(\mathbf{R}_{k+1}) - \boldsymbol{\mathcal{F}}(\mathbf{R}_{k}) \rVert_2<tol_F (1+\lVert\boldsymbol{\mathcal{F}}(\mathbf{R}_{k})\rVert_{2})$, or $\lVert \mathbf{R}_{k+1} - \mathbf{R}_{k} \rVert_{2}<tol_R (1+\lVert \mathbf{R}_{k}\rVert_{2})$, or $k\ge k_{max}$}
\RETURN tuned parameters $\mathbf{R}_k$ ($\boldsymbol{\nu}$ for $type=HxS$ and $\mathbf{p}$ for $type=NN$ 
)
\end{algorithmic}
\end{algorithm}

The aforementioned optimization problem is reformulated in \cref{app:OptPIHS_LM} into non-linear residual minimization problems.~The  particular form of the latter is used by the LM algorithm to minimize a residual vector $\boldsymbol{\mathcal{F}}(\mathbf{R})\in\mathbb{R}^{\zeta}$ w.r.t. to a vector of parameters $\mathbf{R}\in\mathbb{R}^{\eta}$ using the Jacobian matrix $\mathbf{J}(\mathbf{R})=\nabla_{\mathbf{R}}\boldsymbol{\mathcal{F}(\mathbf{R})}\in\mathbb{R}^{\zeta\times \eta}$; hereby, the dependence of the residuals and the Jacobian matrix on the parameters is explicitly denoted.~In correspondence to the above unifying notation, we obtain:
\begin{itemize}
    \item for the hybrid scheme minimization problem in \cref{sb:hyb_scheme_ders}, the parameter vector $\mathbf{R}\equiv\boldsymbol{\nu}\in \mathbb{R}^\eta$ in \cref{eq:pars_PIHxS} where $\eta=N\left( L(M+2) + \binom{M+h}{h} +1 \right)$, the residual vector $\boldsymbol{\mathcal{F}}(\mathbf{R}) \equiv \boldsymbol{\mathcal{F}}^{HxS}(\boldsymbol{\nu})\in\mathbb{R}^\zeta$ in \cref{eq:res_HS} where $\zeta=N(Q+R+1)$, and the Jacobian matrix $\mathbf{J}(\mathbf{R})\equiv\mathbf{J}(\boldsymbol{\nu})=\nabla_{\boldsymbol{\nu}} \boldsymbol{\mathcal{F}}^{HxS}(\boldsymbol{\nu})$ which can be computed by \cref{eq:resHxS_ders1,eq:resHxS_ders2,eq:resHxS_ders3}.  
    \item for the SLFNNs minimization problem in \cref{sb:NN_scheme_ders}, the parameter vector $\mathbf{R}\equiv\mathbf{p}\in \mathbb{R}^\eta$ in \cref{eq:pars_PINN} where $\eta=N(L(M+2)+1)$, the residual vector $\boldsymbol{\mathcal{F}}(\mathbf{R}) \equiv \boldsymbol{\mathcal{F}}^{NN}(\mathbf{p})\in\mathbb{R}^\zeta$ in \cref{eq:res_NN} where $\zeta=N(Q+1)$, and the Jacobian matrix $\mathbf{J}(\mathbf{R})\equiv\mathbf{J}(\mathbf{p})=\nabla_{\mathbf{p}} \boldsymbol{\mathcal{F}}^{NN}(\mathbf{p})$ which can be computed by \cref{eq:resNN_ders1,eq:resNN_ders2}.  
\end{itemize}
For the implementation of the LM iterative algorithm, an initial guess of parameters, say $\mathbf{R}_0$, and an initial damping factor, say $\lambda_0$, are first set; see \cref{sbsb:train} for the initialization.~Then, at the $k$-th iteration, the residual vector $\boldsymbol{\mathcal{F}}(\mathbf{R}_k)$ and the Jacobian matrix $\mathbf{J}(\mathbf{R}_k)$ are computed on the basis of the parameter values $\mathbf{R}_k$.~Having constructed the Jacobian matrix, the LM algorithm next estimates the Hessian matrix of the $k$-th iteration by $\mathbf{J}(\mathbf{R}_k)^\top\mathbf{J}(\mathbf{R}_k)$, which is used for determining the search direction $\mathbf{d}_k\in\mathbb{R}^\eta$ from the solution of the damped linear system:
\begin{equation}
\left( \mathbf{J}(\mathbf{R}_k)^\top\mathbf{J}(\mathbf{R}_k)  + \lambda_k \mathbf{I}_\eta \right) \mathbf{d}_k = - \mathbf{J}(\mathbf{R}_k)^\top \boldsymbol{\mathcal{F}}(\mathbf{R}_k),
    \label{eq:LMsd}
\end{equation}
where $\lambda_k$ is the damping factor of the $k$-th iteration, and $\mathbf{I}_\eta$ is the $\eta$-dim. unitary matrix.~Finally, the parameters $\mathbf{R}_{k+1}$ and the damping factor $\lambda_{k+1}$ at the next iteration are updated according to the following (success or fail) condition:
\begin{itemize}
\item if the parameters to be updated are reducing the residuals, then updated them and decrease the damping factor, that is if $\lVert \boldsymbol{\mathcal{F}} (\mathbf{R}_k+\mathbf{d}_k)\rVert_{2}<\lVert\boldsymbol{\mathcal{F}}(\mathbf{R}_k)\rVert_{2}$, then $\mathbf{R}_{k+1} = \mathbf{R}_k+\mathbf{d}_k$ and $\lambda_{k+1}=\lambda_{k}/10$, or
\item if the parameters to be updated are not reducing the residuals, then do not update them, but instead increase the damping factor, that is if $\lVert \boldsymbol{\mathcal{F}}(\mathbf{R}_k+\mathbf{d}_k)\rVert_2\ge\lVert \boldsymbol{\mathcal{F}}(\mathbf{R}_k)\rVert_2$, then $\mathbf{R}_{k+1} = \mathbf{R}_k$ and $\lambda_{k+1}=10\lambda_k$.
\end{itemize}
The stopping criteria for the termination of the LM algorithm are based on relative convergence of the residuals, relative step tolerance and maximum iteration limits; i.e., the LM algorithm stops when either $\lVert \boldsymbol{\mathcal{F}}(\mathbf{R}_{k+1}) - \boldsymbol{\mathcal{F}}(\mathbf{R}_{k}) \rVert_2<tol_F (1+\lVert\boldsymbol{\mathcal{F}}(\mathbf{R}_{k})\rVert_{2})$, or $\lVert \mathbf{R}_{k+1} - \mathbf{R}_{k} \rVert_{2}<tol_R (1+\lVert \mathbf{R}_{k}\rVert_{2})$, or $k\ge k_{max}$.~For the selection of $tol_F$, $tol_R$ and $k_{max}$, see \cref{sbsb:train}.~A detailed pseudocode of the LM iterative algorithm is provided in \cref{alg:PI_LM}, which is focused on the minimization problems resulting from both the hybrid and the NNs schemes.

\renewcommand{\theequation}{D.\arabic{equation}}
\renewcommand{\thefigure}{D.\arabic{figure}}
\setcounter{equation}{0}
\setcounter{figure}{0}
\section{Polynomial Regression Accuracy: Linear vs. Non-linear Optimization}
\label[appendix]{app:reg_LvsNL}
This example demonstrates how non-linear optimization techniques, such as the LM algorithm, introduce numerical errors in polynomial regression compared to linear least-squares methods like the Moore-Penrose pseudo-inverse, particularly for high-degree approximations.~Consider the Gaussian profile:
\begin{equation}
    f(x) = 1-exp(-10x^2),  \label{eq:GaussEx}
\end{equation}
within the region $x\in[-0.3,0.3]$ around $f(0)=0$.~Therein, we seek a power series approximation of degree $h$ in the form of \cref{eq:IM_polySexp}, expressed as $\pi(x,\boldsymbol{a};\mathcal{H}_1) =  \sum_{i=0}^h \alpha_{i} x^i$.~Here $\boldsymbol{a}=[\alpha_1,\ldots,\alpha_h]^\top$ represents the power series coefficients and $\mathcal{H}_1$ denotes the hyperparameters, including the degree $h$ of the power series.

Given a set of features $x^{(q)}\in [-0.3,0.3]$ and their corresponding target values $y^{(q)}=f(x^{(q)})$ for $q=1,\ldots,Q=200$, we solve the regression optimization problem: 
\begin{equation}
    \argmin_{\boldsymbol{a}} \mathcal{L}(\boldsymbol{a};\mathcal{H}_1) = \sum_{q=1}^Q\lVert y^{(q)} - \pi(x^{(q)},\boldsymbol{a};\mathcal{H}_1)\rVert^2. \label{eq:PolyRegOptim}
\end{equation}
to determine the coefficients $\boldsymbol{a}$.~This least-squares problem can be approached using both non-linear and linear techniques.~For the non-linear solution, we employ the LM algorithm, following the iterative procedure described in \cref{app:LMalg}.~The residuals $\mathcal{F}_q^{PS}$ and their derivatives w.r.t. the power series coefficients are given by the expressions:
\begin{equation}
\mathcal{F}_q^{PS}=y^{(q)}-\pi(x^{(q)},\boldsymbol{a};\mathcal{H}_1), \qquad \dfrac{\partial \mathcal{F}_q^{PS}}{\partial \alpha_i} = - \dfrac{\pi(x^{(q)},\boldsymbol{a};\mathcal{H}_1)}{\partial \alpha_i} = -\left(x^{(q)} \right)^i.
\end{equation}

Alternatively, the problem in \cref{eq:PolyRegOptim} can be formulated as a linear system $\mathbf{V} \boldsymbol{a} = \mathbf{y}$, where $\mathbf{V}\in\mathbb{R}^{Q\times h}$ is the so-called Vandermonde matrix and $\mathbf{y}=[y^{(1)},\ldots,y^{(Q)}]^\top$.~The solution is then obtained via the Moore-Penrose pseudo-inverse as $\boldsymbol{a} = \mathbf{V}^\dagger \mathbf{y}$.

\begin{figure}[!h]
\centering
    \subfigure[Degree $h=10$]{\includegraphics[width=0.47\textwidth]{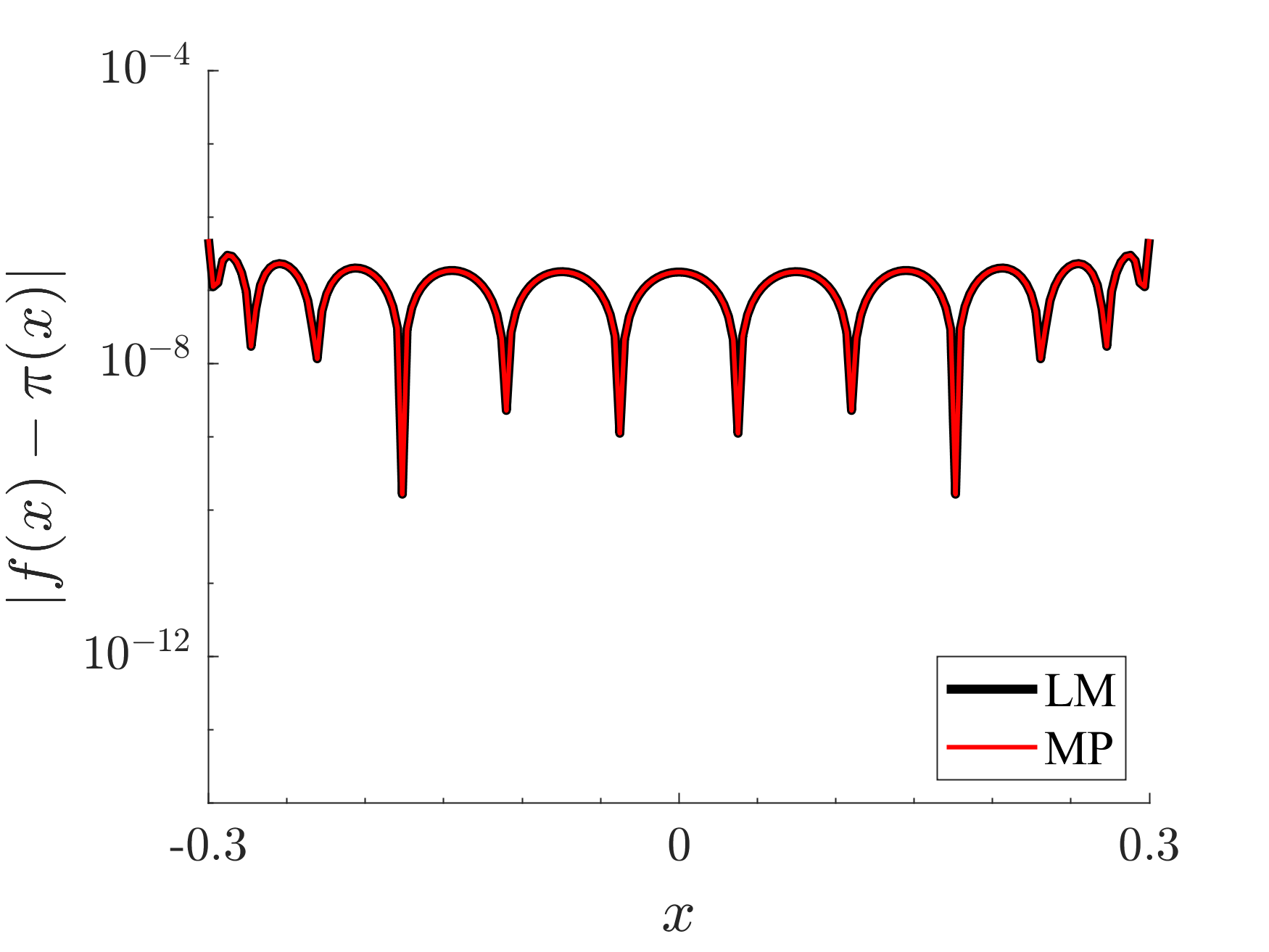}} \hspace{10pt}
    \subfigure[Degree $h=20$]{\includegraphics[width=0.47\textwidth]{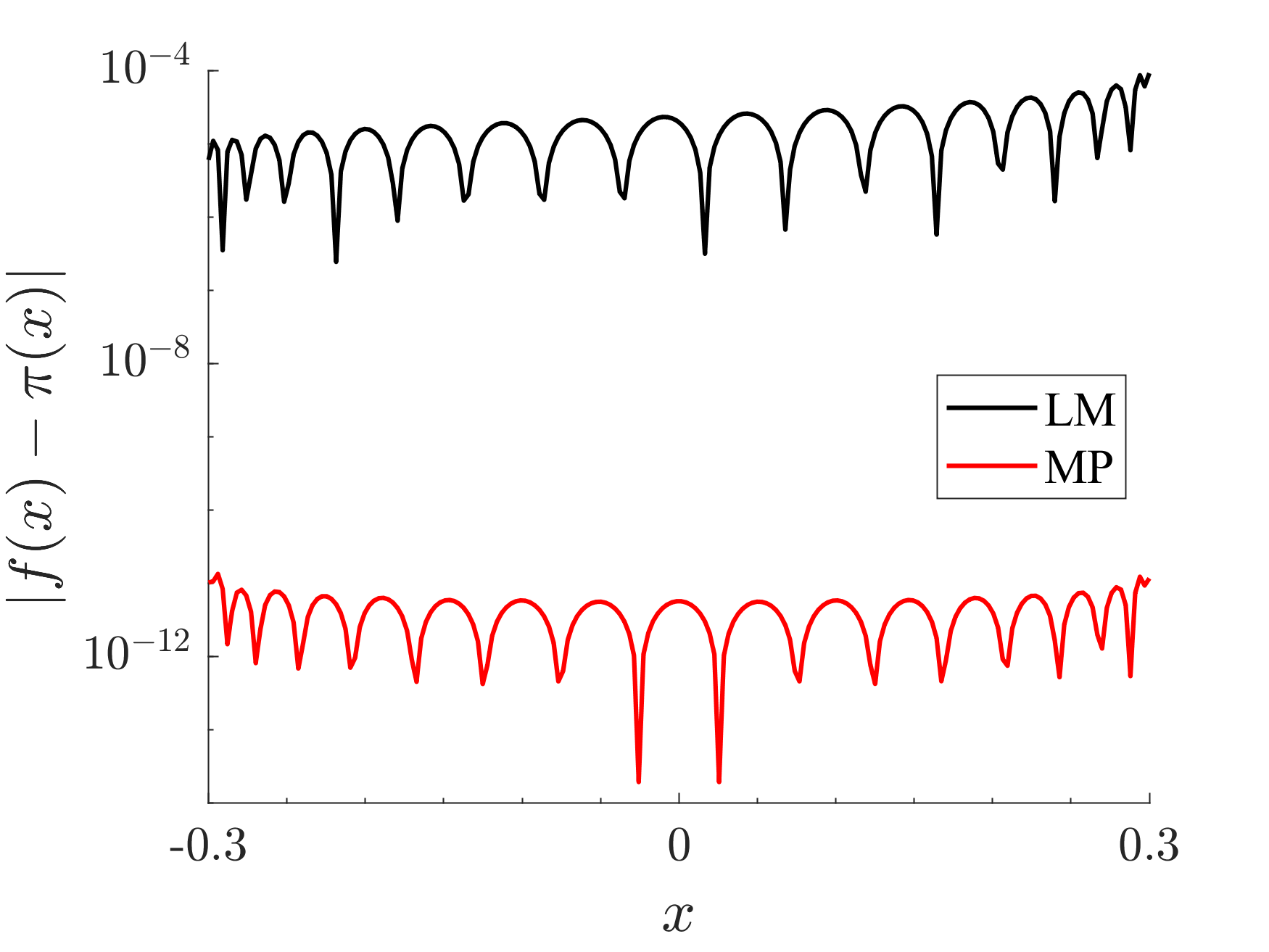}} 
\caption{Absolute errors of the power series approximations $\pi(x)$ of the function $f(x)$ in \cref{eq:GaussEx}, for polynomials of degree $h=10$ and $h=20$.~LM and MP denote the use of the Levenberg-Marquardt non-linear optimization algorithm and the Moore-Penrose pseudo-inverse linear least-squares method for obtaining the power series coefficients.}
\label{fig:regGauss}
\end{figure}

We solved \cref{eq:PolyRegOptim} using both non-linear and linear methods for polynomial degrees $h=10$ and $h=20$.~\Cref{fig:regGauss} compares the absolute errors of the resulting approximations.~For $h=10$, both methods achieve comparable approximation accuracy.
However, at $h=20$, the non-linear optimization introduces numerical errors, leading to inferior approximation compared to the linear method.~This observation aligns with \cref{rm:2}, underscoring the advantage of linear techniques for high-degree polynomial regression in our hybrid scheme.

\renewcommand{\theequation}{E.\arabic{equation}}
\setcounter{equation}{0}
\section{Implementation of the physics-informed approaches for the IM approximation}
\label[appendix]{app:Imp}
This appendix provides all the necessary details for implementing the hybrid scheme introduced in \cref{sb:PI_opt_hyb}, along with the standalone NNs scheme.~Specifically, we describe the training process and discuss the evaluation of the numerical accuracy of the learned IM approximations.~An overview of the training and evaluation process is presented in \cref{alg:IM_imp}, the steps of which are explained in detail in the following paragraphs.

\begin{algorithm}[!h] 
\footnotesize
\caption{Outline of training and evaluation of the IM functionals approximated via the PI hybrid and NN schemes; bold comments denote paragraphs of \cref{app:Imp} where the specific step is discussed.} 
\label{alg:IM_imp}
\begin{algorithmic}[1]
\REQUIRE Discrete system in the form of \cref{eq:genDDS1} 
\ENSURE Assumptions of \cref{th:th1Kaz} around equilibrium \COMMENT{$\mathbf{A}$, $\mathbf{B}$, $\mathbf{C}$, $\mathbf{f}(\mathbf{x},\mathbf{y})$ and $\mathbf{g}(\mathbf{y})$ of \cref{eq:genDDS2}}
\STATE Define $\Omega \subset \mathcal{V}$ to approximate the IM functional in \cref{eq:IM}
\STATE Sample collocation points $\mathbf{y}^{(q)}\in \Omega$ \COMMENT{See \textbf{Training data acquisition}}
\STATE Sample testing data points $(\mathbf{x}^{(s)},\mathbf{y}^{(s)}) \in \mathcal{M}$ with $\mathbf{y}^{(s)}\in \Omega$  \COMMENT{See \textbf{Testing data acquisition}}
\STATE Set type of IM approximation (hybrid or pure NN) and select architectures and hyperparameters $\mathcal{H}_1$, $\mathcal{H}_2$, and/or $\mathbf{r}$ \COMMENT{See \textbf{Architectures and Hyperparameters}}
\STATE Set collocation points $\mathbf{y}^{(r)}\in \partial \mathcal{D}$ for $r=1,\ldots,R$
\FOR{multiple training realizations}
\STATE Initialize parameters $\boldsymbol{a}$ and $\mathbf{p}$ \COMMENT{See \textbf{Parsimonious initialization of parameters}}
\item[] \textcolor{gray}{$\triangleright$~\textit{Solve the related PI optimization problem using the LM algorithm (see \cref{app:LMalg})}}
\REPEAT 
\STATE Form non-linear residuals of loss function \COMMENT{\cref{eq:min_hybrid}}
\STATE Form Jacobian matrix of the residuals w.r.t. the parameters
\STATE Perform the LM iteration to update the parameters and the damping factor
\UNTIL{Convergence or maximum epochs reached}
\STATE Compute training metrics on collocation points $\mathbf{y}^{(q)}$ \COMMENT{See \textbf{Training metrics}}
\STATE Compute testing metrics on testing points $(\mathbf{x}^{(s)},\mathbf{y}^{(s)})$\COMMENT{See \textbf{Numerical accuracy metrics}} 
\ENDFOR
\STATE Average over multiple realizations and compute percentiles in training and testing metrics
\STATE Visualize numerical accuracy with the best-performing (on testing)  parameter set
\end{algorithmic}
\end{algorithm}

\subsection{Training process}
\label[appendix]{sbsb:train}

The training process was conducted over 100 training realizations to quantify the uncertainty in training performance.~Below, we outline its key components: training data acquisition,  the architectures and hyperparameters selected, the initialization of parameters, and the metrics used to evaluate convergence.~These components correspond to steps 2, 4, 7, and 13 of \cref{alg:IM_imp}, respectively.~The minimization of non-linear residuals via the LM algorithm is presented in \cref{app:LMalg}. 

\paragraph{Training data acquisition.} To learn the IM approximations using the proposed hybrid scheme in \cref{sb:PI_opt_hyb} (or the standalone NN scheme) a set of $Q$ collocation points $\mathbf{y}^{(q)}$ for $q=1,\ldots,Q$ is required.~Since the IM $\mathcal{M}$ in \cref{eq:IM} is defined in a local neighborhood $\mathcal{V}\subset\mathbb{R}^M$ around the equilibrium $\mathbf{y}_0$ (see \cref{th:th1Kaz}), the collocation points are collected within a desired domain $\Omega \subset \mathcal{V}$, chosen to avoid singular points.

For $M=1$ cases, one can select an interval $\Omega=[a_1,b_1]$ containing $y_0=0$.~However, the generalization of $\Omega$ in higher dimensions is not straightforward.~For example, there is no guarantee that the IM exists for all $\mathbf{y}$ in the hyperrectangle $\Omega = [a_1,b_1]\times\cdots\times[a_M,b_M]$.~To address this, for $M>1$, we collected collocation points from numerically derived trajectories, following the approach in \cite{patsatzis2024slow,patsatzis2024slowB}.~Specifically, we selected $n_{IC}$ initial conditions randomly chosen outside a predefined hyperrectangle $\Omega$ and generated trajectories of the system in \cref{eq:genDDS1}.~After discarding a transient of $k_{trans}$ time steps (to ensure the dynamics have approached the IM), we recorded all subsequent time steps until reaching the equilibrium, using a cutoff tolerance of $0.001$ for each variable to avoid oversampling near the equilibrium.~From the collected points in the $M>1$ cases, we applied an arc-length parameterization to each trajectory and selected points at equal intervals along the cumulative arc length, to ensure more uniform sampling than random sampling in $\Omega$ and to prevent oversampling near the equilibrium due to stiffness.~After verifying $\mathbf{y}\in \Omega$ for the selected points, we retained $Q$ of them to form the training dataset $\mathbf{y}^{(q)}$ for $q=1,\ldots,Q$.~For the $M=1$-dim.~cases, we simply sampled $Q$ uniform points from $\mathcal{U}([a_1,b_1])$.

To enable fair comparison of computational training times between the hybrid and NN schemes, we chose the number of collocation points $Q$ such that the number of residuals is 20 times the number of unknown parameters of the--common between the two schemes--NN component.

In addition, for the continuity term in the loss function in \cref{eq:min_hybrid}, it is further required to collect $R$ collocation points on the boundary $\partial \mathcal{D}$.~For the $M=1$-dim. case, we simply select the two points $\mathbf{y}^{(r)}=\pm\mathbf{r}$.~For higher dimensions, we uniformly sample $5^{M-1}$ points along each hyperplane of $\partial \mathcal{D}=\{\mathbf{y}\in\mathbb{R}^M: \lvert y_m\rvert =r_m, \text{ for some } m=1,\ldots,M; \lvert y_j\rvert \leq r_j, \forall j\neq m\}$; e.g., for $M=2$, we sample $R=4\cdot5=20$ points across the four edges of the rectangle $\partial \mathcal{D}$, while for $M=3$, we sample $R=6\cdot 5^2=150$ points across the six surfaces of the cube $\partial \mathcal{D}$.

\paragraph{Architectures and Hyperparameters.} For the polynomial series components of the hybrid scheme presented in \cref{eq:IM_polySexp}, we considered power series, Legendre and Chebyshev (of the second kind) polynomials, the degree $h$ of which varied across the examples studied, depending on the desired accuracy.~For the NN component, we employed the SLFNN architecture presented in \cref{eq:IM_NN}.~Increasing the depth to two hidden layers did not yield significant improvements in accuracy, so we used a single hidden layer with $L$ neurons (varying per case study) and a logistic sigmoid activation function.

An important hyperparameter for the PI hybrid scheme is the radius $\mathbf{r}$ of the hyperrectangle $\mathcal{D}$.~Each element $r_m$ of $\mathbf{r}$ is free to vary within the interval $[0, max(\lvert a_m\rvert,\lvert b_m\rvert)]$, where $[a_m,b_m]$ is the range of the $y_m$ variable in the training set $\Omega$.~For Legendre and Chebyshev polynomial series, the restriction $r_m\le1$ applies to ensure the polynomials are well-defined.~The choice of $r_m$ determines the hyperrectangle $\mathcal{D}$ where the polynomial series counterpart is active.~Since no prior information about the underlying IM functional is generally available, $r_m$ requires careful tuning.~In all the examples, we tuned $r_m$ starting from values leaving approximately half of the collocation points outside $\mathcal{D}$.~While it is out of the scope of this work, we emphasize that further analysis is needed to determine optimal values of $r_m$, taking into account the polynomial degree $h$.

For the optimization problem, we determined the three balancing terms $\omega_{\Omega},\omega_0,\omega_{\partial \mathcal{D}}$ in \cref{eq:min_hybrid} using residual balancing \cite{wang2021understanding}.~For the implementation of the LM algorithm, we configured the initial damping factor to $\lambda_0=10^{-2}$ and the hyperparameters related to the stopping criteria to a maximum number of epochs $k_{max}=1000$, a relative function tolerance $tol_F=10^{-8}$ and a relative step tolerance $tol_R=10^{-4}$; see \cref{app:LMalg} for further details.

\paragraph{Parsimonious initialization of parameters.} The solution of the optimization problem in \cref{eq:min_hybrid} requires an initial guess of the learnable parameters: the coefficients $\boldsymbol{a}$ of the polynomial series and the parameters $\mathbf{p}$ of the NN.~While a random initial guess can be selected, we observed that such an approach leads to poor convergence for approximately 20\% of the training realizations.~This behavior arises in high-degree polynomial series, because random initialization can result in high-order terms dominating the loss function, causing the optimizer to converge to suboptimal local minima of the non-linear loss function.

To address this, we employ a parsimonious initialization strategy for the polynomial coefficients $\boldsymbol{a}$, ensuring that each monomial is $\mathcal{O}(1)$.~In particular, we initialize the non-zero coefficients $\alpha_n^{i_1,\ldots,i_M}$ in \cref{eq:IM_polySexp} by sampling random values uniformly distributed in the interval $[-1/c_n^{i_1,\ldots,i_M},~1/c_n^{i_1,\ldots,i_M}]$, where 
\begin{equation*}
    c_n^{i_1,\ldots,i_M} = \max_{1\le q \le Q}  \{ \big \lvert P^x_{i_1}(y_1^{(q)}) \cdots P^x_{i_M}(y_M^{(q)}) \big \rvert \}, 
\end{equation*}
where $i_1,\ldots,i_M=1,\ldots,h$ restricted by $i_1 + \ldots + i_M \le h$ and $y_m^{(q)}$ denotes the $m$-th element of the $q$-th collocation point $\mathbf{y}^{(q)}$.~For the NN parameters $\mathbf{p}$, we use a uniform Xavier/Glorot initialization \cite{glorot2010understanding}.~To further avoid numerically exploding or vanishing gradients, we ensure that $ \mathbf{W}^{(n)} \mathbf{y}  + \mathbf{b}^{(n)}=\mathcal{O}(1)$ by further scaling the Xavier-initialized input weights $w_{l,m}^{(n)}$ by $max_q(y_m^{(q)})-min_q(y_m^{(q)})$.~The parameters $\boldsymbol{a}$ and $\mathbf{p}$, initialized using the above procedure, consistently demonstrated good convergence by improving the conditioning of the optimization problem and enabling the LM algorithm to perform reliably.

\paragraph{Training metrics.} To evaluate the convergence of the optimization problem in \cref{eq:min_hybrid}, we record the values of the respective loss functions at the end of training for each of the 100 realizations of different initial parameters, along with the computational times required for training.~The value of each loss function corresponds to the squared $l_2$ error, $\lVert \boldsymbol{\mathcal{F}} \rVert_2$, of the non-linear residuals $\boldsymbol{\mathcal{F}}$ minimized through the LM algorithm.~To quantify the uncertainty and assess the consistency of convergence, we compute the mean and 5-95\% percentiles of these training metrics across the 100 realizations.

\subsection{Evaluation of numerical approximation accuracy}
\label[appendix]{sbsb:numAc}
For assessing the numerical approximation accuracy of the IM approximations obtained by the hybrid and standalone NNs scheme, we compare data points that lie exclusively on the true IM $\mathcal{M}$ in \cref{eq:IM} to their corresponding approximations.~Specifically, for a given $\mathbf{y} \in \Omega$, the true IM mapping $\boldsymbol{\pi}(\mathbf{y})$ provides the corresponding $\mathbf{x}$ (i.e., $\mathbf{x}=\boldsymbol{\pi}(\mathbf{y})$), while the IM functional approximation yields the estimate $\tilde{\mathbf{x}}=\tilde{\boldsymbol{\pi}}(\mathbf{y})$.~By comparing $\mathbf{x}$ and $\tilde{\mathbf{x}}$, we quantify the accuracy of the IM approximation.~In what follows, we discuss the construction of the testing data sets and the error metrics used for quantifying numerical accuracy.

\paragraph{Testing data acquisition.} As previously discussed, the construction of the testing data set requires data points that lie exclusively on the underlying IM $\mathcal{M}$ in \cref{eq:IM}.~However, an explicit expression for the true IM mapping is generally unknown for dynamical systems of the form \cref{eq:genDDS1}.~To address this, we collected testing data from numerically derived trajectories.~Following the same procedure as for acquiring the training data, we selected $n_{IC}$ random initial conditions outside $\Omega$ within example-specific  bounds, generated trajectories, discarded a transient of $k_{trans}$ time steps and finally recorded subsequent time steps until reaching the equilibrium, using a cutoff of $0.001$ for each variable.~We randomly selected $S=10,000$ from the recorded points and used both $\mathbf{y}^{(s)}$ and $\mathbf{x}^{(s)}$ to form the testing sets, $(\mathbf{x}^{(s)},\mathbf{y}^{(s)})\in\mathcal{M}$, ensuring $\mathbf{y}^{(s)}\in \Omega$ for $s=1,\ldots,S$.~This construction procedure was followed for every example, except for the one in \cref{app:lnEx}, where the IM mapping $\pi(y) = \ln(1+y)$ was analytically known and used for generate the testing data.

\paragraph{Numerical accuracy metrics.} To assess the numerical approximation accuracy, we evaluated the error between $\mathbf{x}^{(s)}$ and $\tilde{\mathbf{x}}^{(s)}=\tilde{\boldsymbol{\pi}}(\mathbf{y}^{(s)})$ on the testing data set.~Specifically, we computed the $L^1$, $L^2$ and $L^\infty$ norms of the relative errors $\lVert x^{(s)}-\tilde{\pi}(y^{(s)})\rVert/ \lVert x^{(s)}\rVert$ across all testing points $s=1,\ldots,S$ for the IM approximations.~For the examples with dimension higher than $N=1$, we aggregated these metrics by taking the mean over all $n=1,\ldots,N$ components.~To quantify the numerical accuracy, we report the mean and 5-95\% percentiles of the above metrics across the 100 parameter sets learned during training.

For visualization purposes, we select the parameter set with the best $L^2$ norm error and display the relative errors $\lvert (\mathbf{x}_n^{(s)} - \tilde{\boldsymbol{\pi}}_n(\mathbf{y}^{(s)}))/ \mathbf{x}_n^{(s)} \rvert$ across all points to highlight regions of high and low numerical approximation accuracy.

\renewcommand{\theequation}{F.\arabic{equation}}
\renewcommand{\thefigure}{F.\arabic{figure}}
\setcounter{equation}{0}
\setcounter{figure}{0}
\section{Example with analytic invariant manifold}
\label[appendix]{app:lnEx}
In this appendix, we consider a non-linear discrete dynamical system with an exactly known IM to validate our PI hybrid and standalone NN schemes.~In particular, we consider the system in the form of \cref{eq:genDDS1} with $N=1$ and $M=1$, as:
\begin{align}
    x(k+1) & = \beta x(k) + y(k), \nonumber \\
    y(k+1) & = (1+y(k))^\beta e^{y(k)} - 1,
    \label{eq:TP_DDS1}
\end{align}
where $\beta \neq 0$ and $y\ge-1$.~For $\beta<0$, a singularity occurs at $y=-1$.~The equilibrium of the system is $(x_0,y_0)=(0,0)$ and \Cref{ass1} is satisfied for $\beta \neq -1$.~Under this assumption, we rewrite the system in the form of \cref{eq:genDDS2} as:
\begin{align}
    x(k+1) & = \beta x(k) + y(k), \nonumber \\
    y(k+1) & = (1+\beta) y(k) + (1+y(k))^\beta e^{y(k)} - (1+\beta) y(k) - 1, 
    \label{eq:TP_DDS2}
\end{align}
where all assumptions regarding the non-linear functions are satisfied at the equilibrium point.~According to \cref{th:th1Kaz} and since $\mathbf{A}=1+\beta$ and $\mathbf{B}=\beta$, an analytic local IM $\mathcal{M}$ exists around a neighborhood of the equilibrium in the form of \cref{eq:IM}, provided that the non-resonance condition $(1+\beta)^d\neq \beta$ holds for any $d\in \mathbb{N}$.~It can be shown that the IM mapping $\pi$, given by:
\begin{equation}
    x = \pi(y) = \ln(1+y),
    \label{eq:lnIM_anal}
\end{equation}
is the exact analytic solution of the NFEs in \cref{eq:NFEs_IM} for the system in \cref{eq:TP_DDS2}.~Here, we set $\beta=-0.4$, ensuring that the equilibrium is stable (trajectories are attracted to it) and the non-resonance condition is satisfied.~Note that, as previously mentioned and evident from \cref{eq:lnIM_anal}, a singularity exists at $y=-1$.

Although the exact  IM mapping is known for this example, we approximated it analytically using the PSE approach proposed in \cite{kazantzis2001invariant}, as described in \cref{app:IM_PS}. The resulting PSE for the IM approximation coincides with the well-known $h$-th degree Taylor expansion of the logarithmic function:
\begin{equation}
    \tilde{x} = \tilde{\pi}^{PSE}(y) =  \sum_{i=1}^h \alpha_{PSE}^i P_i^{P}(y),  \qquad \alpha_{PSE}^i = (-1)^{i+1}/i,
    \label{eq:lnIM_PSanal}
\end{equation}
where $P_i^{P}(y)=y^i$.~The radius of convergence of this power series is $\lvert y \rvert<r=1$.~Additionally, since the exact IM mapping is available in \cref{eq:lnIM_anal}, we  approximated it using Legendre and Chebychev polynomials, yielding:
\begin{equation}
    \tilde{x}  = \tilde{\pi}^{LSE}(y) =  \sum_{i=1}^h \alpha_{LSE}^i P_i^{L}(y), \qquad \alpha_{LSE}^i = \dfrac{2i+1}{2} \int_{-1}^1 \ln(1+y) P_i^{L}(y) dy,
    \label{eq:lnIM_LSanal}
\end{equation}
for the Legendre approximation, and 
\begin{equation}
\tilde{x}  = \tilde{\pi}^{CSE}(y) = \sum_{i=1}^h \alpha_{CSE}^i P_i^{C}(y),\qquad \alpha_{CSE}^i = \dfrac{2}{\pi} \int_{-1}^1 \ln(1+y) P_i^{C}(y) \dfrac{dy}{\sqrt{(1-y^2)}}.
    \label{eq:lnIM_CSanal}
\end{equation}
for the Chebyshev approximation.~The Legendre and Chebyshev polynomials, $P_i^{L}(y)$ and $P_i^{C}(y)$, are defined by the recursive formulas in \cref{eq:LS_rec,eq:CS_rec}, respectively.

In this example, we not only have the exact expression of the IM mapping in \cref{eq:lnIM_anal}, but also the analytic expressions of the power, Legendre and Chebyshev polynomial series in \cref{eq:lnIM_PSanal,eq:lnIM_LSanal,eq:lnIM_CSanal}, computed for degree $h=20$.~The former will be used next for evaluating the accuracy of the PI schemes, while the latter will serve as approximations for comparison purposes.

\subsection{Numerical results}
For the numerical approximation of the IM with the proposed PI hybrid and standalone NN schemes, we considered training data in the domain $\Omega=[-0.9,2]$, as the system exhibits a singularity at $y=-1$.~To examine the impact of polynomial degree, we considered the PI hybrid schemes with polynomials of degree $h=10$ (baseline) and $h=20$, and NNs with $L=10$ neurons in the hidden layer.~Following the procedure detailed in \cref{sbsb:train}, we uniformly sampled $Q=620$ collocation points along $\Omega$.~All the hyperparameters are set as described in \cref{app:Imp}.~We trained the PI hybrid schemes with power, Legendre and Chebyshev polynomial series (PI-HPS, PI-HLS, PI-HCS) and a radius of $r=1$, while we also considered a PI hybrid scheme with power series of radius $r=0.5$.~The convergence results of all nine trained PI approximations are provided in \cref{tb:conv_lnEx} of Supplement~\ref{supp1}, as obtained over 100 training realizations with different initial parameters; see \cref{sbsb:train} for details.

To evaluate the numerical approximation accuracy of the learned IM approximations, we constructed a testing data set based on the true IM mapping in \cref{eq:lnIM_anal}.~Following \cref{app:Imp}, we collected testing points $x^{(s)} = \pi(y^{(s)}) = \ln(1+y^{(s)})$ with $y^{(s)}$ uniformly sampled from $\Omega$, and computed the $L^1$, $L^2$, and $L^\infty$ norms of the relative errors $\lVert x^{(s)}-\tilde{\pi}(y^{(s)})\rVert/ \lVert x^{(s)}\rVert$.~\Cref{tb:acc_lnEx} reports the means and 5-95\% percentiles of these errors over the 100 training realizations, alongside the same error metrics for the power, Legendre and Chebyshev series expansions (PSE, LSE, CSE) derived from the true IM in \cref{eq:lnIM_PSanal,eq:lnIM_LSanal,eq:lnIM_CSanal}.
\begin{table}[!h]
    \centering
    \caption{Numerical approximation accuracy of the IM approximations $\tilde{\pi}(y)$ for the example in \cref{app:lnEx} on the test data.~Relative errors ($L^1$, $L^2$, and $L^\infty$ norms) are reported for the power, Legendre and Chebyshev series expansions (PSE, LSE, CSE) derived from the true IM in \cref{eq:lnIM_PSanal,eq:lnIM_LSanal,eq:lnIM_CSanal}, the standalone PINN and the PI hybrid schemes (PI-HxS with $x=P,L,C$ denoting power, Legendre, and Chebyshev polynomial series); for hyperparameters see \cref{tb:conv_lnEx}.~Mean values and 5–95\% percentiles are computed over the 100 parameter sets obtained during training for the PI schemes.~The testing data errors are evaluated over $S=10,000$ points.}
    \resizebox{\textwidth}{!}{
    \begin{tabular}{l l c c c c c c}
    &	\textbf{Scheme}	&	\multicolumn{2}{c}{$\lVert x^{(s)}-\tilde{\pi}(y^{(s)})\rVert_1/ \lVert x^{(s)}\rVert_1$}	&	\multicolumn{2}{c}{$\lVert x^{(s)}-\tilde{\pi}(y^{(s)})\rVert_2/ \lVert x^{(s)}\rVert_2$}	&	\multicolumn{2}{c}{$\lVert x^{(s)}-\tilde{\pi}(y^{(s)})\rVert_\infty/ \lVert x^{(s)}\rVert_\infty$}	\\[2pt]
    &		&	mean	&	5-95\% percentiles			&	mean	&	5-95\% percentiles			&	mean	&	5-95\% percentiles			\\
    \noalign{\smallskip}\hline\noalign{\smallskip}
    &	PINN	&	$2.61\times 10^{-5}$	& [$6.72\times 10^{-6}-7.11\times 10^{-5}$] &	$5.34\times 10^{-5}$	& [$1.44\times 10^{-5}-1.38\times 10^{-4}$] &	$2.32\times 10^{-4}$	& [$6.13\times 10^{-5}-6.01\times 10^{-4}$] \\
    \noalign{\smallskip}\hline\noalign{\smallskip}
    \multirow{7}{*}{\rotatebox[origin=c]{90}{$h=10$}} 
    &	PSE	&	$5.91\times 10^{0}$	&				&	$1.48\times 10^{1}$	&				&	$2.86\times 10^{1}$	&				\\
    &	LSE	&	$1.09\times 10^{3}$	&				&	$2.94\times 10^{3}$	&				&	$6.08\times 10^{3}$	&				\\
    &	CSE	&	$5.58\times 10^{2}$	&				&	$1.50\times 10^{3}$	&				&	$3.09\times 10^{3}$	&				\\
    &	PI-HPS, $r=1$	&	$5.77\times 10^{-4}$	& [$5.50\times 10^{-4}-6.07\times 10^{-4}$] &	$9.88\times 10^{-4}$	& [$9.67\times 10^{-4}-1.00\times 10^{-3}$] &	$3.63\times 10^{-3}$	& [$3.54\times 10^{-3}-3.66\times 10^{-3}$] \\
    &	PI-HPS, $r=0.5$	&	$3.92\times 10^{-5}$	& [$1.59\times 10^{-5}-8.53\times 10^{-5}$] &	$9.43\times 10^{-5}$	& [$3.75\times 10^{-5}-2.12\times 10^{-4}$] &	$4.28\times 10^{-4}$	& [$1.58\times 10^{-4}-9.45\times 10^{-4}$] \\
    &	PI-HLS, $r=1$	&	$5.76\times 10^{-4}$	& [$5.51\times 10^{-4}-6.08\times 10^{-4}$] &	$9.80\times 10^{-4}$	& [$9.68\times 10^{-4}-9.97\times 10^{-4}$] &	$3.59\times 10^{-3}$	& [$3.54\times 10^{-3}-3.66\times 10^{-3}$] \\
    &	PI-HCS, $r=1$	&	$5.74\times 10^{-4}$	& [$5.50\times 10^{-4}-6.06\times 10^{-4}$] &	$9.81\times 10^{-4}$	& [$9.70\times 10^{-4}-1.00\times 10^{-3}$] &	$3.60\times 10^{-3}$	& [$3.54\times 10^{-3}-3.66\times 10^{-3}$] \\[2pt]
    \noalign{\smallskip}\hline\noalign{\smallskip}
    \multirow{7}{*}{\rotatebox[origin=c]{90}{$h=20$}} 
    &	PSE	&	$1.64\times 10^{3}$	&				&	$5.57\times 10^{3}$	&				&	$1.49\times 10^{4}$	&				\\
    &	LSE	&	$1.13\times 10^{8}$	&				&	$4.11\times 10^{8}$	&				&	$1.18\times 10^{9}$	&				\\
    &	CSE	&	$4.33\times 10^{7}$	&				&	$1.58\times 10^{8}$	&				&	$4.51\times 10^{8}$	&				\\
    &	PI-HPS, $r=1$	&	$3.56\times 10^{-5}$	& [$7.83\times 10^{-6}-7.66\times 10^{-5}$] &	$8.68\times 10^{-5}$	& [$1.95\times 10^{-5}-1.76\times 10^{-4}$] &	$3.94\times 10^{-4}$	& [$8.67\times 10^{-5}-7.67\times 10^{-4}$] \\
    &	PI-HPS, $r=0.5$	&	$1.05\times 10^{-5}$	& [$4.55\times 10^{-6}-2.16\times 10^{-5}$] &	$1.94\times 10^{-5}$	& [$8.70\times 10^{-6}-4.29\times 10^{-5}$] &	$7.11\times 10^{-5}$	& [$3.00\times 10^{-5}-1.66\times 10^{-4}$] \\
    &	PI-HLS, $r=1$	&	$5.64\times 10^{-6}$	& [$4.75\times 10^{-6}-7.65\times 10^{-6}$] &	$7.28\times 10^{-6}$	& [$6.52\times 10^{-6}-9.15\times 10^{-6}$] &	$1.76\times 10^{-5}$	& [$1.73\times 10^{-5}-1.81\times 10^{-5}$] \\
    &	PI-HCS, $r=1$	&	$5.94\times 10^{-6}$	& [$4.17\times 10^{-6}-8.96\times 10^{-6}$] &	$7.70\times 10^{-6}$	& [$6.30\times 10^{-6}-1.10\times 10^{-5}$] &	$1.78\times 10^{-5}$	& [$1.73\times 10^{-5}-1.85\times 10^{-5}$] \\
    \noalign{\smallskip}\hline
    \end{tabular}
    }
    \label{tb:acc_lnEx}
\end{table}
As shown in \cref{tb:acc_lnEx}, the purely polynomial-based IM approximations (PSE, LSE, CSE), derived by the approach in \cite{kazantzis2001invariant,kazantzis2005model}, provide inaccurate results.~This inaccuracy worsens as $h$ increases, and is attributed to testing data in $\Omega=[-0.9,2]$ lying far from the equilibrium.~In this example, the neighborhood of the equilibrium is bounded by the radius of convergence $r=1$ of the true IM polynomial series expansions.~In contrast, the PI hybrid schemes and the PINN deliver accurate IM approximations; particularly, the PI hybrid scheme with power series and $r=0.5$ (PI-HPS, $r=0.5$) for low-degree polynomials and the PINN.~As $h$ increases, the accuracy of all PI hybrid schemes improves, and outperforms the approximation provided by PINN.

\begin{figure}[!h]
\centering
    \subfigure[Degree $h=10$ ]{\includegraphics[width=0.47\textwidth]{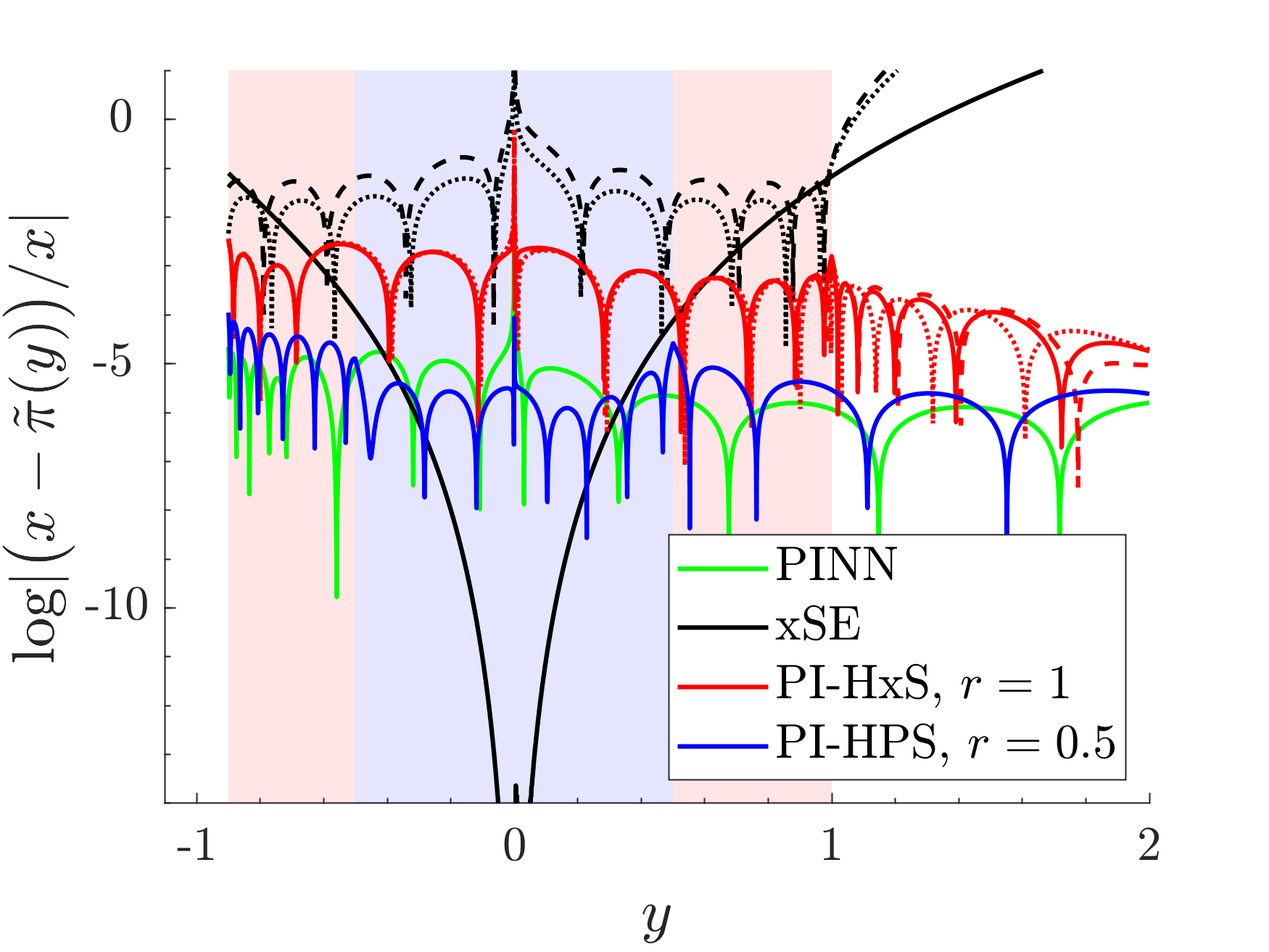}} \hspace{10pt}
    \subfigure[Degree $h=20$]{\includegraphics[width=0.47\textwidth]{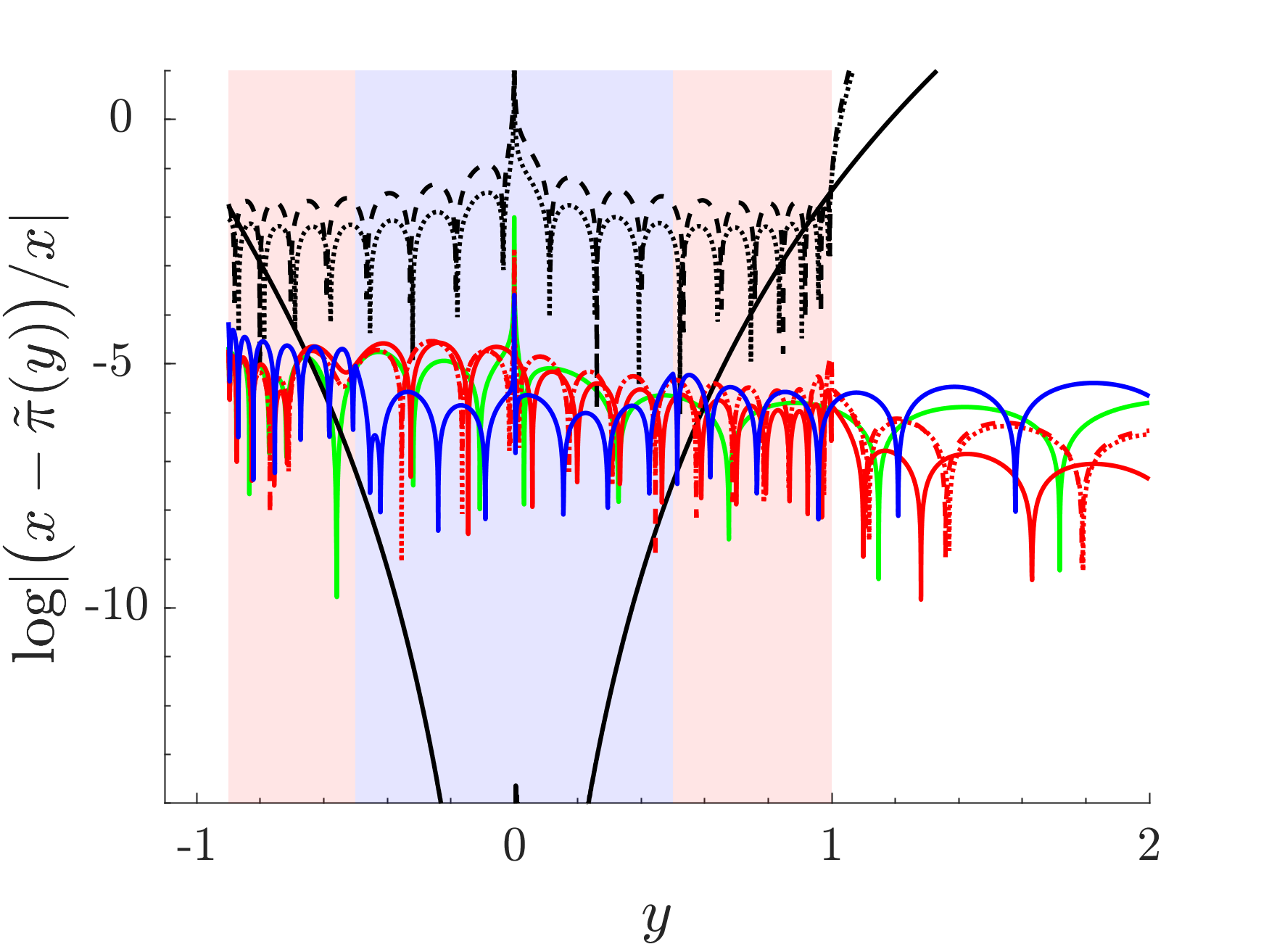}} 
\caption{Relative errors of IM approximations $\tilde{\pi}(y)$ compared to the true IM $x=\ln(1+y)$ for the example in \cref{app:lnEx} on the test data.~Panels show the polynomial series expansions xSE (black), the PI hybrid schemes PI-HxS (red, blue) and the PINN approximations (green).~Power, Legendre, and Chebyshev polynomial series ($x=P,L,C$) are distinguished by solid, dashed, and dotted curves, respectively.~Red and blue background indicates the hyperrectangle $\mathcal{D}$ for the PI hybrid schemes.~Panels (a) and (b) correspond to polynomial degrees $h=10$ and $h=20$, respectively, and $L=10$ neurons are used in the NNs.}
\label{fig:acc_lnEx}
\end{figure}

We highlight here that the error metrics in \cref{tb:acc_lnEx} are \textit{global} over the testing set $\Omega$.~To examine \textit{local} accuracy, we visualize the point-wise relative errors $\lvert x^{(s)}-\tilde{\pi}(y^{(s)}) \rvert/\lvert x^{(s)}\rvert$ in \cref{fig:acc_lnEx}, where panels (a,b) display the approximations in \cref{tb:acc_lnEx} for $h=10$ and $h=20$, respectively.
It is evident that the inaccuracy of all the polynomial series expansions arises from high errors far from the equilibrium.~\Cref{fig:acc_lnEx} confirms that the errors of the PI hybrid and PINN schemes are almost homogeneous in $\Omega$.~Interestingly, for $h=10$, the PI hybrid schemes exhibit lower accuracy than the PINN scheme, as reflected in \cref{tb:acc_lnEx}.~However, their accuracy improves and matches that of the PINN scheme, as $h$ increases.~It is also shown that the choice of polynomials does not affect the approximation accuracy of the hybrid schemes.~More importantly, with decreased $r$, the approximation provided by the PI hybrid schemes is more accurate near the equilibrium than that provided by the PINN.

\clearpage
\newpage

\begin{center}
\huge{Supplementary Material} \\[10pt]
\Large{Invariant Manifolds of Discrete-time Input-Driven Dynamical Systems via Hybrid Physics-Informed Neural Networks} \\[10pt]
\hrule
\large{Dimitrios G. Patsatzis, Nikolaos Kazantzis, Ioannis G. Kevrekidis, Lucia Russo, Constantinos Siettos}
\end{center}
\vspace{20pt}

\makeatletter
\renewcommand{\theHtable}{S1.\arabic{table}}
\renewcommand{\theHsection}{S1.\arabic{section}}
\makeatother

\newcounter{suppsection}
\renewcommand{\thesuppsection}{S\arabic{suppsection}}
\renewcommand{\thesection}{\thesuppsection}
\setcounter{suppsection}{1}
\renewcommand{\thetable}{\thesection.\arabic{table}}
\setcounter{table}{0}
\section{Convergence results of PI hybrid and standalone NN schemes}
\label{supp1}

\begin{table}[!h]
\centering
\caption{Convergence results of the proposed PI schemes for the enzymatic bioreactor problem in \cref{sb:Kax05ex} on the training set.~Mean values and 5–95\% percentiles of the loss function $\mathcal{L}(\cdot)$ and computational training times (in seconds) are reported for 100 randomly initialized realizations.~Schemes include standalone PINNs and PI hybrid schemes (PI-HxS with $x=P,L,C$ denoting power, Legendre, and Chebyshev polynomial series).~The radii of the PI-HPS $r$ vary from $0.5$ to $4$.~Hyperparameters: polynomial degree $h=10$ and $L=10$ neurons for the NN and PI-HxS schemes.}
\resizebox{\textwidth}{!}{
\begin{tabular}{l c c c c}
\hline\noalign{\smallskip}
\textbf{PI schemes} & \multicolumn{2}{c}{\textbf{loss $\mathcal{L}(\cdot)$}} & \multicolumn{2}{c}{\textbf{comp. time (in seconds)}} \\
 & mean & 5--95\% percentiles & mean & 5--95\% percentiles \\
\noalign{\smallskip}\hline\noalign{\smallskip}
PINN   & $3.63 \times 10^{-10}$ & [$9.71\times10^{-12}-7.46\times10^{-10}$] & $8.21\times10^{0}$  & [$3.72\times10^{0}-1.01\times10^{1}$] \\
PI-HPS, $r=0.5$      & $3.01\times10^{-10}$ & [$1.24\times10^{-11}-3.01\times10^{-9}$]  & $4.93\times10^{1}$  & [$6.30\times10^{0}-8.49\times10^{1}$] \\
PI-HPS, $r=1$        & $1.37\times10^{-9}$  & [$1.22\times10^{-11}-3.00\times10^{-9}$]  & $2.23\times10^{1}$  & [$1.20\times10^{1}-7.04\times10^{1}$] \\
PI-HPS, $r=2$        & $5.20\times10^{-10}$ & [$2.52\times10^{-11}-1.97\times10^{-9}$]  & $7.17\times10^{1}$  & [$1.67\times10^{1}-9.44\times10^{1}$] \\
PI-HPS, $r=4$        & $3.60\times10^{-11}$ & [$3.60\times10^{-11}-3.60\times10^{-11}$] & $3.86\times10^{-2}$ & [$2.17\times10^{-2}-6.15\times10^{-2}$] \\
PI-HLS, $r=1$        & $2.31\times10^{-10}$ & [$1.10\times10^{-11}-4.24\times10^{-9}$]  & $8.16\times10^{1}$  & [$4.97\times10^{1}-9.45\times10^{1}$] \\
PI-HCS, $r=1$        & $2.84\times10^{-10}$ & [$8.35\times10^{-12}-4.70\times10^{-9}$]  & $6.60\times10^{1}$  & [$1.95\times10^{1}-9.36\times10^{1}$] \\
\noalign{\smallskip}\hline
\end{tabular}
}
\label{tb:conv_Kaz05ex}
\end{table}

\begin{table}[!h]
\centering
\caption{Convergence results of the proposed PI schemes for the car-following problem with an autonomou leader in \cref{sb:CFM} on the training set.~Mean values and 5–95\% percentiles of the loss function $\mathcal{L}(\cdot)$ and computational training times (in seconds) are reported for 100 randomly initialized realizations.~Schemes include standalone PINNs and PI hybrid schemes (PI-HxS with $x=P,L,C$ denoting power, Legendre, and Chebyshev polynomial series).~Hyperparameters: polynomial degree $h=3$, radius $r=1$ for the PI-HxS schemes and $L=20$ neurons for the NN and PI-HxS schemes.}
\resizebox{\textwidth}{!}{
\begin{tabular}{l c c c c}
\hline\noalign{\smallskip}
\textbf{PI schemes} & \multicolumn{2}{c}{\textbf{loss $\mathcal{L}(\cdot)$}} & \multicolumn{2}{c}{\textbf{comp. time (in seconds)}} \\
 & mean & 5--95\% percentiles & mean & 5--95\% percentiles \\
\noalign{\smallskip}\hline\noalign{\smallskip}
PINN 	        &	$7.50\times 10^{-5}$	& [$4.30\times 10^{-5}-1.18\times 10^{-4}$]&	$5.99\times 10^{3}$	& [$5.74\times 10^{3}-6.38\times 10^{3}$] \\
PI-HPS, $r=1$	&	$9.94\times 10^{-5}$	& [$6.76\times 10^{-5}-1.37\times 10^{-4}$]&	$8.16\times 10^{3}$	& [$7.80\times 10^{3}-8.81\times 10^{3}$] \\
PI-HLS, $r=1$	&	$1.08\times 10^{-4}$	& [$6.97\times 10^{-5}-1.37\times 10^{-4}$]&	$8.31\times 10^{3}$	& [$7.82\times 10^{3}-8.99\times 10^{3}$] \\
PI-HCS, $r=1$  	&	$1.08\times 10^{-4}$	& [$7.58\times 10^{-5}-1.38\times 10^{-4}$]&	$8.02\times 10^{3}$	& [$7.85\times 10^{3}-8.56\times 10^{3}$] \\
\noalign{\smallskip}\hline
\end{tabular}
}
\label{tb:conv_CFM}
\end{table}

\begin{table}[!h]
    \centering
    \caption{Convergence results of the proposed PI schemes for the example in \cref{app:lnEx} on the training set.~Mean values and 5–95\% percentiles of the loss function $\mathcal{L}(\cdot)$ and computational training times (in seconds) are reported for 100 randomly initialized realizations.~Schemes include standalone PINN and PI hybrid schemes (PI-HxS with $x=P,L,C$ denoting power, Legendre, and Chebyshev polynomial series).~Hyperparameters: polynomial degrees $h=10$ and $h=20$, radius $r=1$ and $r=0.5$ for the HPS schemes and $L=10$ neurons for the NN and HxS schemes.}
    \resizebox{\textwidth}{!}{
    \begin{tabular}{l l c c c c }
    \hline\noalign{\smallskip}
    &	\textbf{PI schemes} &	\multicolumn{2}{c}{\textbf{loss $\mathcal{L}(\cdot)$}}&	\multicolumn{2}{c}{\textbf{comp. time (in seconds)}}		\\
    &		&	mean	&	5-95\% percentiles			&	mean	&	5-95\% percentiles			\\
    \noalign{\smallskip}\hline\noalign{\smallskip}
    &	PINN	       &	$1.95\times10^{-7}$	& [$1.44\times10^{-8}-1.28\times10^{-6}$] &	$6.06\times10^{0}$	& [$3.90\times10^{0}-8.40\times10^{0}$] \\
    \noalign{\smallskip}\hline\noalign{\smallskip}
    \multirow{4}{*}{\rotatebox[origin=c]{90}{$h=10$}} 
    &	PI HPS, $r=0.5$ &	$6.74\times10^{-7}$	& [$1.04\times10^{-7}-3.31\times10^{-6}$] &	$4.47\times10^{1}$	& [$1.87\times10^{1}-4.96\times10^{1}$] \\
    &	PI HPS, $r=1$	&	$1.66\times10^{-4}$	& [$3.37\times10^{-6}-9.31\times10^{-4}$] &	$4.42\times10^{1}$	& [$1.82\times10^{1}-4.93\times10^{1}$] \\
    &	PI HLS, $r=1$	&	$2.00\times10^{-4}$	& [$4.69\times10^{-6}-8.96\times10^{-4}$] &	$4.27\times10^{1}$	& [$4.09\times10^{1}-4.56\times10^{1}$] \\
    &	PI HCS, $r=1$	&	$1.52\times10^{-4}$	& [$2.87\times10^{-6}-7.46\times10^{-4}$] &	$4.12\times10^{1}$	& [$2.08\times10^{1}-4.56\times10^{1}$] \\
    \noalign{\smallskip}\hline\noalign{\smallskip}
    \multirow{4}{*}{\rotatebox[origin=c]{90}{$h=20$}}	
    &	PI HPS, $r=0.5$&	$5.65\times10^{-7}$	& [$2.91\times10^{-8}-2.30\times10^{-6}$] &	    $5.01\times10^{1}$	& [$4.17\times10^{1}-5.45\times10^{1}$] \\
    &	PI HPS, $r=1$	&	$3.04\times10^{-8}$	& [$3.59\times10^{-9}-1.60\times10^{-7}$] &	    $5.51\times10^{1}$	& [$3.47\times10^{1}-6.18\times10^{1}$] \\
    &	PI HLS, $r=1$	&	$1.00\times10^{-8}$	& [$8.32\times10^{-10}-3.94\times10^{-8}$] &	$1.18\times10^{1}$	& [$3.82\times10^{0}-2.26\times10^{1}$] \\
    &	PI HCS, $r=1$	&	$7.06\times10^{-9}$	& [$9.19\times10^{-10}-2.63\times10^{-8}$] &	$1.44\times10^{1}$	& [$4.16\times10^{0}-2.45\times10^{1}$] \\
    \noalign{\smallskip}\hline
    \end{tabular}
    }
    \label{tb:conv_lnEx}
\end{table}

\clearpage
\newpage
\makeatletter
\renewcommand{\theHfigure}{S2.\arabic{figure}}
\renewcommand{\theHsection}{S2.\arabic{section}}
\makeatother

\setcounter{suppsection}{2}
\renewcommand{\thefigure}{\thesection.\arabic{figure}}
\setcounter{figure}{0}
\section{Relative errors of IM approximation for the car-following problem with an autonomous leader}
\label{supp2}

\begin{figure}[htp]
\centering
    \begin{tikzpicture}
        \node[anchor=north west,inner sep=0pt] at (0,0){\includegraphics[width=0.45\textwidth]{CFM/Err_d3_PSE_x1.png}};
        \node[font=\sffamily\bfseries\large] at (0.5ex,-2ex) {a};
    \end{tikzpicture}
    \hspace{20pt}
    \begin{tikzpicture}
        \node[anchor=north west,inner sep=0pt] at (0,0){\includegraphics[width=0.45\textwidth]{CFM/Err_PINN_x1.png}};
        \node[font=\sffamily\bfseries\large] at (0.5ex,-2ex) {b};
    \end{tikzpicture}
    \\
    \begin{tikzpicture}
        \node[anchor=north west,inner sep=0pt] at (0,0){\includegraphics[width=0.45\textwidth]{CFM/Err_d3_PIHPS_x1.png}};
        \node[font=\sffamily\bfseries\large] at (0.5ex,-2ex) {c};
    \end{tikzpicture}
    \hspace{20pt}
    \begin{tikzpicture}
        \node[anchor=north west,inner sep=0pt] at (0,0){\includegraphics[width=0.45\textwidth]{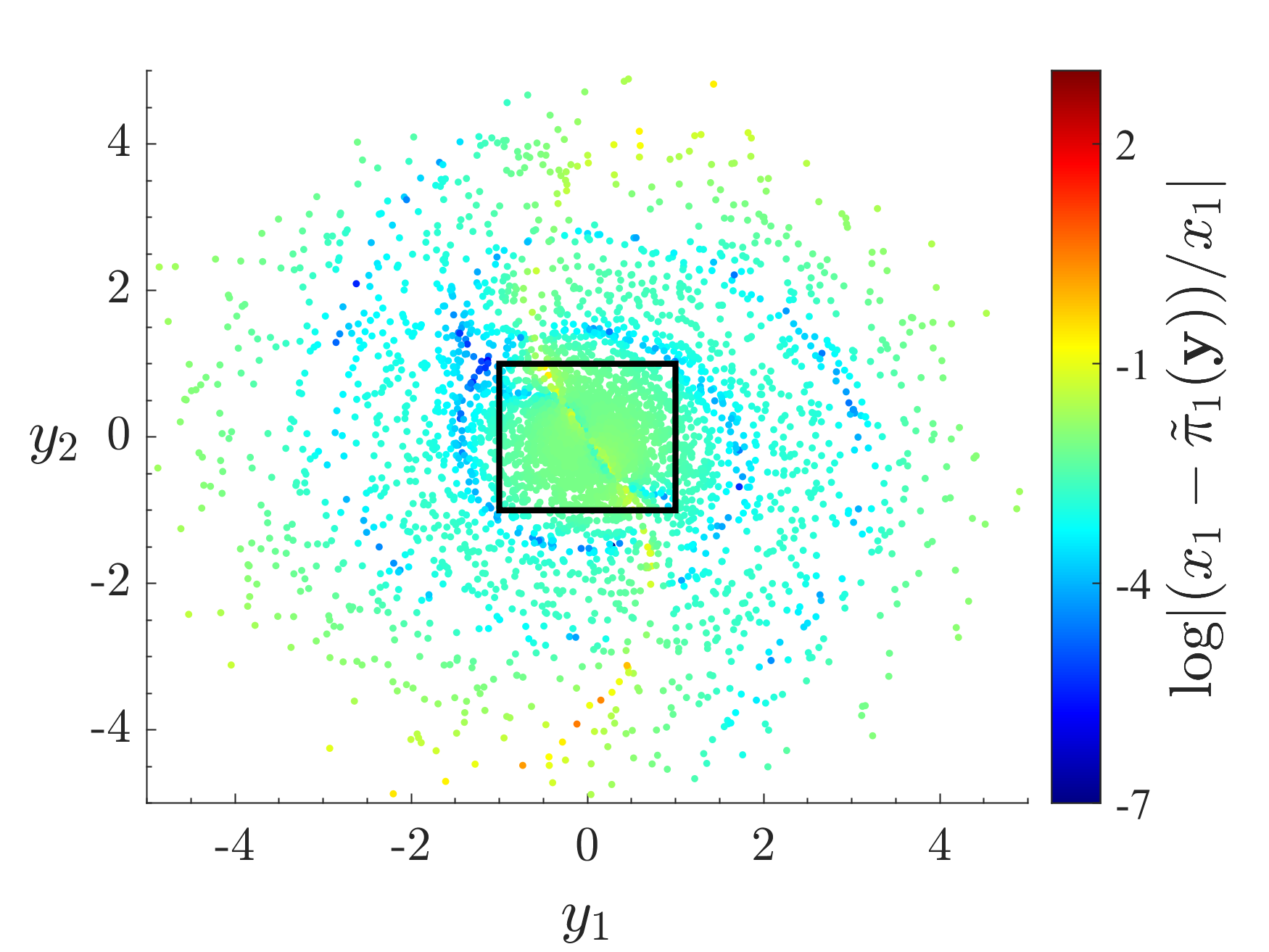}};
        \node[font=\sffamily\bfseries\large] at (0.5ex,-2ex) {d};
    \end{tikzpicture}
    \\
    \hspace{0.45\textwidth} \hspace{28pt}
    \begin{tikzpicture}
        \node[anchor=north west,inner sep=0pt] at (0,0){\includegraphics[width=0.45\textwidth]{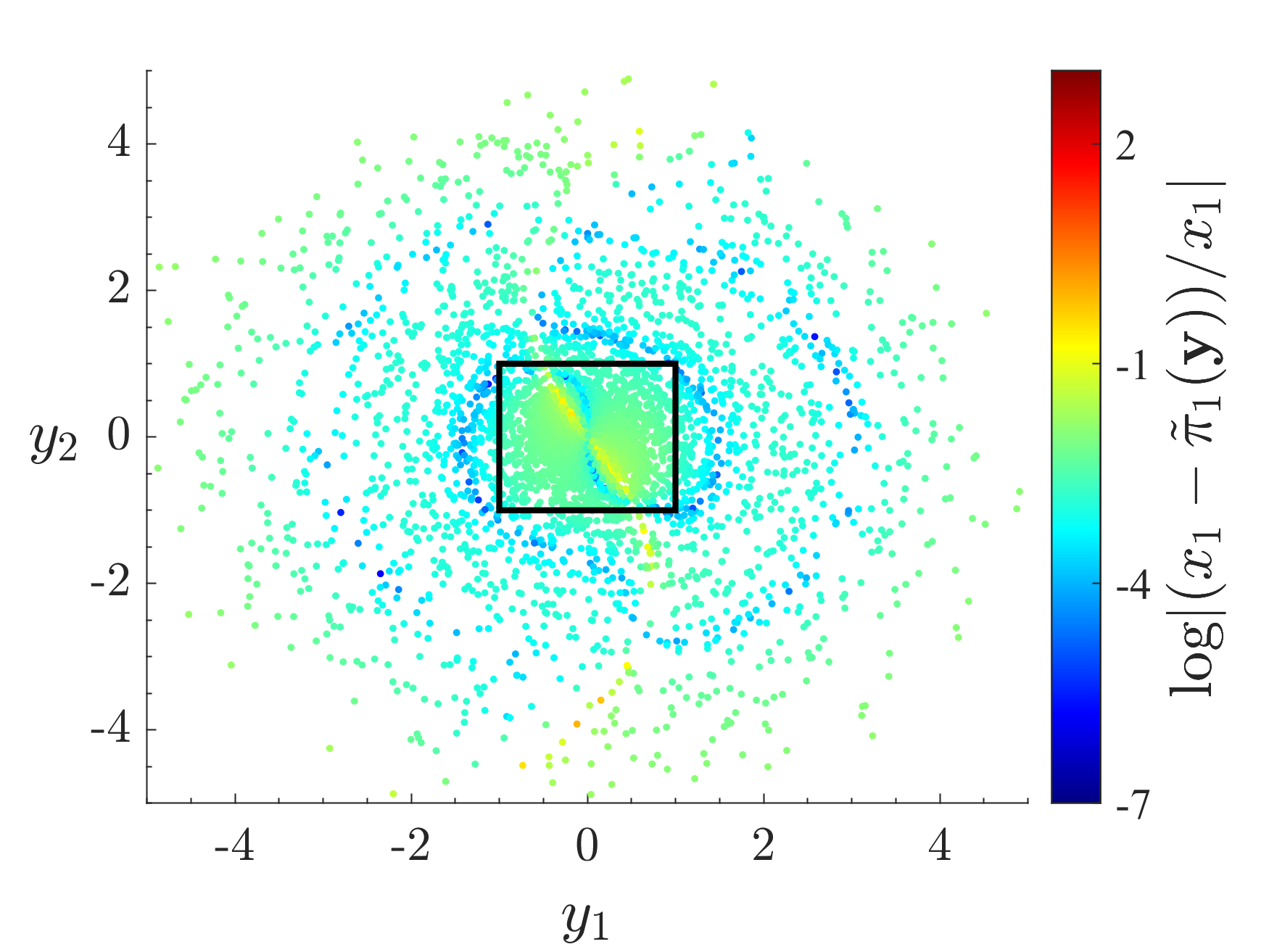}};
        \node[font=\sffamily\bfseries\large] at (0.5ex,-2ex) {e};
    \end{tikzpicture} 
\caption{Relative errors of the $x_1$ component of the IM approximations $\tilde{\boldsymbol{\pi}}(\mathbf{y})$ for the car-following problem with an autonomous leader in \Cref{sb:CFM}, on the test data lying on the IM.~Panels show (a) the PSE, (b) the PINN, and (c,d,e) the PI hybrid schemes with Power, Legendre, and Chebyshev polynomial series.~The black rectangle in panels (c-e) denotes the rectangle $\mathcal{D}$ of the PI hybrid schemes.~Polynomial degree is $h=3$ for all cases and $L=20$ neurons are used in the NNs.}
\label{fig:acc_CFM_x1_full}
\end{figure}

\begin{figure}[htp]
\centering
    \begin{tikzpicture}
        \node[anchor=north west,inner sep=0pt] at (0,0){\includegraphics[width=0.45\textwidth]{CFM/Err_d3_PSE_x5.png}};
        \node[font=\sffamily\bfseries\large] at (0.5ex,-2ex) {a};
    \end{tikzpicture}
    \hspace{20pt}
    \begin{tikzpicture}
        \node[anchor=north west,inner sep=0pt] at (0,0){\includegraphics[width=0.45\textwidth]{CFM/Err_PINN_x5.png}};
        \node[font=\sffamily\bfseries\large] at (0.5ex,-2ex) {b};
    \end{tikzpicture}
    \\
    \begin{tikzpicture}
        \node[anchor=north west,inner sep=0pt] at (0,0){\includegraphics[width=0.45\textwidth]{CFM/Err_d3_PIHPS_x5.png}};
        \node[font=\sffamily\bfseries\large] at (0.5ex,-2ex) {c};
    \end{tikzpicture}
    \hspace{20pt}
    \begin{tikzpicture}
        \node[anchor=north west,inner sep=0pt] at (0,0){\includegraphics[width=0.45\textwidth]{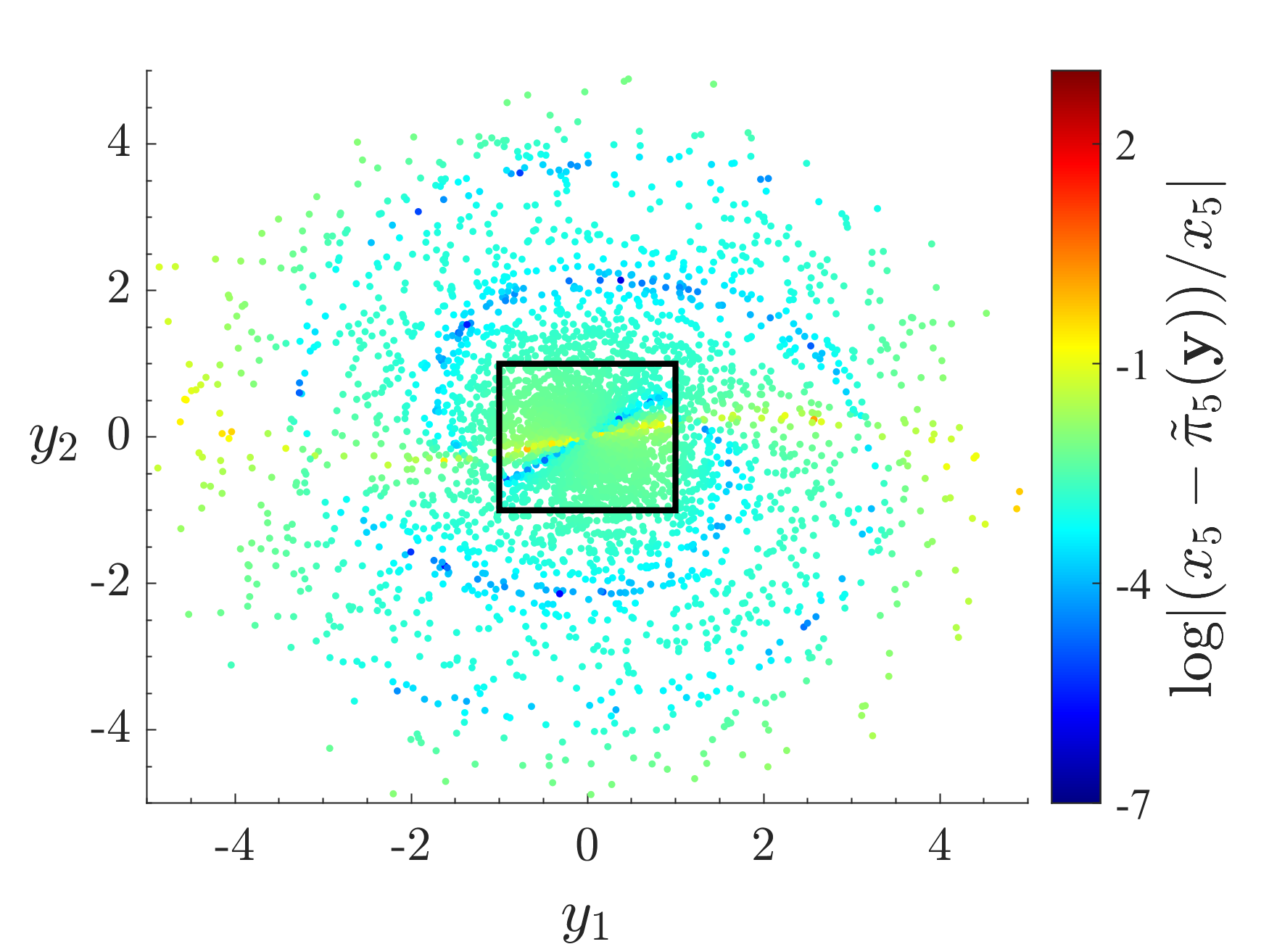}};
        \node[font=\sffamily\bfseries\large] at (0.5ex,-2ex) {d};
    \end{tikzpicture}
    \\
    \hspace{0.45\textwidth} \hspace{28pt}
    \begin{tikzpicture}
        \node[anchor=north west,inner sep=0pt] at (0,0){\includegraphics[width=0.45\textwidth]{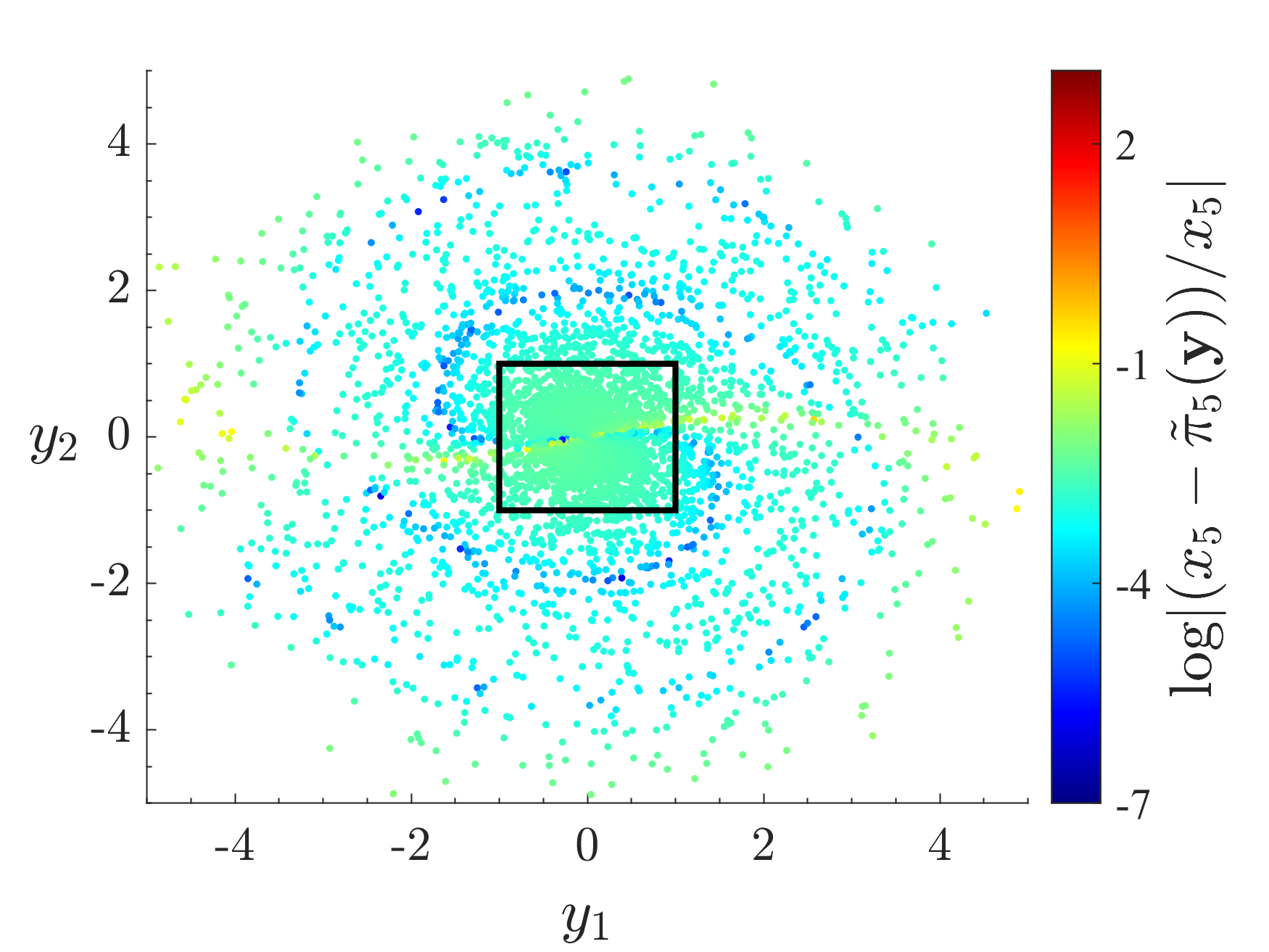}};
        \node[font=\sffamily\bfseries\large] at (0.5ex,-2ex) {e};
    \end{tikzpicture} 
\caption{Relative errors of the $x_5$ component of the IM approximations $\tilde{\boldsymbol{\pi}}(\mathbf{y})$ for the car-following problem with an autonomous leader in \Cref{sb:CFM}, on the test data lying on the IM.~Panels show (a) the PSE, (b) the PINN, and (c,d,e) the PI hybrid schemes with Power, Legendre, and Chebyshev polynomial series.~The black rectangle in panels (c-e) denotes the rectangle $\mathcal{D}$ of the PI hybrid schemes.~Polynomial degree is $h=3$ for all cases and $L=20$ neurons are used in the NNs.}
\label{fig:acc_CFM_x5_full}
\end{figure}

\end{document}